\documentclass[final]{siamltex}

\usepackage{showkeys}
\usepackage{amsmath}
\usepackage{amsfonts}
\usepackage{amssymb}
\usepackage{a4wide}
\usepackage{graphicx}
\usepackage{array}
\usepackage[colorlinks=true,pdftex]{hyperref}
\usepackage{morefloats}

\newcommand{\del}{\partial}

\title{Numerical study of blow-up in solutions to generalized Kadomtsev-Petviashvili equations}
\author{C.~Klein\thanks{Institut de Math\'ematiques de Bourgogne,
		Universit\'e de Bourgogne, 9 avenue Alain Savary, 21078 Dijon
		Cedex, France
    ({\tt christian.klein@u-bourgogne.fr})}
\and
R.~Peter\thanks{Institut de Math\'ematiques de Bourgogne,
		Universit\'e de Bourgogne, 9 avenue Alain Savary, 21078 Dijon
		Cedex, France
    ({\tt ralf.peter@u-bourgogne.fr})}
}
\date{\today}

\begin{document}
\maketitle
\begin{abstract}
We present a numerical study of solutions to the
generalized Kadomtsev-Petviashvili  equations with critical and 
supercritical nonlinearity for localized initial data with a 
single minimum and single maximum. In the cases with blow-up, we use a 
dynamic rescaling to identify the type of the singularity. We present 
a discussion of the observed blow-up scenarios.

\end{abstract}	  
\section{Introduction}
The celebrated Korteweg-de Vries (KdV) equation modelling 
wave phenomena in $1+1$ dimensions in the limit of long wavelengths 
appears in many domains of application as hydrodynamics, acoustics, 
nonlinear optics, plasma physics, \ldots. It 
has a $2+1$ dimensional generalization for essentially 
one-dimensional wave phenomena with weak transverse effects, the 
Kadomtsev-Petviashvili (KP) equations \cite{KP} with a similar range 
of applications. Remarkably both KdV and KP are completely 
integrable. Therefore many explicit solutions as solitons are known 
for these equations. The Cauchy problem for the KP equations and the stability 
of certain exact solutions are actively studied, see for instance 
\cite{KS12} for a recent review. 

Since both KdV and KP appear in approximations for long wavelengths, 
the dispersion in these equations is too strong compared to what is 
observed in applications. A possible cure to this short coming is to 
shift the balance between nonlinearity and dispersion towards the 
nonlinearity. This leads for KdV to generalized KdV (gKdV) equations, 
and for KP to generalized KP (gKP) equations
\begin{equation}
 u_{t} + u^n u_{x} + u_{xxx}
            + \lambda\del_x^{-1} u_{yy} = 0   
    \label{gKP}
\end{equation}
where $n\in\mathbb{Q}$, $n\geq 1$ and where $\lambda=\pm1$; the antiderivative 
$\del_x^{-1}$ is defined via its Fourier symbol $-i/k_{x}$. Note that 
gKP for $n=2$ appears  as a model for the evolution of 
sound waves in antiferromagnetic materials, see \cite{FT85}. We will 
always consider asymptotically (for $|x|,|y|\to\infty$) decreasing 
solutions in the following. The case $\lambda=1$ is 
denoted gKP II, the one with $\lambda=-1$ gKP I. gKP I has a 
\textit{focusing} effect, gKP II a \textit{defocusing} one. If $u$ is 
independent of $y$, the gKP equations reduce to the gKdV equations; 
for $n=1$, one obtains the standard KP equations. The equations are 
only completely integrable in the latter case, but have conserved 
quantities as the \textit{mass}, the $L_{2}$ norm of $u$, 
\begin{equation}
    M[u] = ||u||_{2},
    \label{mass}
\end{equation}
and the 
\textit{energy}
\begin{equation}
            E[u] = \int_{\mathbb{R}^{2}}^{}\left(\frac{1}{2}u_x^2
           - \frac{1}{(n+1)(n+2)}
	   u^{n+2}-\frac{\lambda}{2}(\partial_{x}^{-1}u_{y})^{2}\right)dx\,dy
    \label{energy}.
\end{equation}
Note that the 
equations (\ref{gKP}) are invariant under rescalings of the form $x\to x/s$, 
$y\to y/s^{2}$, $t\to t/s^{3}$ and $u\to s^{2/n}u$ with $s=const$. 
Since the $L_{2}$ 
norm of $u$ is invariant under this rescaling for $n=4/3$, this case is also 
referred to as $L_{2}$ critical. For gKdV the critical 
exponent in this sense is $n=4$.

The appearance of an antiderivative in (\ref{gKP}) has the consequence that 
\begin{equation}
    \int_{\mathbb{R}}^{}\partial_{yy}u(x,y,t)\,dx=0,\quad \forall t>0
    \label{const},
\end{equation}
which can be easily seen by differentiating (\ref{gKP}) with respect 
to $x$ and integrating over the real line. This is 
even true if this constraint is not satisfied for the initial data 
$u_{0}(x,y)$. In \cite{FS,MST} it was shown that the solution to a 
Cauchy problem not satisfying the constraint will not be smooth in 
time for $t=0$. Numerical experiments in \cite{KSM} indicate that 
the solution develops an 
infinite `trench' the integral over which just ensures that 
(\ref{const}) is fulfilled. This non-regularity of the solution for 
$t\to0$ is numerically a problem in the sense that it delimits the 
achievable accuracy. To avoid such problems, we always consider here initial 
data that are $x$-derivatives of a rapidly decreasing function and 
thus satisfy (\ref{const}); as a concrete example we will study in 
this paper the initial data
\begin{equation}
    u_{0}(x,y)=\beta\partial_{xx}\exp(-(x^{2}+y^{2}))
    \label{initial}
\end{equation}
with $\beta>1$.

Another consequence of the antiderivative in (\ref{gKP}) is that  
solutions to Cauchy problems with 
initial data $u_{0}(x,y)$ in the Schwarz space $\mathcal{S}(\mathbb{R}^{2})$ of 
rapidly decreasing functions will not stay in 
this space unless $u_{0}(x,y)$ satisfies an 
infinite number of constraints. This can be already seen on 
the level of the linearized KP equation, see e.g.~\cite{BPP,KSM}, where 
the Green's function implies a slow algebraic decrease in $y$ for 
$|y| \to\infty$. This leads to the formation of \emph{tails} with an 
algebraic decrease towards infinity for generic localized initial data. 
The amplitude of these  grows with time (see for instance \cite{KSM}).
 
The gKP I equations are known to have localized travelling (in 
$x$-direction) wave 
solutions called \textit{lumps} for $n<4$. This corresponds to solutions to gKP of 
the form $u(x,y,t)=Q(x-ct,y)$, $c=const$, where (the nontrivial) $Q(z,y)$ 
satisfies
\begin{equation}
    -cQ_{zz}+\frac{1}{n+1}(Q^{n+1})_{zz}+Q_{zzzz}+\lambda Q_{yy}=0
    \label{soliton}.
\end{equation}
The only explicit form of these 
solutions is known for KP I in terms of rational functions, see 
\cite{MZBM}.  For general nonlinearities, solutions  of (\ref{soliton}) were studied 
in \cite{BS97a,BS97b}. It was shown that they
have an algebraic fall-off in 
$x$ and $y$ for $|x|,|y| \to \infty$, i.e.,   
 that the solutions cannot decrease more rapidly 
than $1/(x^{2}+y^{2})$. The precise fall off rate was given in \cite{Gra08}. 
The lumps are unstable for $n> 4/3$, see  \cite{BS97a,BS97b}. There 
are no such solutions for gKP II.
For KP I it was argued in \cite{AF} for small initial 
data that the solutions for large $t$ decompose into lumps and 
radiation, a fact stressing the importance of lumps. It is conjectured 
that this behavior is also true for general initial data. Due to the 
instability of lumps for gKP I equations with $n\geq 4/3$, 
an evolution of initial 
data into an array of lumps cannot be expected for this case for 
large $t$. But it 
is interesting to study which structures can be observed in the time 
evolution of localized initial data and which blow up. This is one of the goals of 
this paper. For gKdV, it was found numerically in \cite{KP13} that 
initial data with sufficient mass can lead to the formation of 
several solitons, and that the one appearing first will finally blow 
up.

It is well known that there can be blow-up, i.e., a loss of 
regularity with respect to the initial data,  in 
solutions to gKdV and gKP I equations with $n\geq 4$ and $n\geq 4/3$ respectively
for certain initial data, see \cite{Liu2001}. For the latter it was 
shown in \cite{MST} that the $L_{2}$ norm of $u_{y}$ blows up in 
these cases. It is unclear whether there can be 
blow-up in solutions to gKP II equations. A first numerical study of these issues 
has been presented in \cite{KS12}\footnote{Note that in that paper, 
the gKP equations, which were numerically studied, had a factor $1/2$ 
in front of the nonlinearity which was not mentioned in the text. 
Therefore some of the results there differ from what is presented 
here.}. In this paper we present a more 
detailed numerical analysis of these questions with the goal to 
identify the type of blow-up for localized initial data with a single 
minimum. 

The paper is organized as follows: in Section 2 we present the used 
numerical methods and discuss a dynamic rescaling of the gKP 
equations. In section 3 we consider the $L_{2}$ critical case $n=4/3$ 
and discuss various examples. In the cases with blow-up, we try to 
identify the type of 
the singularity. The same analysis is performed in section 4 for the 
case $n=2$ as an example for a supercritical situation. In section 5 
we consider blow-up in gKP II solutions for $n=3$ and $n=4$. 
We add some concluding remarks in section 6.

\section{Numerical methods}
In this section we present the numerical methods to be used in this 
paper to integrate the gKP equations. The main tool will be a direct 
integration of the equations with a Fourier spectral method for the 
spatial coordinates, and a fourth order \textit{exponential time differencing} 
(ETD)
scheme for the time 
integration. We also discuss a dynamic rescaling of the equation in 
order to study blow-up cases in more detail.

\subsection{Direct numerical integration}
Fourier spectral approaches are generally the most efficient method for the 
numerical treatment of smooth periodic 
functions. The reason for this is the in practice exponential 
decrease of the Fourier coefficients of such functions which implies excellent 
approximation properties for the latter. Due to the finite 
numerical precision, rapidly decreasing functions can be seen as 
essentially periodic if the computational domain is chosen large 
enough such that the function and its first derivatives vanish 
at the boundaries of the domain. For such functions the Fourier 
coefficients decrease to machine precision ($10^{-16}$ in our case)
for a moderate number $N$ of 
Fourier modes. The two-dimensional Fourier transform is defined as 
\begin{equation}
    \hat{u}(k_{x},k_{y})=\int_{\mathbb{R}^{2}}^{}u(x,y)e^{-i(k_{x}x+k_{y}y)}dx\, dy.
    \label{fourier}
\end{equation}
For gKP there is the additional advantage that the 
antiderivative is defined most conveniently in Fourier space. The 
Fourier multiplier $-i/k_{x}$ is regularized for $k_{x}\to0$ as 
discussed in \cite{KR11} by adding a small imaginary constant to 
$k_{x}$. The gKP equation takes in Fourier 
space the form 
\begin{equation}
    \hat{u}_{t}=\mathcal{F}(u)+\mathcal{L}\hat{u},
    \label{gKdVfourier}
\end{equation}
where $\mathcal{L}=ik_{x}^{3}-i\lambda k_{y}^{2}/k_{x}$, and where
$\mathcal{F}(u)=-ik_{x}\widehat{u^{n+1}}/(n+1)$. To avoid numerical 
errors due to the regularization of the term  $-i/k_{x}$, we rewrite 
(\ref{gKdVfourier}) as an equation for $\hat{w}$, where 
$\hat{u}=ik_{x}\hat{w}$, i.e., for the integral of $u$ with respect 
to $x$. This is possible since we use an explicit scheme as explained 
below, and since the initial condition $u_{0}$ is chosen to be an 
$x$-derivative of a rapidly decreasing function. In addition the term 
$\mathcal{F}$ in (\ref{gKdVfourier}) is an $x$-derivative itself, 
whereas $\mathcal{L}$ appears in the used scheme only in the exponential 
integrator. Thus we use (\ref{etd}) for $\hat{w}$ and 
$\mathcal{F}/(ik_{x})$ in the actual computation which ensures that 
$u$ is for all time steps an $x$-derivative within numerical precision.

An important practical advantage of Fourier transforms is the fact 
that the derivatives are diagonal in Fourier space. This allows for 
an efficient time integration since for equations of 
the form (\ref{gKdVfourier}), there are many high-order time 
integrators, see 
e.g.~\cite{CoxMatthews2002,KassamTrefethen2005,HO,Klein2008,KleinRoidot2011} and references therein, especially for 
diagonal $\mathcal{L}$ as in the Fourier case. For the numerical 
integration of (\ref{gKdVfourier}), the Fourier transform will be 
approximated via a truncated Fourier series, a \emph{discrete Fourier 
transform} which will be computed with a \emph{fast Fourier 
transform}. Thus the equations  (\ref{gKP}) 
will be approximated via a system of ordinary differential equations 
(ODEs) of finite dimension for the Fourier coefficients. The latter 
system is in the present case 
\emph{stiff}, where the term stiffness is used to indicate that there 
are vastly different timescales in the studied problem. This makes the use of 
standard explicit methods inefficient for stability reasons.

It was shown in \cite{KleinRoidot2011} that 
ETD schemes perform best for KP equations
among the studied stiff integrators, and  that the performance 
is similar for different ETD 
schemes. For these methods one uses a constant time 
step $h=t_{m+1}-t_{m}$ and integrates (\ref{gKdVfourier}) with an 
exponential factor to obtain
\begin{equation}
    \hat{u}(t_{m+1})=e^{\mathcal{L}h}\hat{u}(t_{m})+\int_{0}^{h} 
e^{\mathcal{L}(h-\theta)}\mathcal{F}(\hat{u}(\theta+t_{m}),\theta+t_{m})\,d\theta.
    \label{etd}
\end{equation}
The different ETD schemes differ in the approximation of the integral 
in (\ref{etd}).
We 
use here the method by Cox and Matthews \cite{CoxMatthews2002} of classical order 
four. An important aspect in the implementation of ETD schemes is the accurate 
and efficient computation of the so-called $\phi$-functions,
$$\phi_{i}=\frac{1}{(i-1)!}\int_{0}^{1} 
e^{(1-\tau)z}\theta^{i-1}\,d\theta,\quad i=1,2,\ldots,$$
which appear in all ETD schemes. To avoid cancellation errors, we use 
here
contour integrals in the complex plane as in \cite{KassamTrefethen2005} in the enhanced version 
 \cite{Schme} as discussed in \cite{Klein2008}.

Accuracy of the numerical solution is controlled as in 
\cite{Klein2008,KleinRoidot2011} via the numerically computed mass 
(\ref{mass}) or energy (\ref{energy})
which will depend on time due to unavoidable numerical errors. We use 
the quantity
\begin{equation}
    \Delta=|E(t)/E(0)-1|
    \label{Delta}
\end{equation}
(and similar for $M$) as an indicator of the numerical 
precision. It was shown in \cite{Klein2008,KleinRoidot2011} that the numerical accuracy 
of this quantity overestimates the $L_{\infty}$ norm of the 
difference between numerical and exact solution by two to three 
orders of magnitude. A precondition for the usability of this 
quantity is  a large 
enough number of Fourier modes. In practice the numerical error 
cannot be smaller than the smallest modulus of the Fourier 
coefficients. Note that the energy is more sensitive to a loss of 
regularity since it contains a derivative. Thus we generally consider 
this quantity unless otherwise noted. But the mass conservation is 
typically of the same order as the energy conservation in the 
numerical experiments.

We typically choose the number $N$ of Fourier modes high enough 
that they decrease to machine precision for the 
initial data. The number $N$ thus depends   
on the size of the computational domain, $x\in[-\pi,\pi]L_{x}$, 
$y\in[-\pi,\pi]L_{y}$ where the 
real constants $L_{x}$, $L_{y}$ are chosen large enough to ensure `periodicity' of 
the initial data in the sense discussed above. For gKP there is the 
problem of the algebraic fall off of the solution for $t>0$ even if 
the initial data are rapidly decreasing. Thus we are forced to use a 
larger domain than would be necessary for corresponding initial data 
for gKdV where the solution stays rapidly decreasing if the initial 
data are. In practice we reach a resolution in Fourier space of 
better than $10^{-10}$ for times much smaller than the time of a 
possible blow-up, see also \cite{KleinRoidot2011}. 

The occurrence of a blow-up leads to an increase of 
the Fourier modes for the high wave numbers which eventually 
breaks the code. Since we use an explicit method, we can 
afford sufficiently small time steps to have satisfactory resolution in 
time up to $t\sim t^{*}$. Typically we run out of resolution in 
Fourier space first.  A blow-up is also identified via diverging 
norms of the solution (see below). 

\subsection{Dynamic rescaling}
In \cite{KP13} we have used a dynamic rescaling of the gKdV equation 
to analyze blow-up in more detail with an adaptive approach. This 
method can be also generalized to gKP equations. The basic idea is to 
use the scaling invariance of the gKP equation discussed in the 
previous section, but now with a time dependent scaling factor. 
As for gKdV we consider the 
coordinate change with a factor $L(t)$
\begin{equation}
    \xi = \frac{x-x_{m}}{L},\quad \eta=\frac{y-y_{m}}{L^{2}},\quad 
    \frac{d\tau}{dt}=\frac{1}{L^{3}},\quad U = L^{2/n}u
    \label{gKP4}.
\end{equation}
This leads for (\ref{gKP}) to
\begin{equation}
    U_{\tau}-a\left(\frac{2}{n}U+\xi U_{\xi}+2\eta 
    U_{\eta}\right)-v_{\xi}U_{\xi}-v_{\eta}U_{\eta}
    +U^{n}U_{\xi}+U_{\xi\xi\xi}+\lambda 
    \int_{-\infty}^{\xi}U_{\eta\eta}\,d\xi=0
    \label{gKP5},
\end{equation}
with 
\begin{equation}
    a=(\ln L)_{\tau}, \quad v_{\xi} = \frac{x_{m,\tau}}{L},\quad 
    v_{\eta}=\frac{y_{m,\tau}}{L^{2}},
    \label{a}
\end{equation}
where the index $\tau$ denotes the derivative with respect to $\tau$.
Thus the space and time scales are changed 
adaptively around blow-up which is reached here for $\tau\to\infty$. 
The asymmetry in $x$ and $y$ with respect to the rescaling with $L$ 
also explains the stronger divergence of the $y$-derivative of $u$ at 
blow-up than of the $x$-derivative as observed in \cite{KS12}.

As for gKdV, equation (\ref{gKP5}) is also important for theoretical 
purposes in describing asymptotically blow-up. In the limit 
$\tau\to\infty$,  the functions $U$, $v_{\xi}$, $v_{\eta}$ 
and $a$ are expected to 
become independent of $\tau$ which is denoted by a 
superscript $\infty$. 
Thus (\ref{gKP5}) reduces  in this limit to 
\begin{equation}
    -a^{\infty}\left(\frac{2}{n}U^{\infty}+\xi 
    U^{\infty}_{\xi}+2\eta U^{\infty}_{\eta}\right)-v_{\xi}^{\infty}U^{\infty}_{\xi}
    -v_{\eta}^{\infty}U^{\infty}_{\eta}
    +(U^{\infty})^{n}U^{\infty}_{\xi}+\epsilon^{2}U^{\infty}_{\xi\xi\xi}
    +\lambda 
    \int_{-\infty}^{\xi}U^{\infty}_{\eta\eta}\,d\xi=0
    \label{ODE}.
\end{equation}
In contrast to gKdV, one does not get an ODE in this case, and there 
is no reason to assume that this partial differential equation 
reduces to an ODE in generic 
cases. 
There are in principle two different scenarios important in this 
context,  an algebraic or an exponential decay of the scaling factor 
$L(\tau)$. In the algebraic case, we have $L(\tau) = 
C_1\tau^{\gamma_1}$ with constants $C_{1}$, $\gamma_1<-1/3$ and thus 
$a^{\infty}=0$ as well as
\begin{equation}
   L(t) \propto (t^* - t)^{1/(3+1/\gamma_{1})}.
   \label{eq:Crit_Lt}
\end{equation}
Then 
equation (\ref{ODE}) reduces for $v_{\eta}^{\infty}=0$ 
to the equation for travelling wave 
solutions of the gKP equations in a commoving frame which has the 
unique nontrivial localized solution $Q$ (\ref{soliton}). Note that the latter 
condition is automatically satisfied for initial data with a symmetry 
with respect to $y\to-y$ as in the examples we consider here. Since 
the gKP equation is invariant under this transformation, we have 
$v_{\eta}=0$ in the studied examples. 

For exponential decay we have $L(\tau) = C_2 e^{a^{\infty}\tau}$ with 
$C_{2}=const$ and $a^{\infty}<0$. Relation (\ref{gKP4}) implies 
in this case
\begin{equation}
   L(t) \propto (t^* - t)^{1/3}.
   \label{eq:SupCrit_Lt}
\end{equation}

For the numerical implementation, the scaling factor $L$ and the 
\emph{speeds} $v_{\xi}$, $v_{\eta}$ have to be chosen in a convenient way. 
A possible choice for  the latter is to fix the single (by 
assumption) global minimum of $U$ at 
$\xi=\eta=0$ which 
implies  $U_{\xi}^{0}=U_{\eta}^{0}=0$, where the superscript 0 
denotes that the function is taken for $\xi=\eta=0$. These conditions 
lead to
\begin{equation}
    \begin{pmatrix}
        U_{\xi\xi}^{0} & U_{\xi\eta}^{0} \\
        U_{\xi\eta}^{0} & U_{\eta\eta}^{0}
    \end{pmatrix}
    \begin{pmatrix}
        v_{\xi} \\
        v_{\eta}
    \end{pmatrix}=
    \begin{pmatrix}
        (U^{0})^{n}U_{\xi\xi}^{0}+U_{\xi\xi\xi\xi}^{0}+\lambda 
	U_{\eta\eta}^{0} \\
        (U^{0})^{n}U_{\xi\eta}^{0}+U_{\xi\xi\xi\eta}^{0}+\lambda 
	\partial_{\xi}^{-1}U_{\eta\eta\eta}^{0}
    \end{pmatrix}
    \label{vxy}.
\end{equation}
In the numerical examples we study in this paper, it turns out that 
$x_{m}$ does not change much, whereas $y_{m}$ is constant for the 
reasons explained above. This is in contrast to gKdV in the 
$L_{2}$ critical case where it was proven in \cite{MMR2012_I} that 
$x_{m}\to\infty$ for $\tau\to\infty$. Since the computation of the 
derivatives in (\ref{vxy}) at each stage is expensive in two 
dimensions, we do not fix $x_{m}$ and $y_{m}$ in the shown examples. 

The coordinate transformation (\ref{gKP5}) implies  
for the $L_{2}$ norm of $u$
\begin{equation}
    ||u||_{2}^{2}=L^{3-4/n}||U||_{2}^{2}
    \label{L2KP}.
\end{equation}
Thus the $L_{2}$ critical case is $n=4/3$ as already mentioned. 
The $L_{2}$ norm of $u_{x}$ scales as
\begin{equation}
        ||u_{x}||_{2}^{2}=L^{1-4/n}||U_{\xi}||_{2}^{2}
    \label{L2xKP},
\end{equation}
which implies that it is invariant under the transformation 
(\ref{gKP4}) for $n=4$. Since the blow-up theorems in \cite{MST} are 
established for the $L_{2}$ norm of $u_{y}$, we consider  
here,
\begin{equation}
        ||u_{y}||_{2}^{2}=L^{-(1+4/n)}||U_{\eta}||_{2}^{2}
    \label{L2yKP}.
\end{equation}
Fixing $||U_{\eta}||$ to be constant, we get 
\begin{equation}
    a=\frac{2n}{(4+n)(n+1)||U_{\eta}||_{2}^{2}}\int_{\mathbb{R}^{2}}
    U^{n+1}U_{\eta\eta\xi}\,d\xi\, d\eta 
    \label{aL2ueta}.
\end{equation}
This will be chosen for the numerical implementation. The quantity 
$L(\tau)$ and the physical time $t$ are computed as in \cite{KP13} 
via the trapezoidal rule. The accuracy of the numerical solution is 
controlled via the computed $L_{2}$ norm of $U$ via (\ref{L2KP}) and 
the energy
\begin{equation}
            E[U] = 
	    \frac{1}{L^{(4/n-1)}}\int_{\mathbb{R}^{2}}^{}\left(\frac{1}{2}U_\xi^2
           - \frac{1}{(n+1)(n+2)}
	   U^{n+2}-\frac{\lambda}{2}(\partial_{\xi}^{-1}U_{\eta})^{2}\right)d\xi\,d\eta
    \label{energyL}.
\end{equation}
Note that the energy is invariant under the rescaling (\ref{gKP4}) as 
the $L_{2}$ norm of $u_{x}$ for $n=4$. In this sense the case $n=4$ 
is energy critical.

The spatial dependence in equation (\ref{gKP5}) as well as the time 
dependence will be treated as for the direct integration of gKP 
outlined in the previous subsection.  
The  numerically problematic terms in (\ref{gKP5}) for the 
Fourier approach are $\xi U_{\xi}$ and $\eta U_{\eta}$. In \cite{KP13} 
the dispersive oscillations with slow decrease towards 
infinity caused for gKdV numerical errors at the boundaries 
via a pollution of the Fourier coefficients at the high 
wave numbers. The problem could be solved in this case by using a very 
high resolution in time to essentially propagate the solution with 
machine precision. But in two spatial dimensions, such a high time 
resolution is computationally expensive. For gKP there is the additional problem of the 
algebraic tails due to the antiderivative. The effects of these tails 
on the Fourier coefficients are much worse than those of the 
oscillations of small amplitude. As will be shown below, they make it 
impossible to compute long enough with the rescaled codes to clearly 
identify blow-up. Nonetheless the rescaling is important since we can 
use it via a postprocessing of the directly integrated solutions to 
identify the function $L(t)$ from computed norms of the solution. 

\subsection{Numerical tests}

\begin{figure}[ht]
   \centering
   \includegraphics[width=0.49\textwidth]{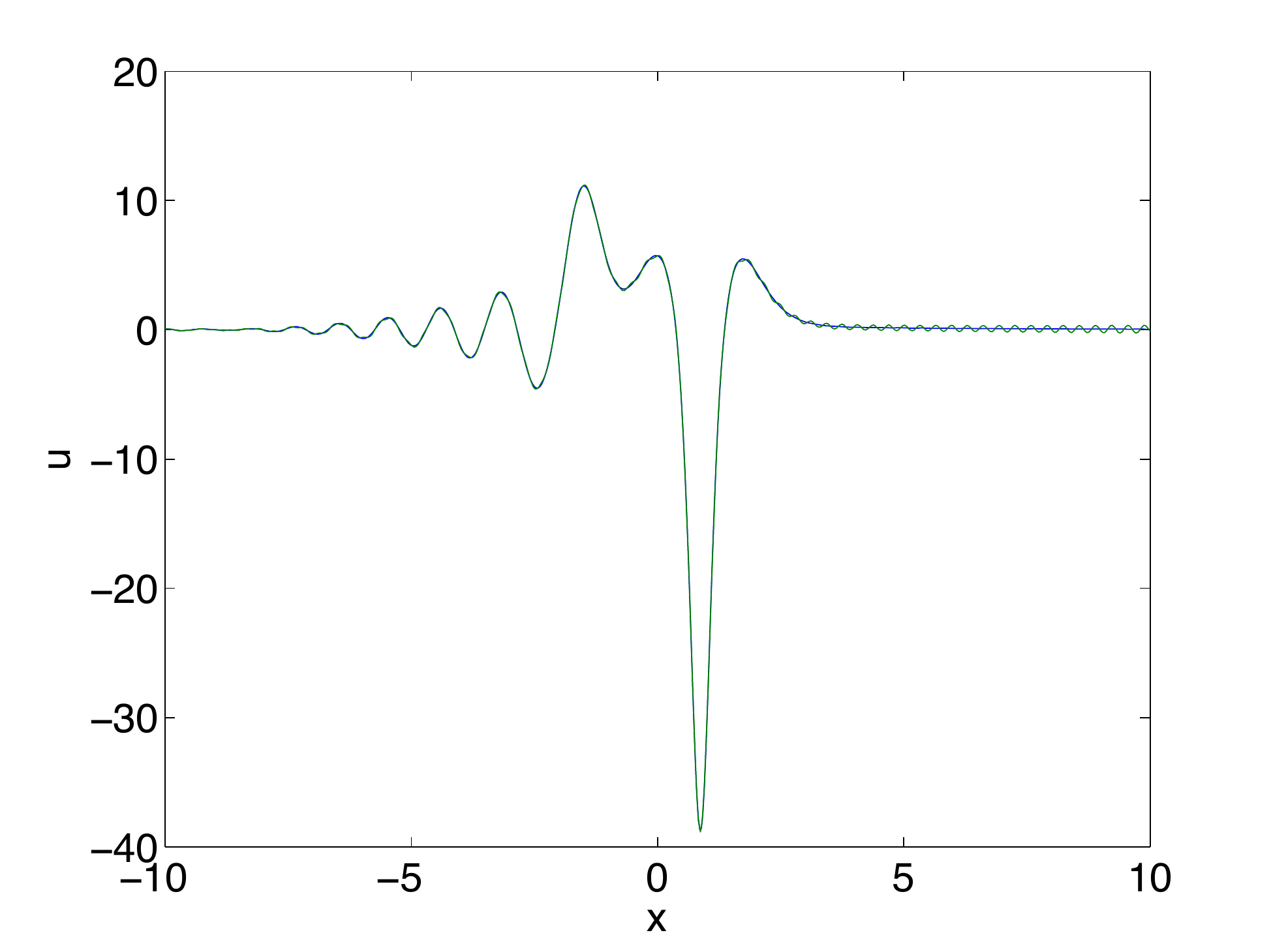}
   \includegraphics[width=0.49\textwidth]{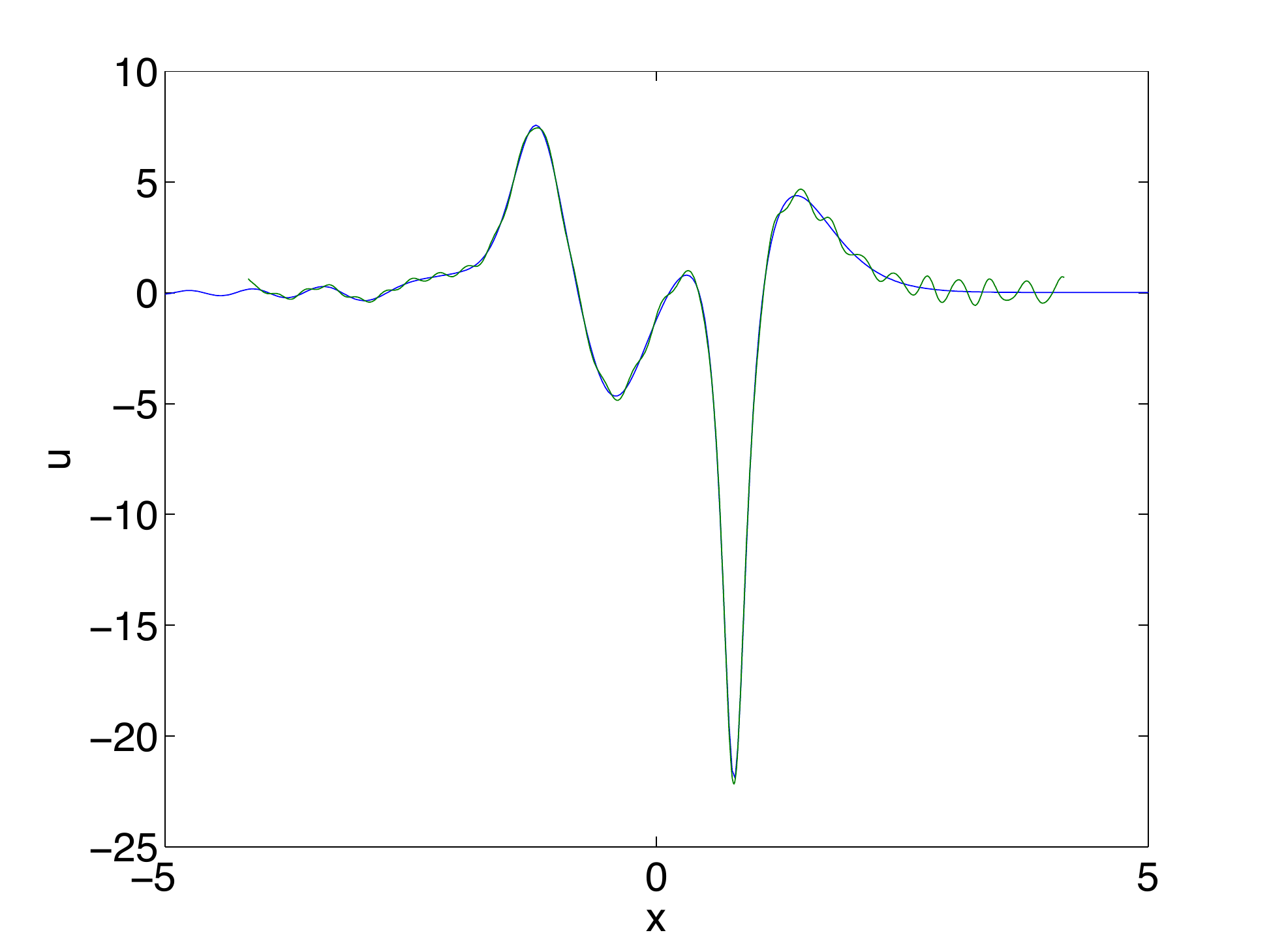}
   \caption{Solution to the rescaled gKP I equation (\ref{gKP5}) for $n = 
   4/3$ and the initial data  
   $U_{0}=12\,\partial_{\xi\xi}\exp(-\xi^{2}-\eta^{2})$ for $\tau=0.1$ 
   on the $x$-axis
   in green and the corresponding solution to the gKP I equation 
   (\ref{gKP}) in blue on the left; on the right the same setting for
    $n=2$ and the 
   initial data $U_{0}=6\,\partial_{\xi\xi}\exp(-(\xi^{2}+\eta^{2}))$ at 
   $\tau=0.5$. }
   \label{2codes}
\end{figure}

In 
addition we can compare the results of both codes  which 
provides an useful test since the codes are independent.
To this end we run the rescaled code for the examples with blow-up in 
the following sections as long as the numerically computed 
mass is conserved to better than $10\%$. Then we directly integrate 
gKP for the same initial data to the physical time corresponding to 
the end value of $\tau$ for the rescaled solution. The resulting 
solutions for the case $n=4/3$ studied in Fig.~\ref{gKPn43tau1} are shown on 
the $x$-axis in Fig.~\ref{2codes}. We present both solutions in one 
figure on the $x$-axis where they differ the most (the solution to 
(\ref{gKP5}) is rescaled back to $u$). Both codes are run 
with $N_{x}=N_{y}=2^{10}$ Fourier modes, and $L_{x}=L_{y}=5$ for the 
direct integration (the values for the rescaled code are given in the 
following sections). It can be seen that the agreement is much better 
than indicated by the conservation of the numerically computed 
mass. The disagreement is mainly due to the imposed periodicity which 
poses a problem for the rescaled code since the dispersive radiation 
reenters the dynamically rescaled domain from the other side. In the 
same way we study the case $n=2$ presented in Fig.~\ref{gKPn26gausstau05}. 
Again the main disagreement is due to the dispersive oscillations 
reentering the domain from the right. But the good agreement between 
both codes except for these oscillations provides the wanted test of 
the numerical approaches. 

\section{The $L_{2}$ critical case $n = 4/3$}
In this section we will  study solutions to the gKP equations with 
the $L_{2}$ critical nonlinearity $n=4/3$. This nonlinearity is not 
very relevant for applications, but mathematically interesting since 
it corresponds to the gKdV case $n=4$ which was intensively studied 
in \cite{MMR2012_I} and references therein. A similar theoretical 
description could be possible in this case. In addition it is 
different from the supercritical cases $n>4/3$, and thus will be 
treated in more detail. The third root of $u$ is computed in standard 
way as $u^{1/3}=\mbox{sign}(u)|u|^{1/3}$ to ensure that the real 
branch of the root is taken.

\begin{figure}[ht]
   \centering
   \includegraphics[width=0.49\textwidth]{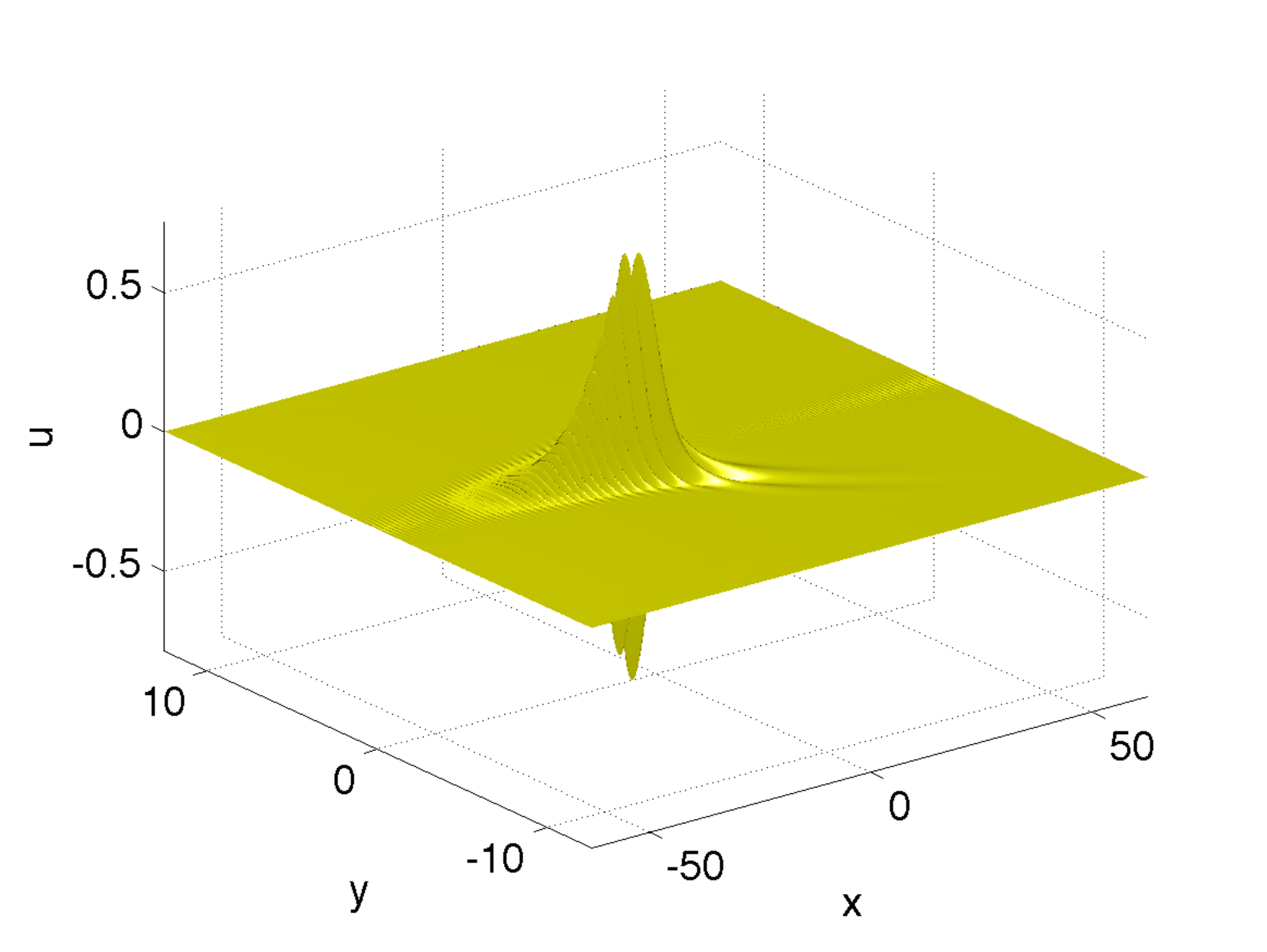}
   \includegraphics[width=0.49\textwidth]{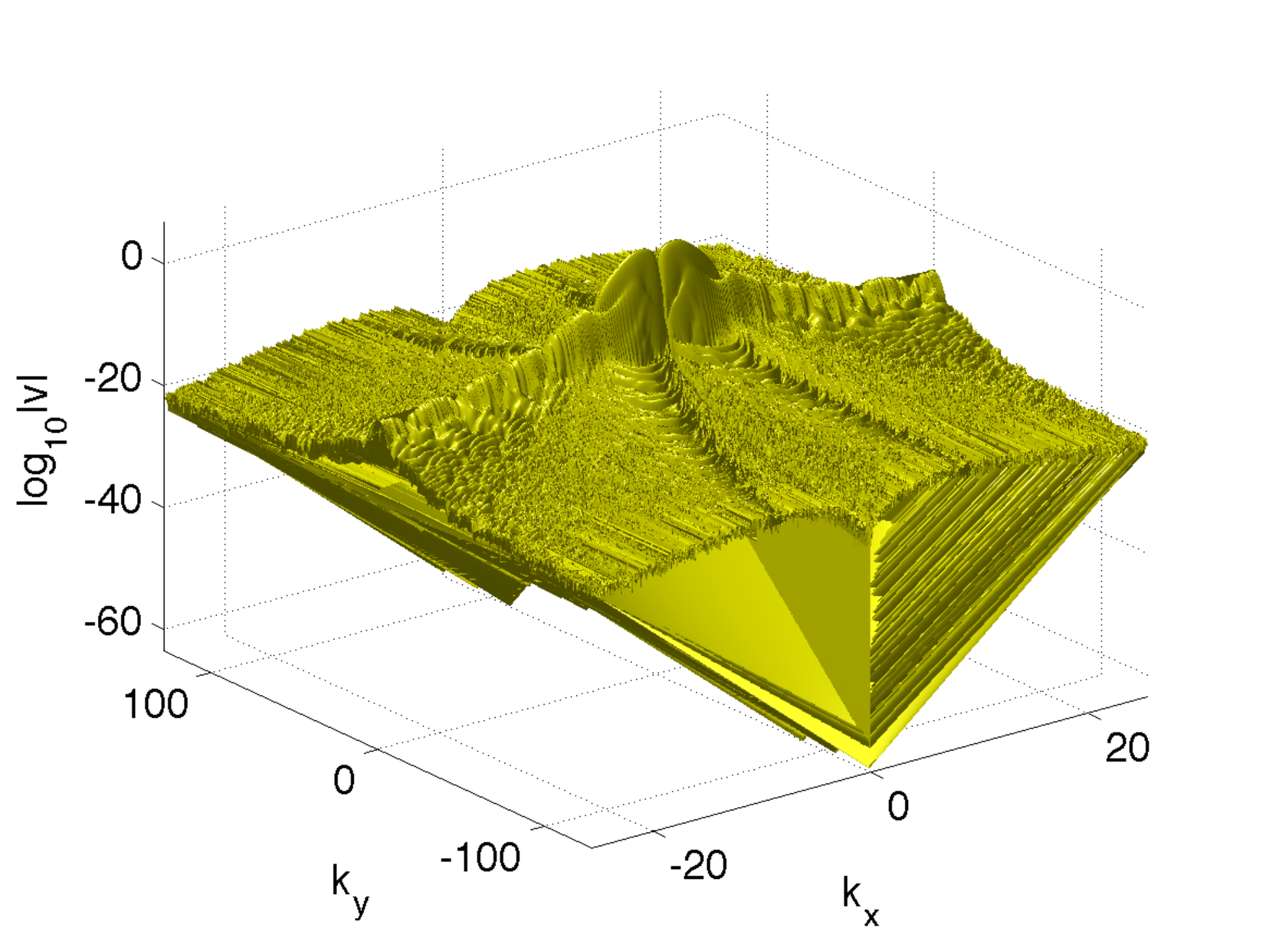}
   \caption{Solution to the gKP I equation (\ref{gKP}) for $n = 
   4/3$ and the initial data  
   $u_{0}=\partial_{xx}\exp(-x^{2}-y^{2})$ for 
   $t=0.5$ on the left; the modulus of the corresponding Fourier 
   coefficients on the right.}\label{gKPn43gaussu}
\end{figure}

We first study for gKP I initial data with positive energy, namely 
$u_{0}=\partial_{xx}\exp(-x^{2}-y^{2})$. The computation is carried 
out with $L_{x}=20$, $L_{y}=4$, $N_{x}=N_{y}=2^{10}$ and $N_{t}=1000$ 
time steps for $0<t\leq 0.5$. The relative computed energy is conserved 
to better than $10^{-5}$.  As can be seen in 
Fig.~\ref{gKPn43gaussu}, oscillations propagating to the left form, and 
the initial hump appears to be just radiated away. Due to the imposed 
periodicity, these oscillation reenter on the right. The algebraic 
tails to the right can be also clearly recognized. A consequence of 
these tails is that a slight Gibbs phenomenon 
appears which can be seen from the Fourier coefficients in 
the same figure. The solution is well resolved in 
$y$-direction, but the Fourier coefficients in $k_{x}$ no longer 
decrease to machine precision (but $10^{-5}$ is more than sufficient 
for our purposes).

\begin{figure}[ht]
   \centering
   \includegraphics[width=0.49\textwidth]{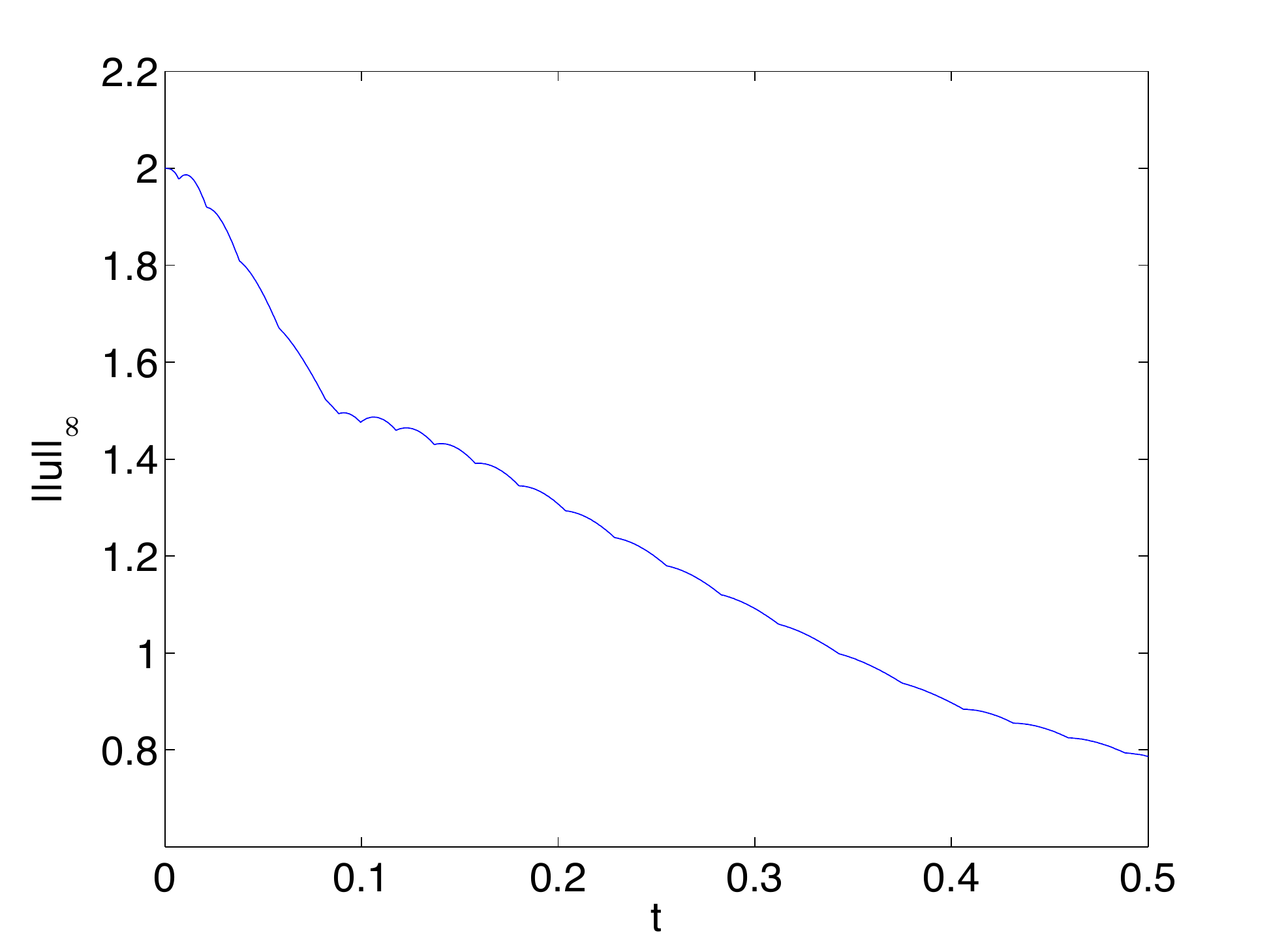}
   \includegraphics[width=0.49\textwidth]{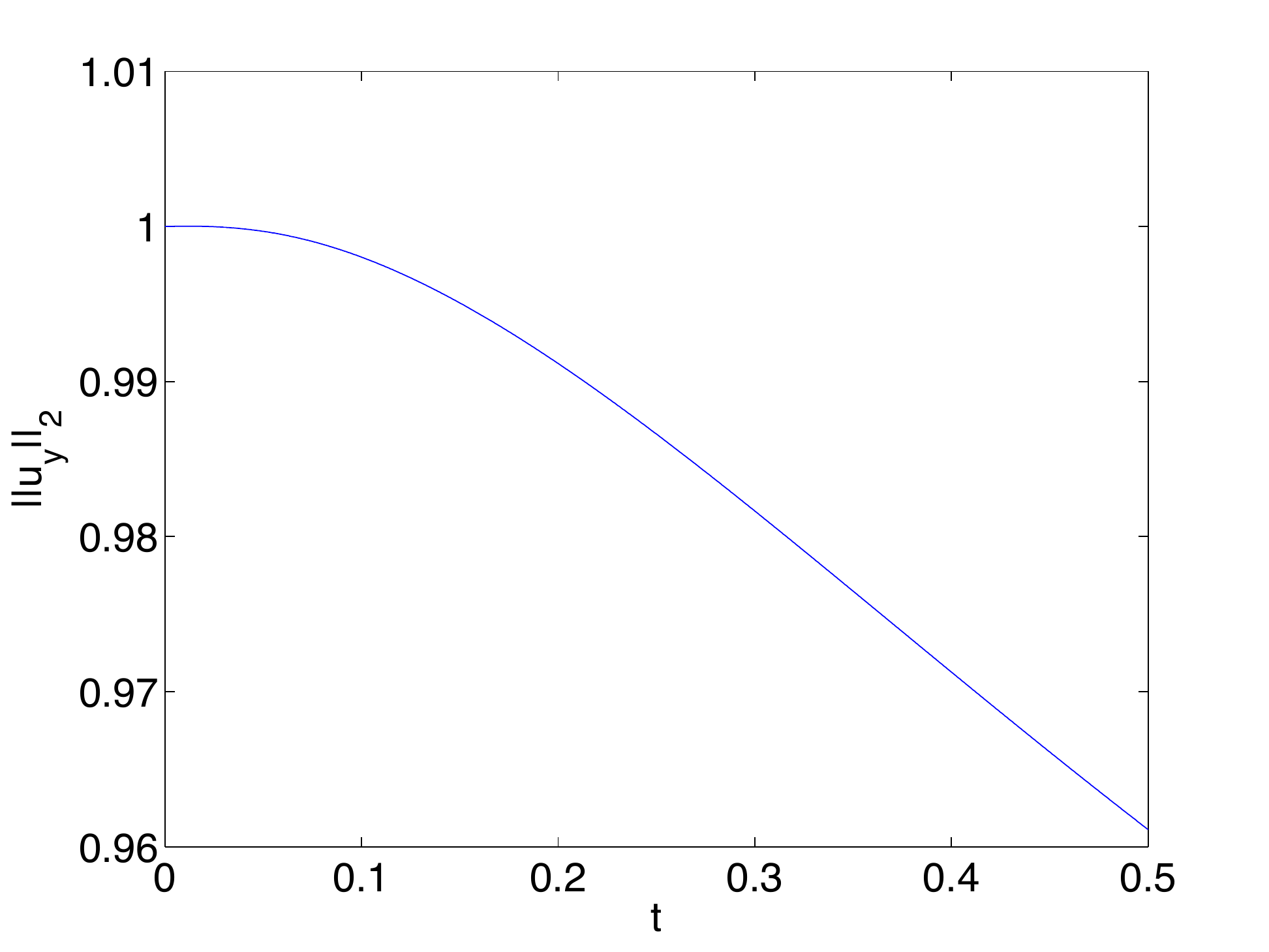}
   \caption{Norms of the 
   solution to the gKP I equation (\ref{gKP}) for $n = 
   4/3$ and the initial data  
   $u_{0}=\partial_{xx}\exp(-x^{2}-y^{2})$ in dependence of $t$; on 
   the left $||u||_{\infty}$, on the right $||u_{y}||_{2}$ normalized 
   to 1 at $t=0$.}
   \label{gKPn43gaussuinf}
\end{figure}

There is no indication of blow-up, and this is even more obvious from 
the norms shown in Fig.~\ref{gKPn43gaussuinf}. Both $||u||_{\infty}$ 
and $||u_{y}||_{2}$ appear to be monotonically decreasing. This 
suggests again that the initial hump will be just radiated away 
towards 
infinity.

\begin{figure}[ht]
   \centering
   \includegraphics[width=\textwidth]{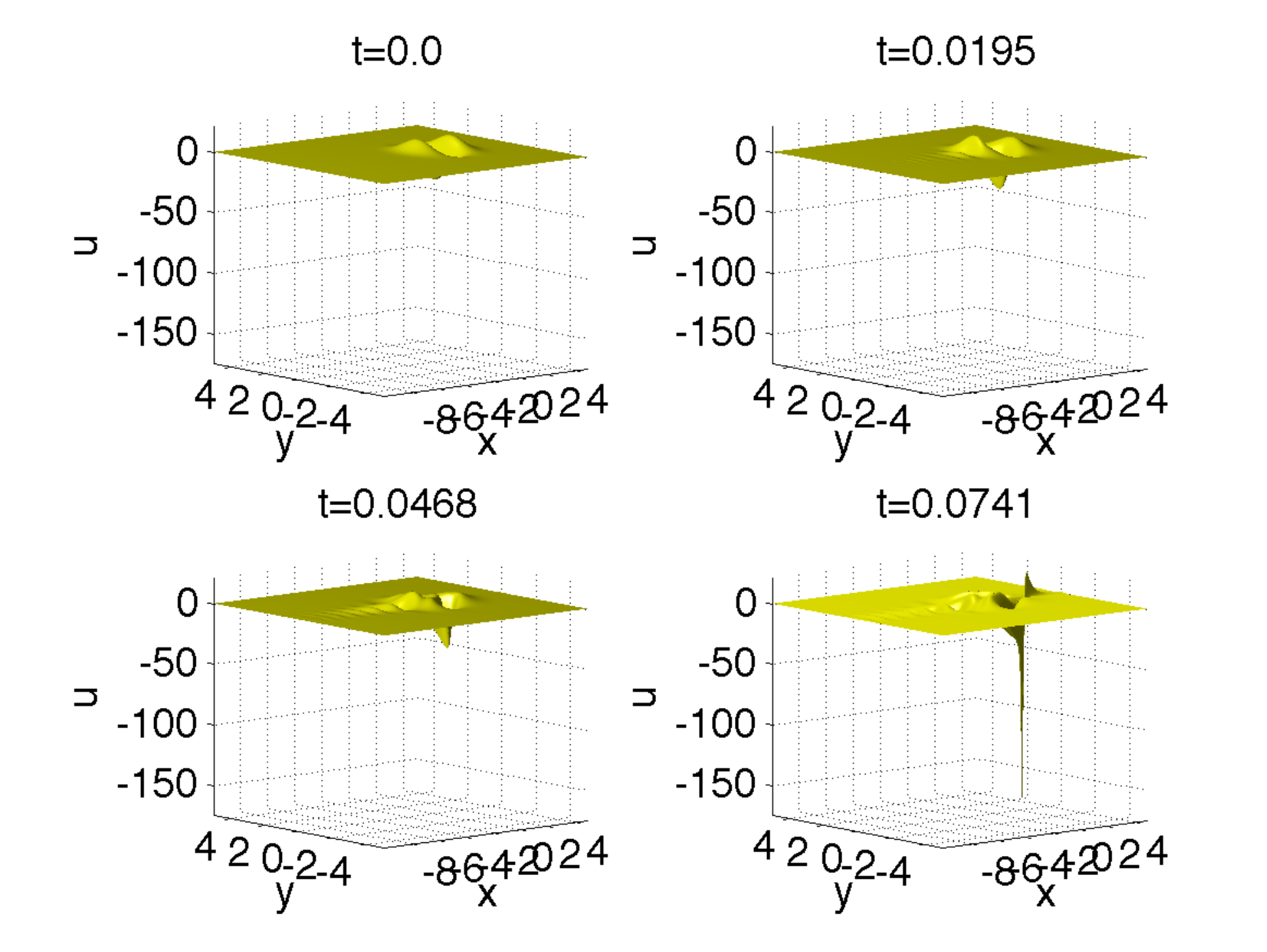}
   \caption{Solution to the gKP I equation (\ref{gKP}) for $n = 
   4/3$ and the initial data  
   $u_{0}=12\,\partial_{xx}\exp(-x^{2}-y^{2})$ for several values of 
   $t$.}\label{fig:n4over3_u}
\end{figure}

The situation is quite different for initial data with negative 
energy, $u_{0}=12\,\partial_{xx}\exp(-x^{2}-y^{2})$. The calculation is 
performed with $L_x = L_{y}=5$, $N_x = 2^{10}$, $N_y = 2^{13}$, and 
$N_t = 50000$ time steps for $0\leq 
t<0.078$. As can be seen in Fig.~\ref{fig:n4over3_u}, there are 
dispersive oscillations forming immediately and propagating to the left. 
Due to the imposed periodicity condition, these oscillations reenter 
at a given time on the right (note that only part of the 
computational domain is shown in the figure). But more importantly 
one can see that the initial minimum gets more and more peaked 
and finally appears to blow up in a point. The code is stopped once 
$\Delta>10^{-3}$. 

\begin{figure}[ht]
   \centering
   \includegraphics[width=.49\textwidth]{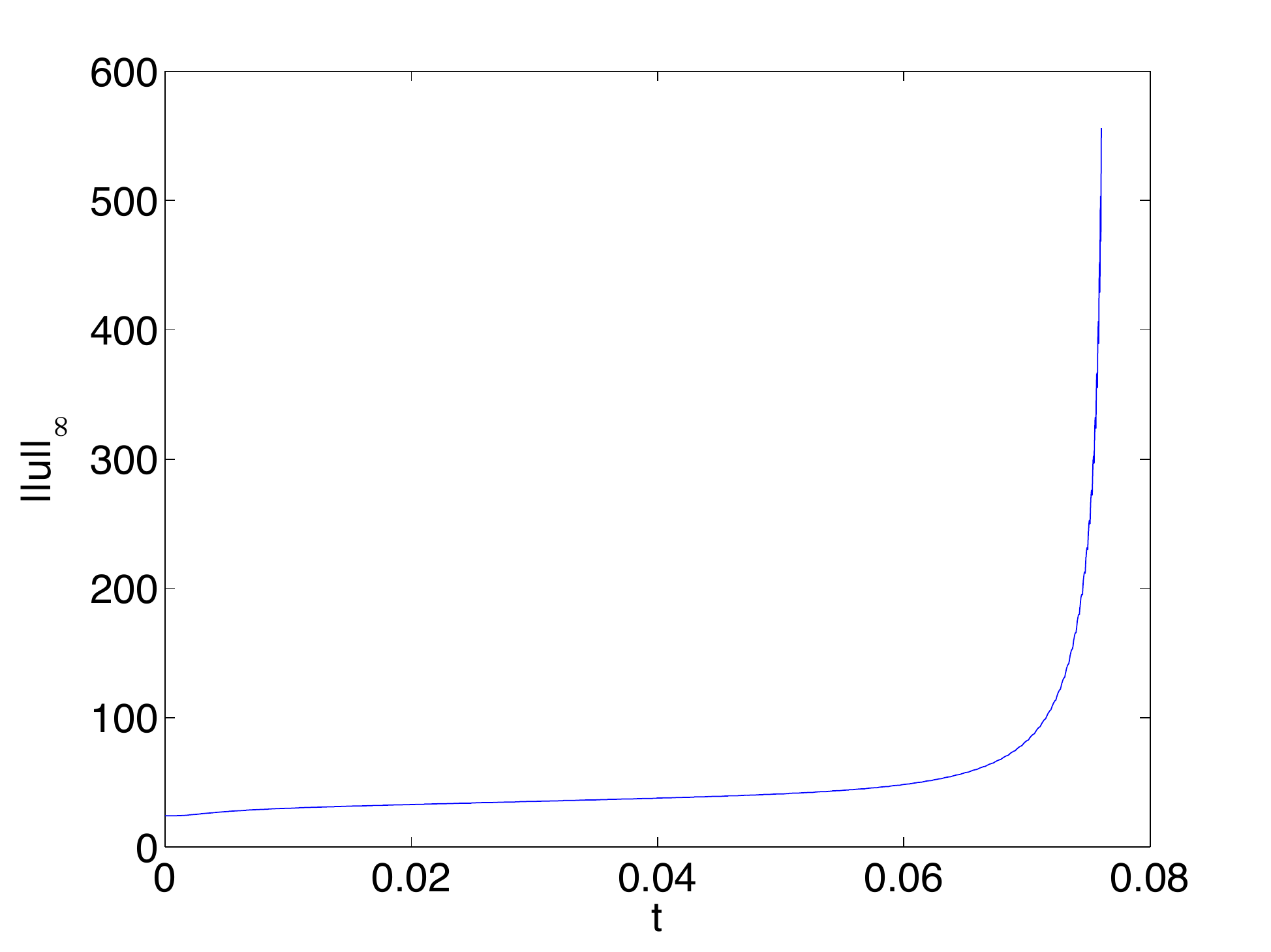}
   \includegraphics[width=.49\textwidth]{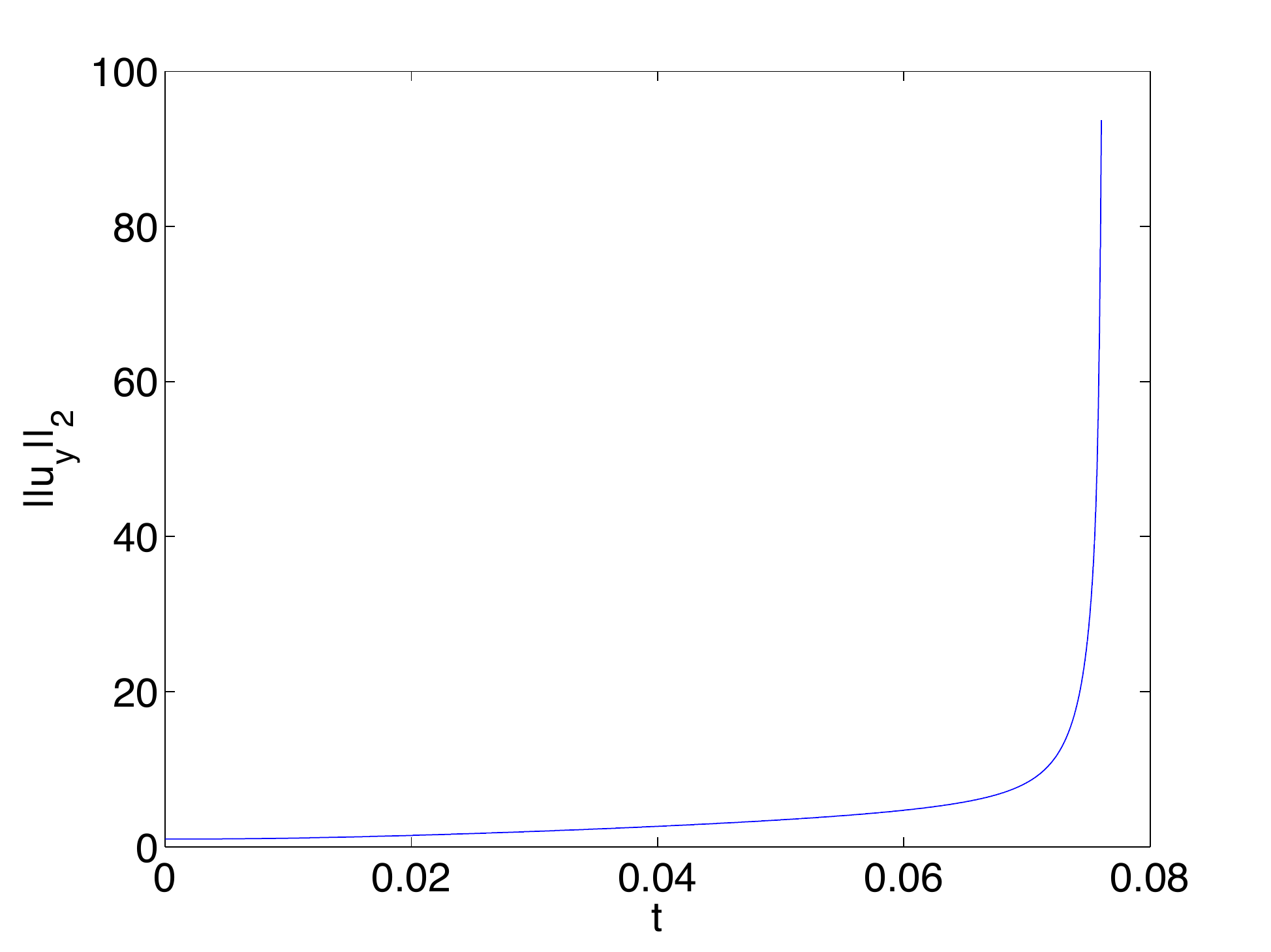}
   \caption{Norms of the solution to the gKP I equation (\ref{gKP}) for $n = 
   4/3$ and the initial data  
   $u_{0}=12\,\partial_{xx}\exp(-x^{2}-y^{2})$ in dependence of $t$; on 
   left the $L_{\infty}$ norm of $u$, on the right the $L_{2}$ norm 
   of $u_{y}$.}\label{fig:n4over3_L2Norm_uy}
\end{figure}

\begin{figure}[ht]
   \centering
   \includegraphics[width=0.49\textwidth]{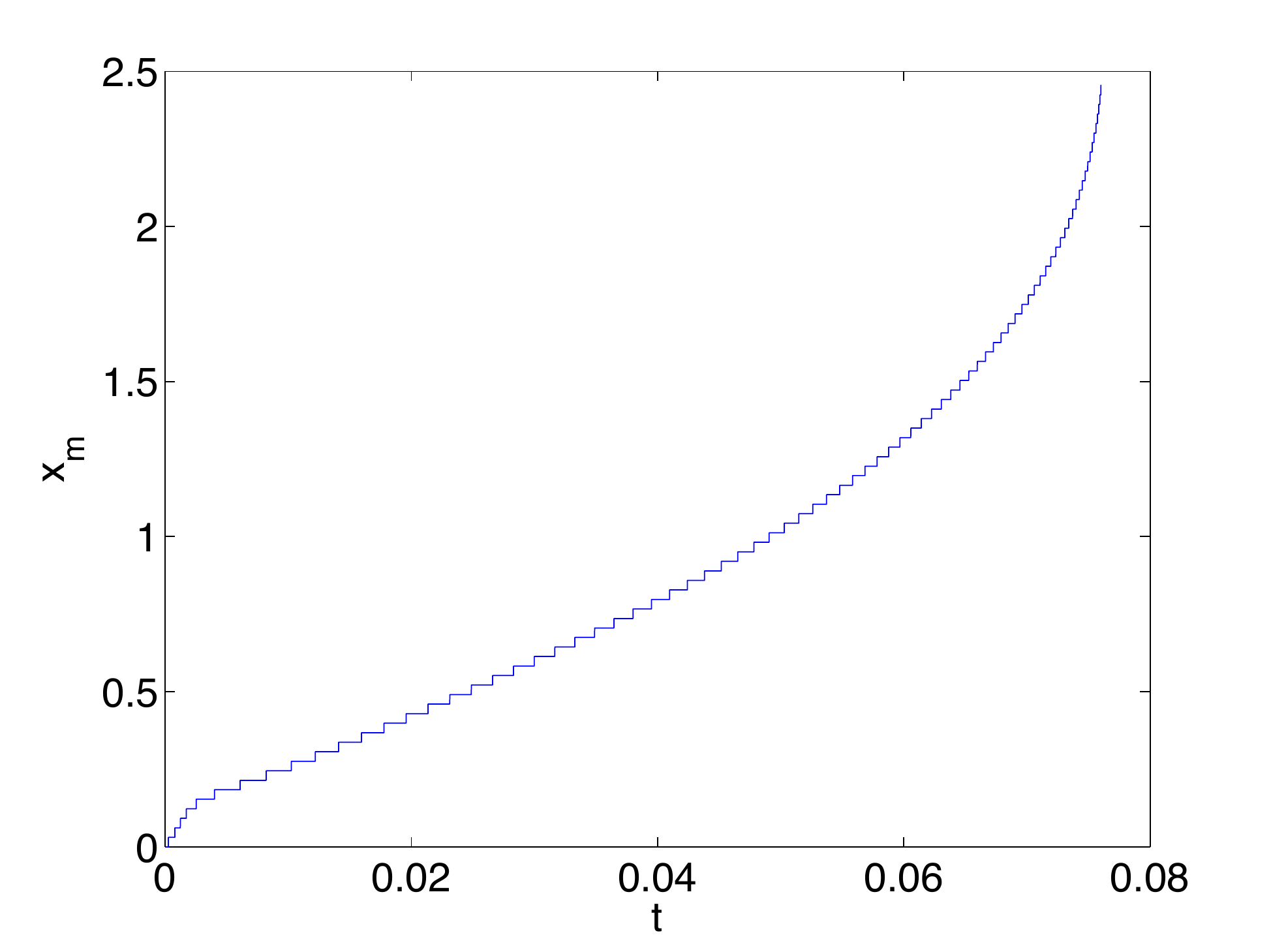}
   \includegraphics[width=0.49\textwidth]{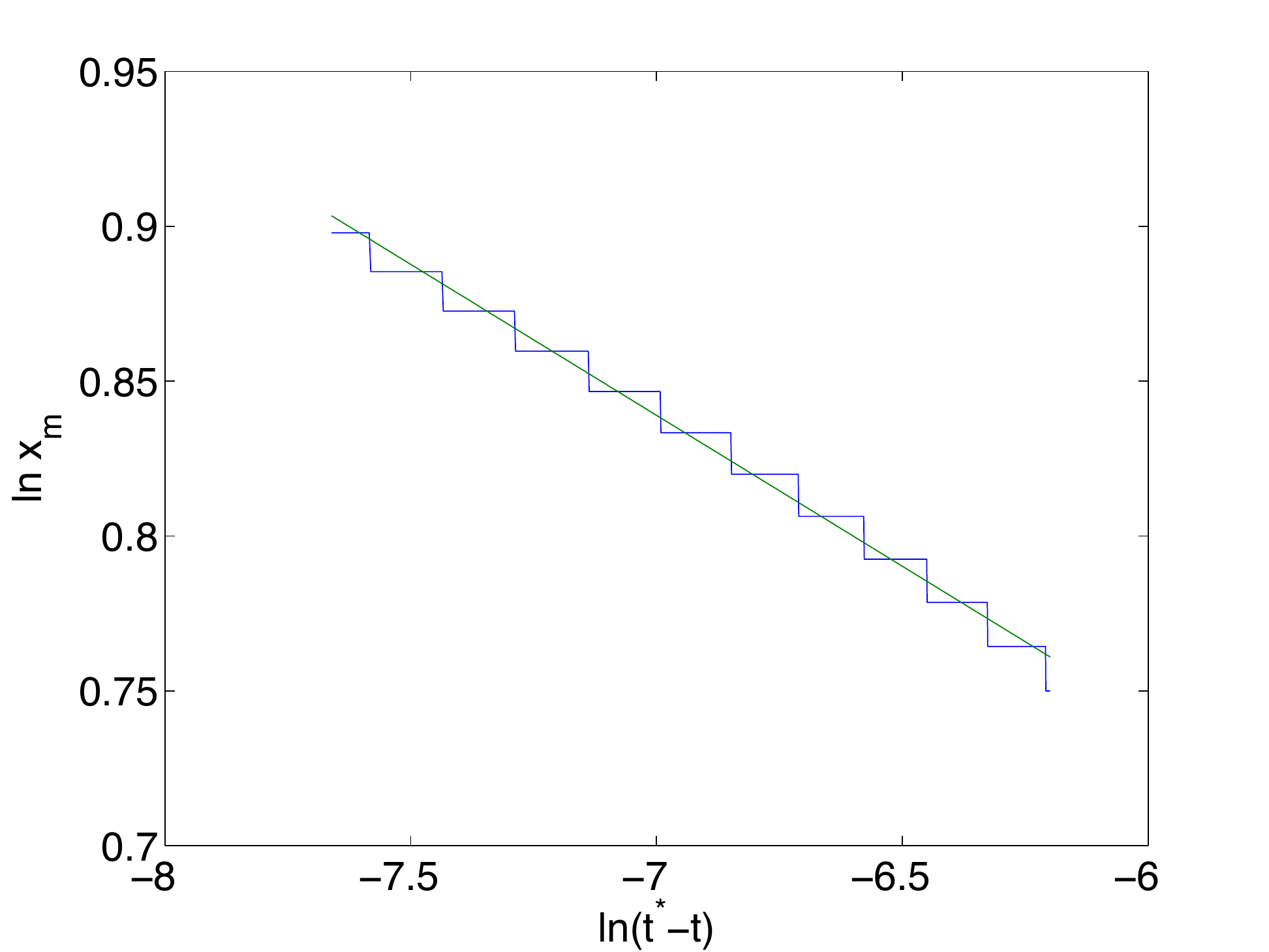}
   \caption{Location $x_{m}$ of the global minimum of the 
   solution to the gKP I equation (\ref{gKP}) for $n = 
   4/3$ and the initial data  
   $u_{0}=12\,\partial_{xx}\exp(-x^{2}-y^{2})$
   in dependence of 
   $t$ on the left; on the right $\ln x_{m}$ for $t\sim t^{*}$ in 
   blue and a least square fit $\alpha_{1}\ln (t^{*}-t)+\alpha_{2}$ 
   in green.}\label{gKPn4312gaussxm}
\end{figure}

The blow-up is even more obvious from the norms $||u||_{\infty}$ and 
$||u_{y}||_{2}$ shown in Fig.~\ref{fig:n4over3_L2Norm_uy} which 
clearly seem to diverge. The Fourier coefficients of the solution in Fig.~\ref{fig:n4over3_u} 
at the last recorded time are shown in Fig.~\ref{gKPn4312gaussxm}. As 
expected from the rescaling (\ref{gKP4}), near a blow-up the 
$y$-derivative diverges more strongly than the $x$-derivative of the 
solution. This implies that the Fourier coefficients 
decrease much more slowly in $k_{y}$ than in $k_{x}$. 
This is why a higher resolution 
in $k_{y}$ was chosen from the beginning. But this behavior will 
persist with even higher resolution in $k_{y}$ due the scaling of the 
solution close to blow-up. On the other hand the resolution 
in time appears to be sufficient since the results for half the 
number of time steps used here do not differ from the ones shown 
within the numerical limits imposed by the Fourier coefficients.
\begin{figure}[ht]
   \centering
   \includegraphics[width=0.6\textwidth]{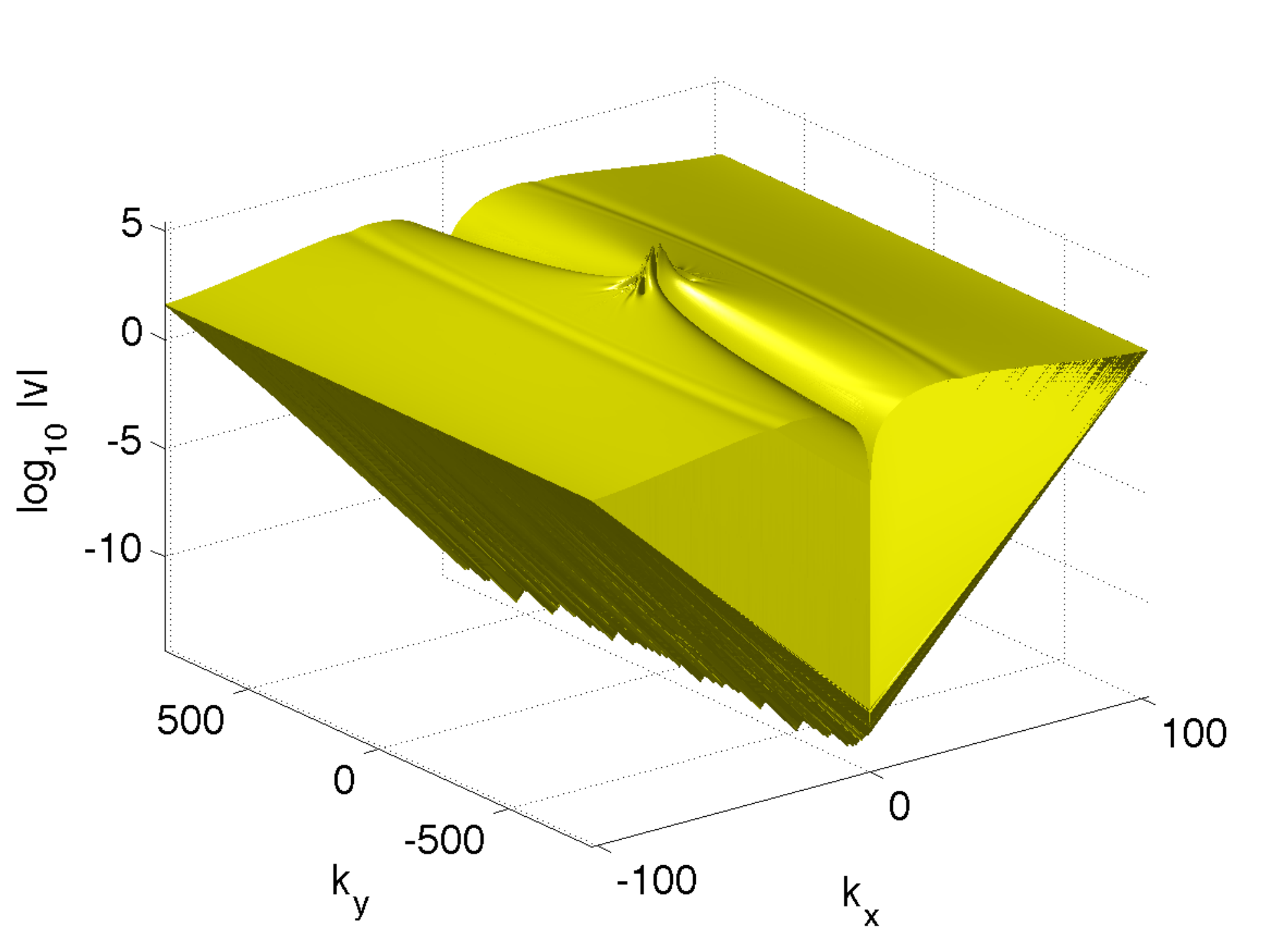}
   \caption{Modulus of the Fourier coefficients for the 
   solution to the gKP I equation (\ref{gKP}) for $n = 
   4/3$ and the initial data  
   $u_{0}=12\,\partial_{xx}\exp(-x^{2}-y^{2})$ for $t=0.0741$.}
   \label{gKPn4312gaussfourier}
\end{figure}

\begin{figure}[ht]
   \centering
   \includegraphics[width=0.49\textwidth]{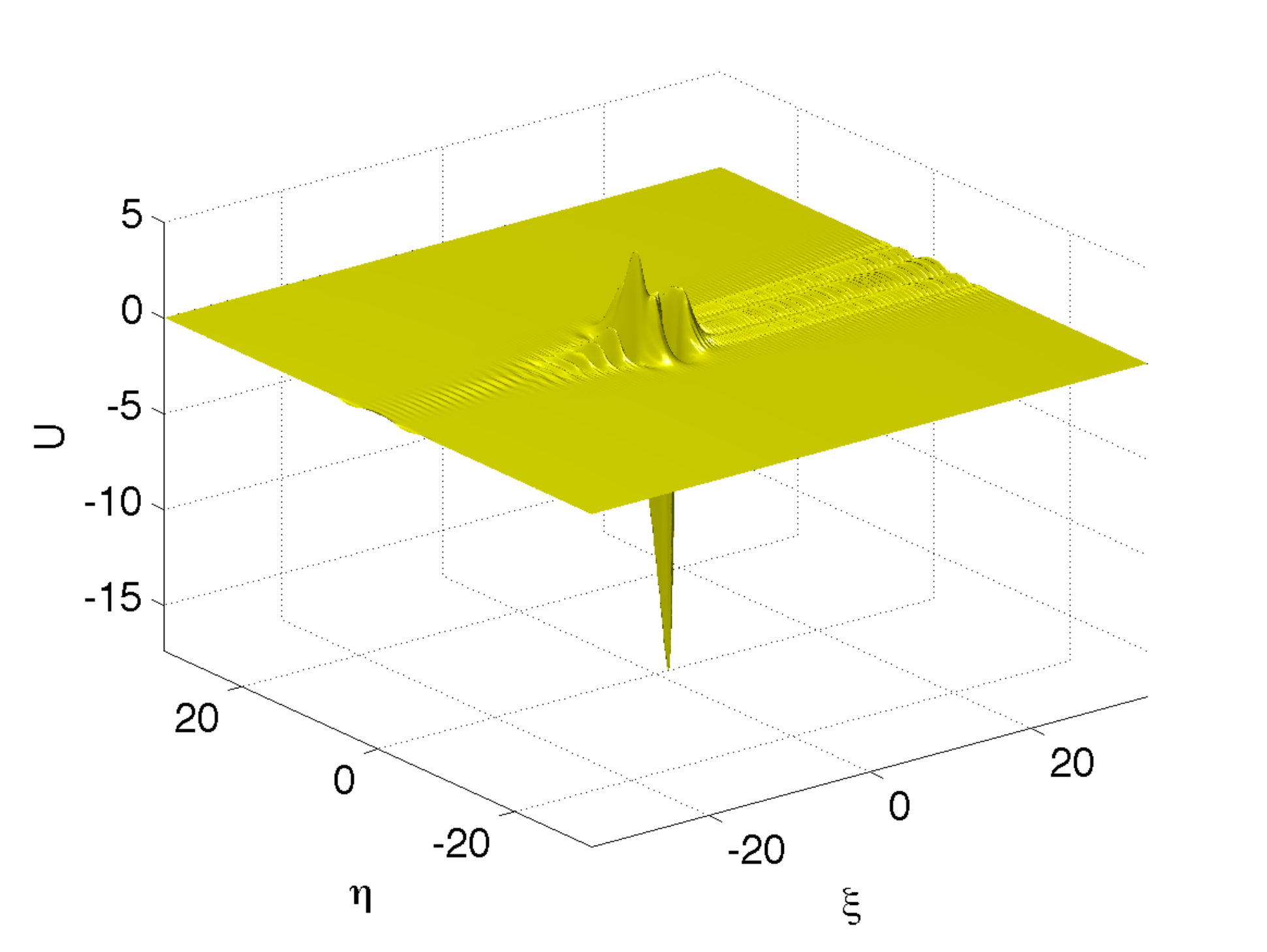}
   \includegraphics[width=0.49\textwidth]{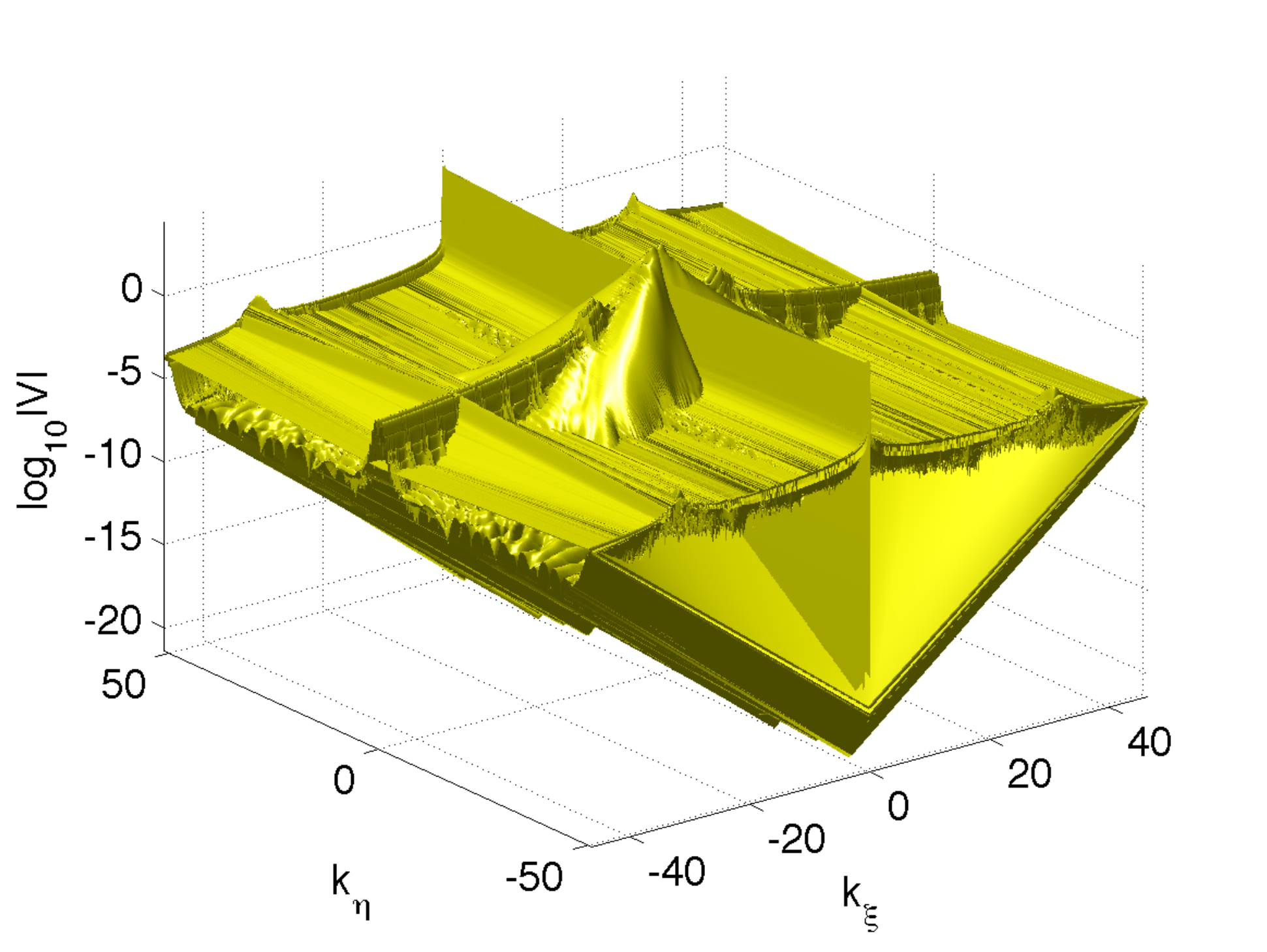}
   \caption{Solution to the rescaled gKP I equation (\ref{gKP5}) for $n = 
   4/3$ and the initial data  
   $U_{0}=12\,\partial_{\xi\xi}\exp(-\xi^{2}-\eta^{2})$ for $\tau=0.1$ 
   on the left, and the modulus of the corresponding Fourier 
   coefficients on the right.}
   \label{gKPn43tau1}
\end{figure}

To understand the type of this blow-up better,  we try to solve the 
rescaled equation (\ref{gKP5}) for these initial data. With 
$L_{\xi}=11$, $L_{\eta}=10$, $N_{\xi}=N_{\eta}=2^{10}$ and $N_{\tau}=10^{4}$ time 
steps for $0<\tau\leq0.1$, we can solve (\ref{gKP5}) with a 
relative conservation of the computed mass of the order of $10^{-1}$. The solution at the final time can be seen in 
Fig.~\ref{gKPn43tau1}.
It can be recognized that the dispersive oscillations and even more so the 
algebraic fall off of the solution lead to a Gibbs phenomenon at the 
boundary of the computational domain. This is reflected in the 
Fourier coefficients in Fig.~\ref{gKPn43tau1}, and by the fact that 
energy conservation is no longer given in this case.  The problematic 
aspect is here the increase in the amplitude of the high wave numbers 
in $k_{\eta}$. These lead eventually to a breakdown of the solution. 
Thus we could not study the solution on larger domains with more 
Fourier modes since even with a high resolution in $\tau$, the code 
would crash. 

\begin{figure}[ht]
   \centering
   \includegraphics[width=.49\textwidth]{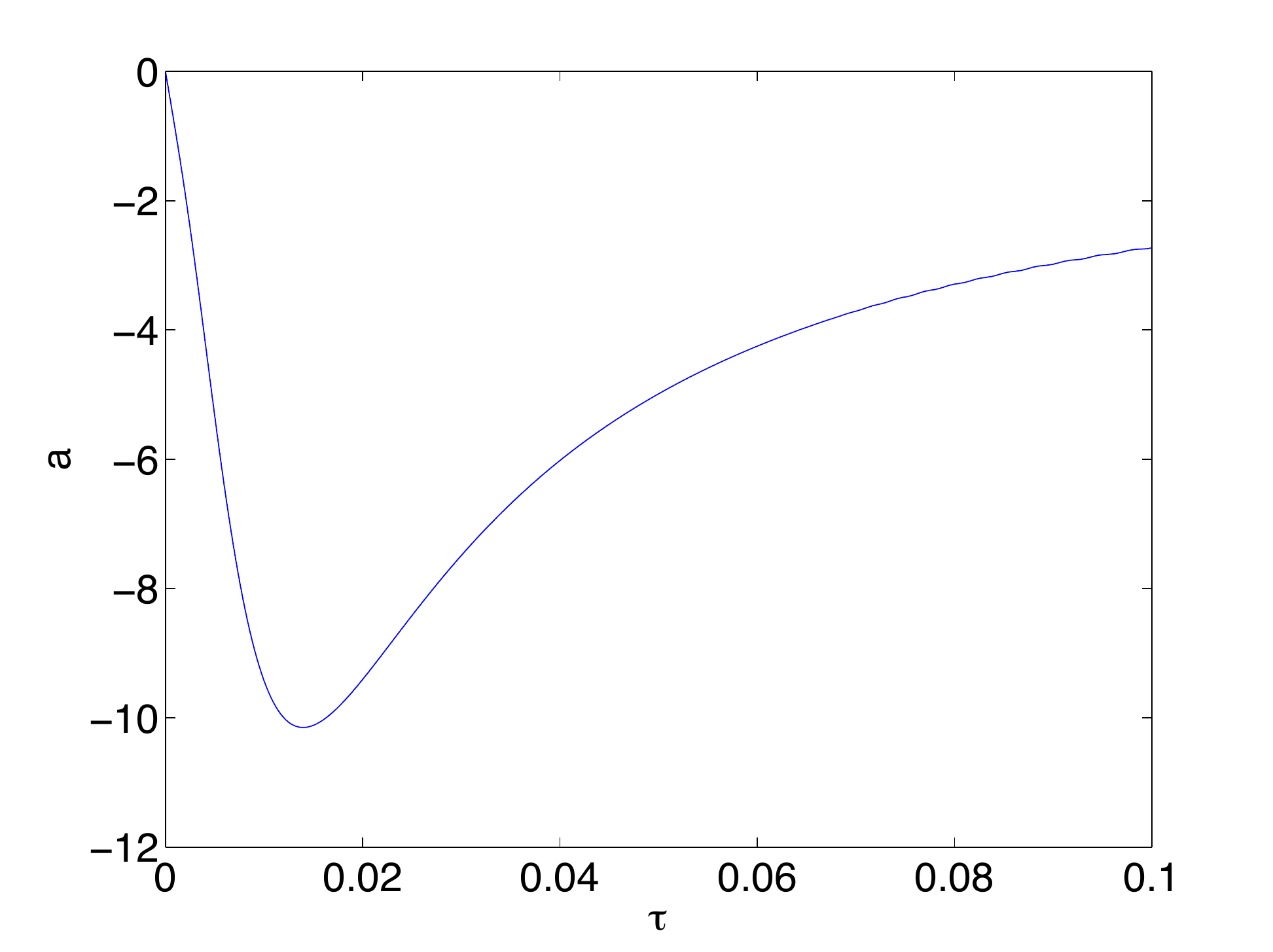}
   \includegraphics[width=.49\textwidth]{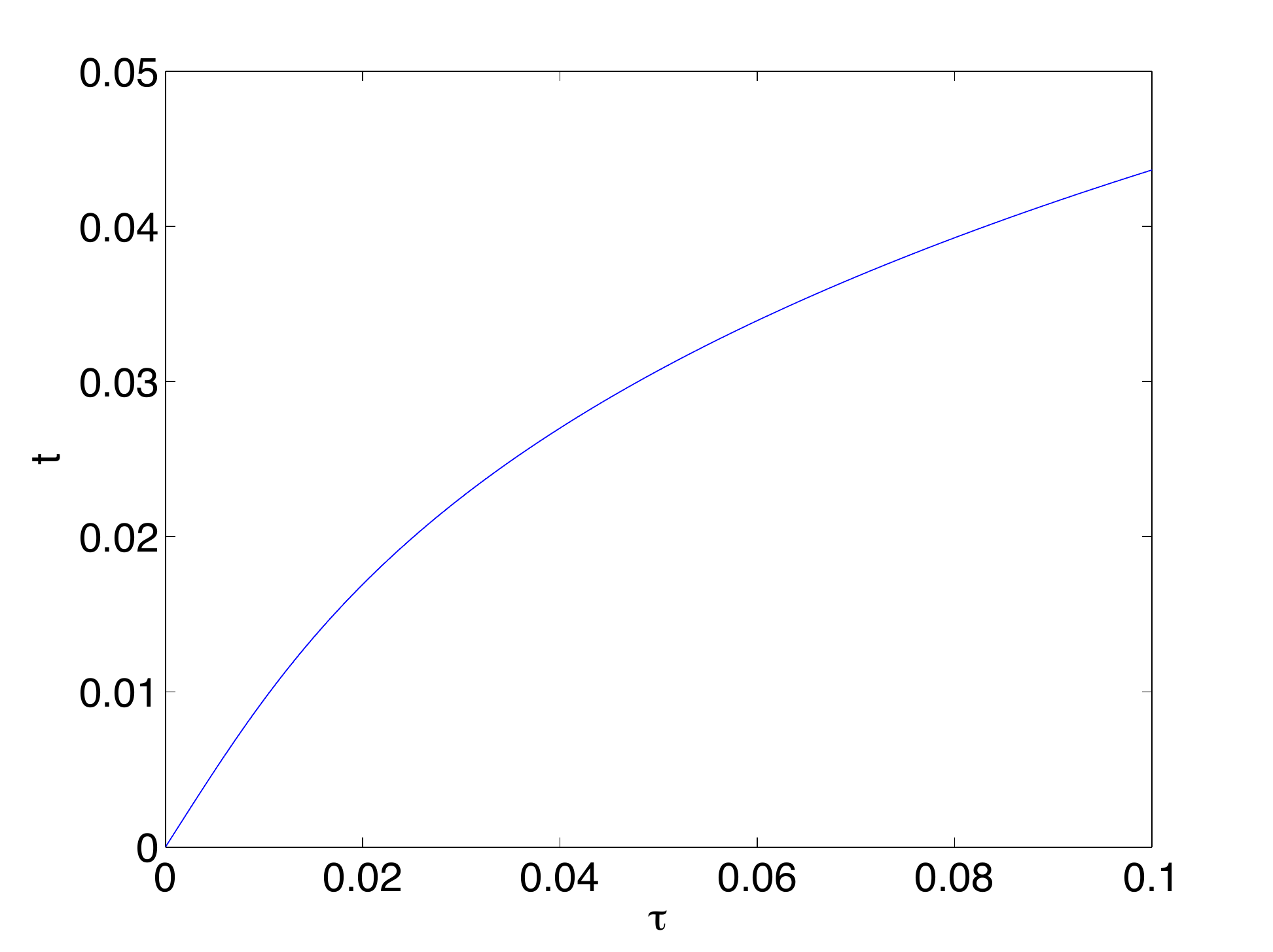}
   \caption{The quantity $a$ from (\ref{a}) and the physical time $t$ 
   in dependence of $\tau$ for the 
   solution to the rescaled gKP I equation (\ref{gKP5}) for $n = 
   4/3$ and the initial data  
   $U_{0}=12\,\partial_{\xi\xi}\exp(-\xi^{2}-\eta^{2})$ in dependence 
   of $\tau$.}
   \label{gKPn43tau1a}
\end{figure}

The function $a$ from (\ref{a}) and the physical time in dependence 
of $t$ can be seen in Fig.~\ref{gKPn43tau1a}. Obviously we do not get 
close enough to the blow-up to decide whether $a$ tends to zero or 
not. The physical time at the end of the computation also shows that 
we are not as close as necessary to $t^{*}$ to use this approach to 
discuss the type of blow-up. 

\begin{figure}[ht]
   \centering
   \includegraphics[width=0.49\textwidth]{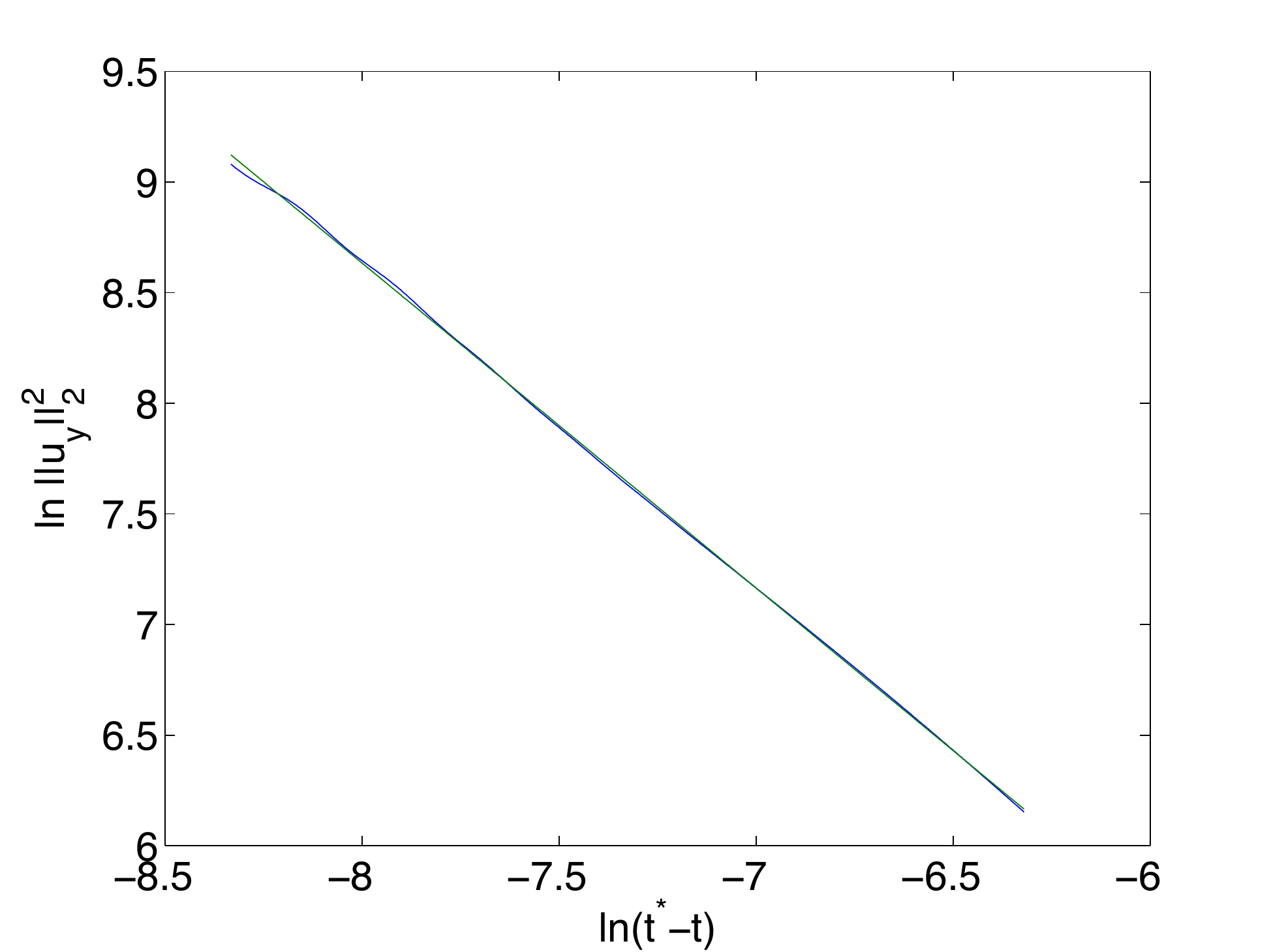}
   \includegraphics[width=0.49\textwidth]{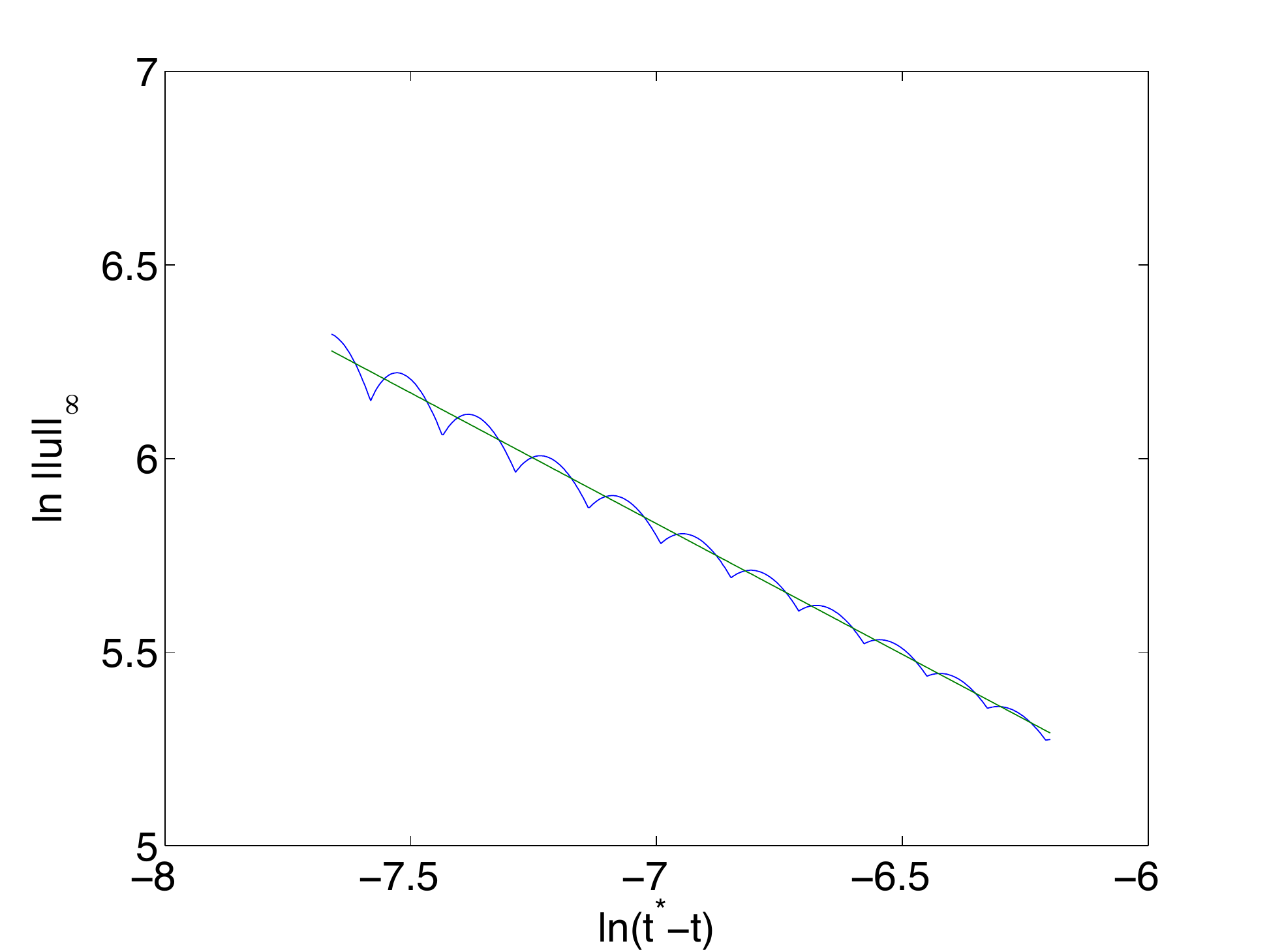}
   \caption{Fitting of $\ln||u_{y}||_{2}^{2}$ and of
 $\ln||u||_{\infty}$ (in blue) to $c \ln(t^{*}-t)+C$ (in green) 
 for the situation shown in 
 Fig.~\ref{fig:n4over3_u}.}\label{gKPn4312gaussfit}
\end{figure}

Therefore it is more promising to integrate gKP directly, and to 
infer the type of blow-up from the norms of the solution. If we 
assume that for $n=4/3$ the asymptotic behavior (\ref{eq:Crit_Lt}) 
is given, we get with  (\ref{gKP4}) as well as (\ref{L2yKP}) (recall 
that $\gamma_{1}=-1$ for gKdV in the case $n=4$) 
\begin{equation}
    ||u_{y}||_{2}^{2}\propto (t^{*}-t)^{-\frac{4}{3+1/\gamma_{1}}},\quad ||u||_{\infty}\propto 
    (t^{*}-t)^{-\frac{3}{2(3+1/\gamma_{1})}}
    \label{n43scal}.
\end{equation}
We fit $\ln ||u||_{\infty}$ for the solution shown in 
Fig.~\ref{fig:n4over3_u} to $C_1+c_{1}\ln(t^* - t)$ and 
$\ln||u_{y}||_{2}^{2}$ to $ C_2+c_{2}\ln(t^* - t)$ via the 
optimization algorithm \cite{fminsearch} distributed with Matlab as 
\textit{fminsearch}. Since the $L_{2}$ norm involves an integral over 
the whole computational domain, it is numerically more stable, and 
the fitting should be thus reliable. However, this is not the case if 
strong gradients appear in 
the solution as in the case of blow-up. It turns out that 
the $L_{\infty}$ norm typically produces results in better  accordance with 
the theory exactly since it only takes into account effects close 
to singularity. We present the  
results for both the $L_{2}$ norm of $u_{y}$ and 
the $L_{\infty}$ norm of the solution to test the 
consistency of the result.  The  fitting is shown in Fig.~\ref{gKPn4312gaussfit} for the last 1000 
computed points. We find for $||u_{y}||^{2}_{2}$ the values 
$t^{*}=0.0763$, $c_{1}=-3.1162$ and  $ C_{1}= -1.4687$, and for         
$||u||_{\infty}$ we get $t^{*}=0.0765$, $c_{2}=-0.6751$ and  $ 
C_{2}= 1.1063$. The good agreement of the blow-up time $t^{*}$ 
shows the consistency of the approaches. 
The oscillations in the $L_{\infty}$ norm suggest as already 
mentioned that the $L_{2}$ norm of $u_{y}$ should be more reliable.
Note, however, that the fitted values for $c_{1}$ and $c_{2}$ are 
with (\ref{n43scal}) compatible with $\gamma_{1}=-1$ (i.e., to $-2$ 
respectively $-3/4$), which 
corresponds to the gKdV case for $n=4$ where  
$L\propto 1/\tau$. But the agreement is better for the $L_{\infty}$ 
norm of $u$, presumable due to a loss of accuracy in the gradients at 
the boundaries of the computational domain which affect 
$||u_{y}||_{2}$, but not $||u||_{\infty}$. 
Nonetheless  the $L_{2}$ critical case seems to be for 
both gKdV and gKP characterized by an algebraic in $\tau$ vanishing 
of $L$. 

The location of the global minimum is traced as a function of time in 
Fig.~\ref{gKPn4312gaussxm}. This shows that $x_{m}$ in fact changes with 
$t$, but much less so than in the $L_{2}$ critical case for gKdV 
where blow-up happens at infinity for almost solitonic initial data. 
Since we  cannot get arbitrarily close to the actual blow-up  
for obvious reasons, it is difficult decide whether there will be blow-up   
at a finite $x^{*}$.  As already mentioned, $y_{m}=0$ for symmetry 
reasons. The plateaus in the plot Fig.~\ref{gKPn4312gaussxm} are due 
to the minimum just being evaluated on grid points. We fit $\ln 
x_{m}$ as the norms of $u$ above for the last 1000 time steps to 
$\alpha_{1}\ln (t^{*}-t)+\alpha_{2}$, where $t^{*}$ is one of the 
values determined above by fitting the norms of $u$. Doing this for 
the values obtained from the $L_{\infty}$ norm of $u$, we get  the 
figure shown in Fig.~\ref{gKPn4312gaussxm} and the values 
$\alpha_{1}=-0.0974$ and 
$\alpha_{2}=0.1569$. Doing the same fitting for the $t^{*}$ obtained 
by fitting the $L_{2}$ norm of $u_{y}$, we get $\alpha_{1}=-0.0735$ 
and $\alpha_{2}=
    0.3000$. Thus the values do not coincide, but they are both 
    negative which would imply a blow-up at infinity. As in the gKdV 
    case \cite{KP13}, the asymptotic regime for the $x_{m}$ appears 
    much later in the computations than the one for the  norms of 
    $u$, and as there, it is difficult to make reliable predictions 
    about the value of $x^{*}$. From the present computations, it 
    cannot be ruled out that the blow-up in the $L_{2}$ critical case 
    happens for infinite $x^{*}$. But the small value of 
    $|\alpha_{1}|$ (much smaller than the critical exponents for the 
    norms discussed above) makes  a finite $x^{*}$ more probable. 
    Recall that for gKdV in the $L_{2}$ critical case  $x_{m}\propto L$, see 
    \cite{MMR2012_I}.

Another interesting question is which initial data are simply 
radiated away and which blow-up. Since no stable solitons are known, 
it is not surprising that no stable structures appear in the 
computations. If we consider initial data as above of the form 
(\ref{initial}), we find that there is no blow-up for values of $\beta\leq 
6$, but that there is blow-up for $\beta\geq 7$. Numerically it is 
difficult to decide at which value of $\beta$ blow-up appears, since 
it is difficult to distinguish a strongly peaked minimum from an exploding 
norm of the solution. But the important observation is that the 
energy for initial data with $\beta=7$ is positive. Thus negative 
energy does not appear to be the necessary condition for blow-up. But 
the mass of the initial data has to be large enough, and the energy 
has to be small enough for a blow-up to appear. Numerically we can 
only give indications where for certain classes of initial data the 
transition between global existence of the solution in time and 
finite time blow-up appears. 

\begin{figure}[ht]
   \centering
   \includegraphics[width=0.7\textwidth]{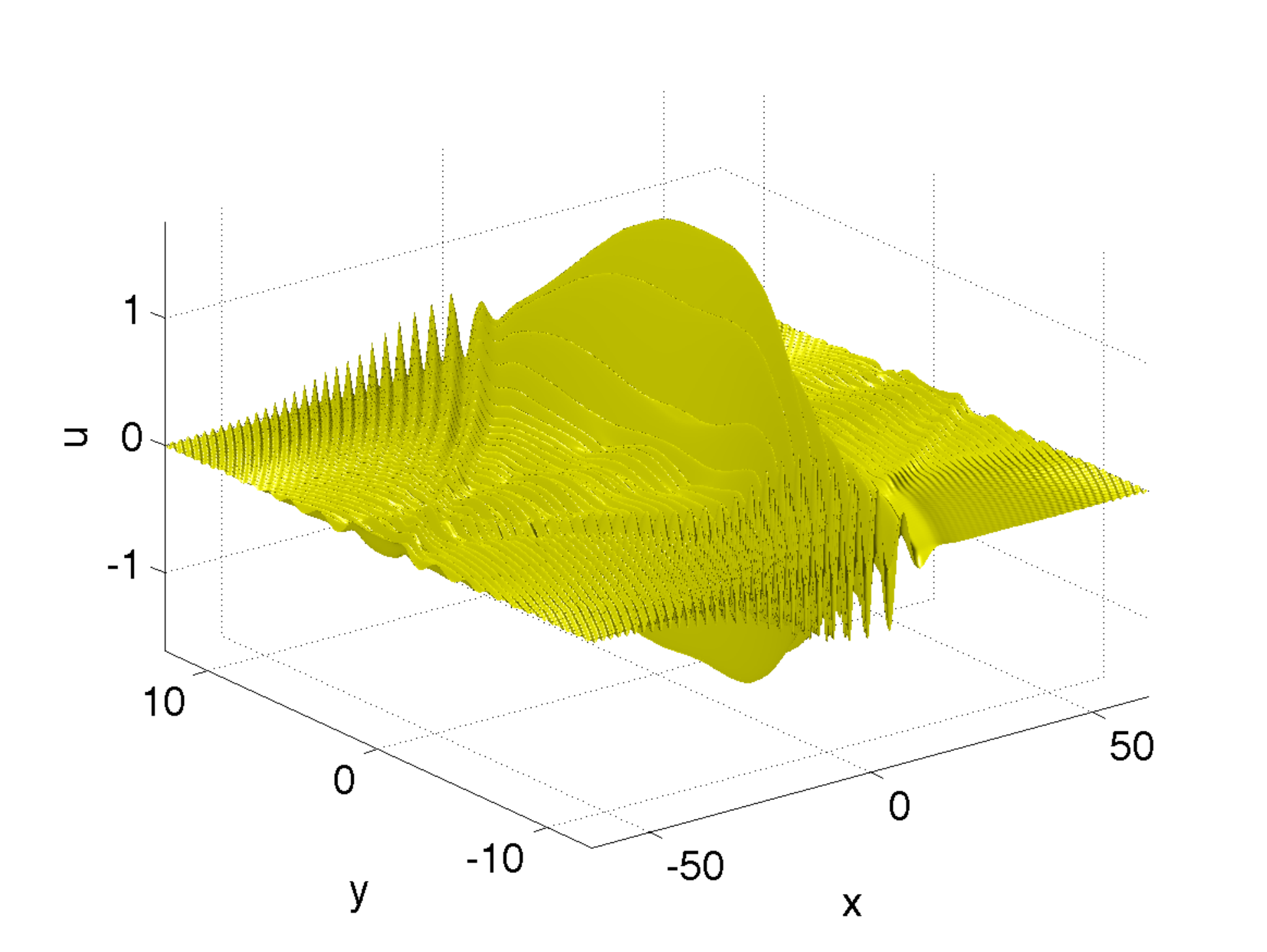}
   \caption{Solution to the gKP II equation (\ref{gKP}) with $\lambda=1$, 
   $n=4/3$ for the 
   initial data $u_{0}=6\,\partial_{xx}\exp(-(x^{2}+y^{2}))$ at 
   $t=2$.}
   \label{gKPIIn436gausst2}
\end{figure}

There are no theoretical results on blow-up for gKP II solutions. In 
\cite{KS12}, numerical results indicate that there is no blow-up in 
this context for $n\leq 2$. We confirm this here by looking at solutions for $n=4/3$ and
$u_{0}=6\,\partial_{xx}\exp(-(x^{2}+y^{2}))$ with $L_{x}=20$, $L_{y}=4$,
$N_{x}=2^{11}$, $N_{y}=2^{10}$ and $N_{t}=10^{3}$ for $t\leq 2$. We 
obtain a $\Delta\sim 10^{-3}$. As can be seen in 
Fig.~\ref{gKPIIn436gausst2}, there is no indication of blow-up, the 
initial data seem to be simply radiated away. Note that the radiation 
propagates in this case to the right, whereas the algebraic tails 
stretch to $-\infty$. Since we compute for a 
long time, the periodicity of the considered situation leads to  a 
reentering of the oscillations from the left side.

\begin{figure}[ht]
   \centering
   \includegraphics[width=0.49\textwidth]{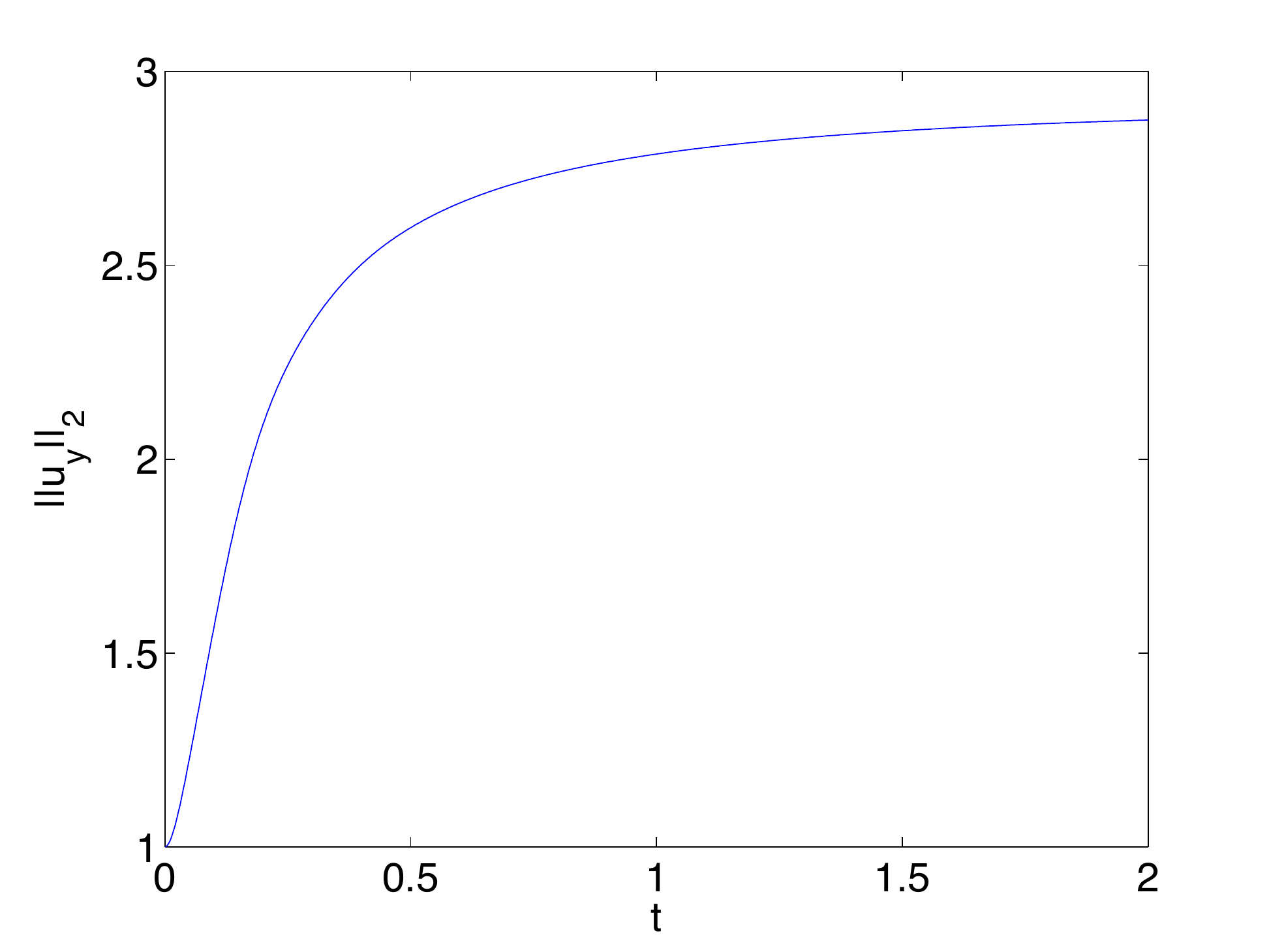}
   \includegraphics[width=0.49\textwidth]{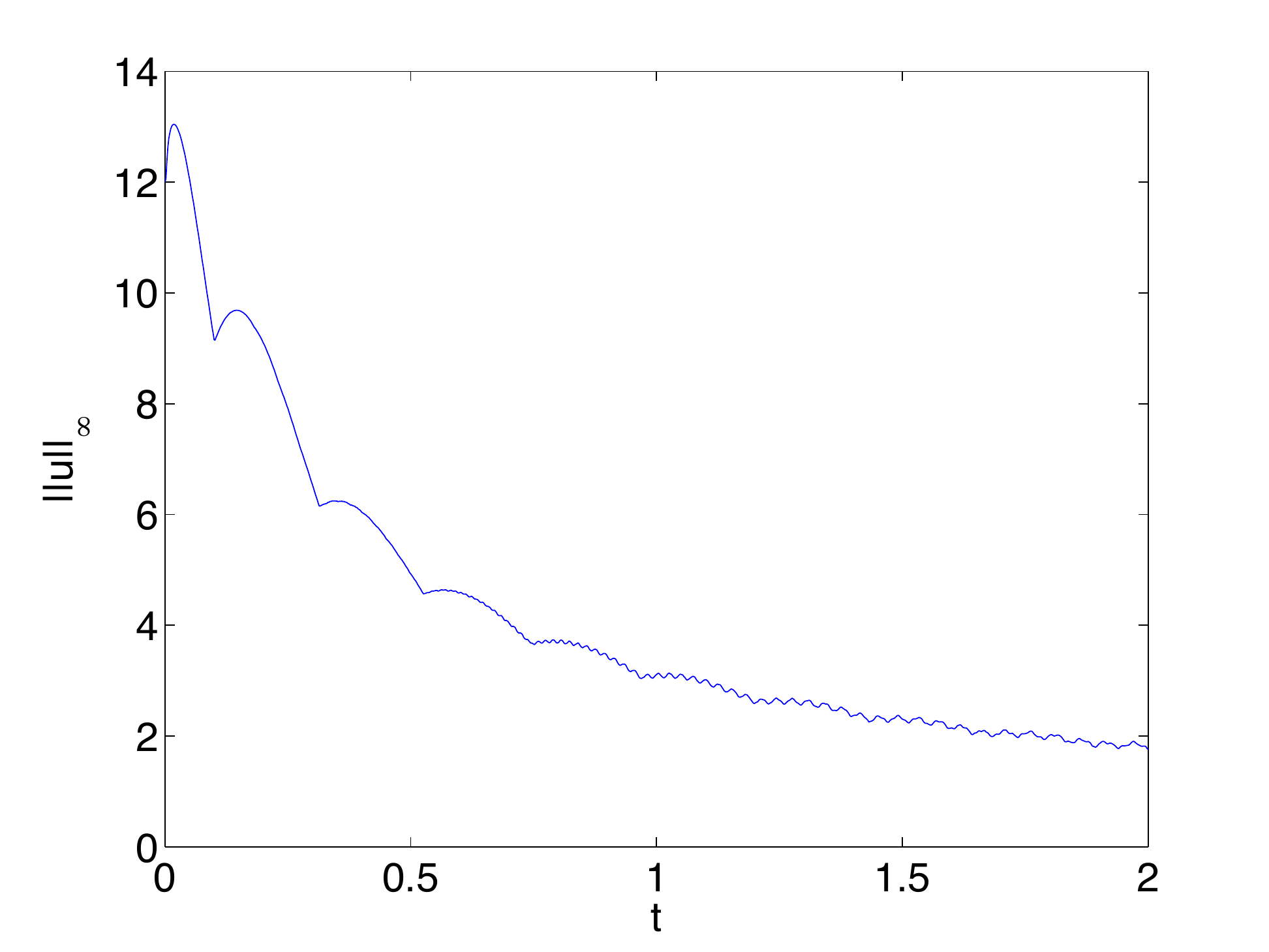}
   \caption{Solution to the equation (\ref{gKP}) with $\lambda=1$, 
   $n=4/3$ for the 
   initial data $u_{0}=6\,\partial_{xx}\exp(-(x^{2}+y^{2}))$; the 
   $L_{2}$ norm of $u_{y}$ on the left, the
   $L_{\infty}$ norm of $u$ on right.}
   \label{gKPIIn436gaussnorm}
\end{figure}

The norms for this solution in Fig.~\ref{gKPIIn436gaussnorm} indicate 
that the $L_{\infty}$ norm decreases monotonically, whereas the 
$L_{2}$ norm of $u_{y}$ appears to reach a constant value, at least 
on the considered time scales. The latter seems to be due to Gibbs phenomena at 
the boundaries because of the algebraic decrease of the solution towards 
infinity. The decrease of the $L_{\infty}$ norm of the solution shows 
in any case a completely different behavior from the gKP I blow-up 
scenarios. Note that the solutions behave qualitatively the same for 
twice the initial data in (\ref{gKPIIn436gausst2}). Thus we did not 
find an indication for blow-up in gKP II solutions for $n=4/3$. It 
seems that the defocusing character stabilizes the solution against 
blow-up, whereas the opposite is true for gKP I, which has a focusing 
effect. 

\section{The supercritical case $n = 2$}
In this section we study the supercritical case $n=2$ which 
appears to be typical for $n>4/3$ for gKP I solutions. It is also relevant in applications 
since it appears in the modelling of sound waves in antiferromagnetic 
systems \cite{FT85}. We will again consider the initial data 
(\ref{initial}) for various values of the parameter $\beta$. 

\begin{figure}[ht]
   \centering
   \includegraphics[width=0.7\textwidth]{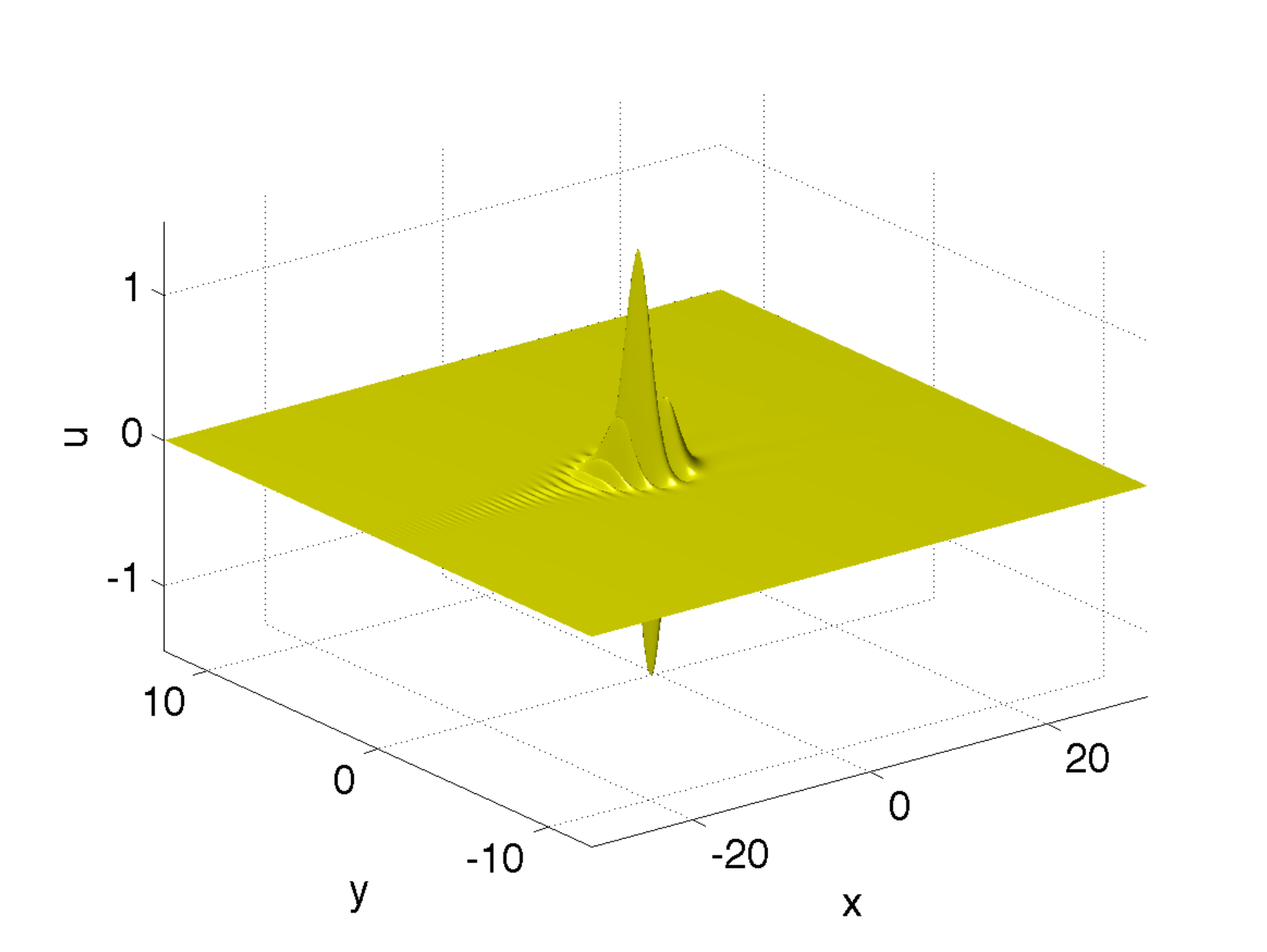}
   \caption{Solution to the gKP I equation (\ref{gKP}) with 
   $n=2$ for the 
   initial data $u_{0}=\partial_{xx}\exp(-(x^{2}+y^{2}))$ at 
   $t=0.1$.}
   \label{gKPn2gausst01}
\end{figure}

First we study gKP I (\ref{gKP}) with $n=2$ and the initial data 
$u_{0}=\partial_{xx}\exp(-(x^{2}+y^{2}))$ for which the energy is 
positive. We use $L_{x}=10$, $L_{y}=4$,
$N_{x}=N_{y}=2^{10}$ and $N_{t}=10^{4}$ for $t\leq 0.1$ and obtain a 
relative conservation of the computed energy of the order of $10^{-8}$.
It can be seen in Fig.~\ref{gKPn2gausst01} that the initial data appear
to be simply radiated away. 

\begin{figure}[ht]
   \centering
   \includegraphics[width=0.49\textwidth]{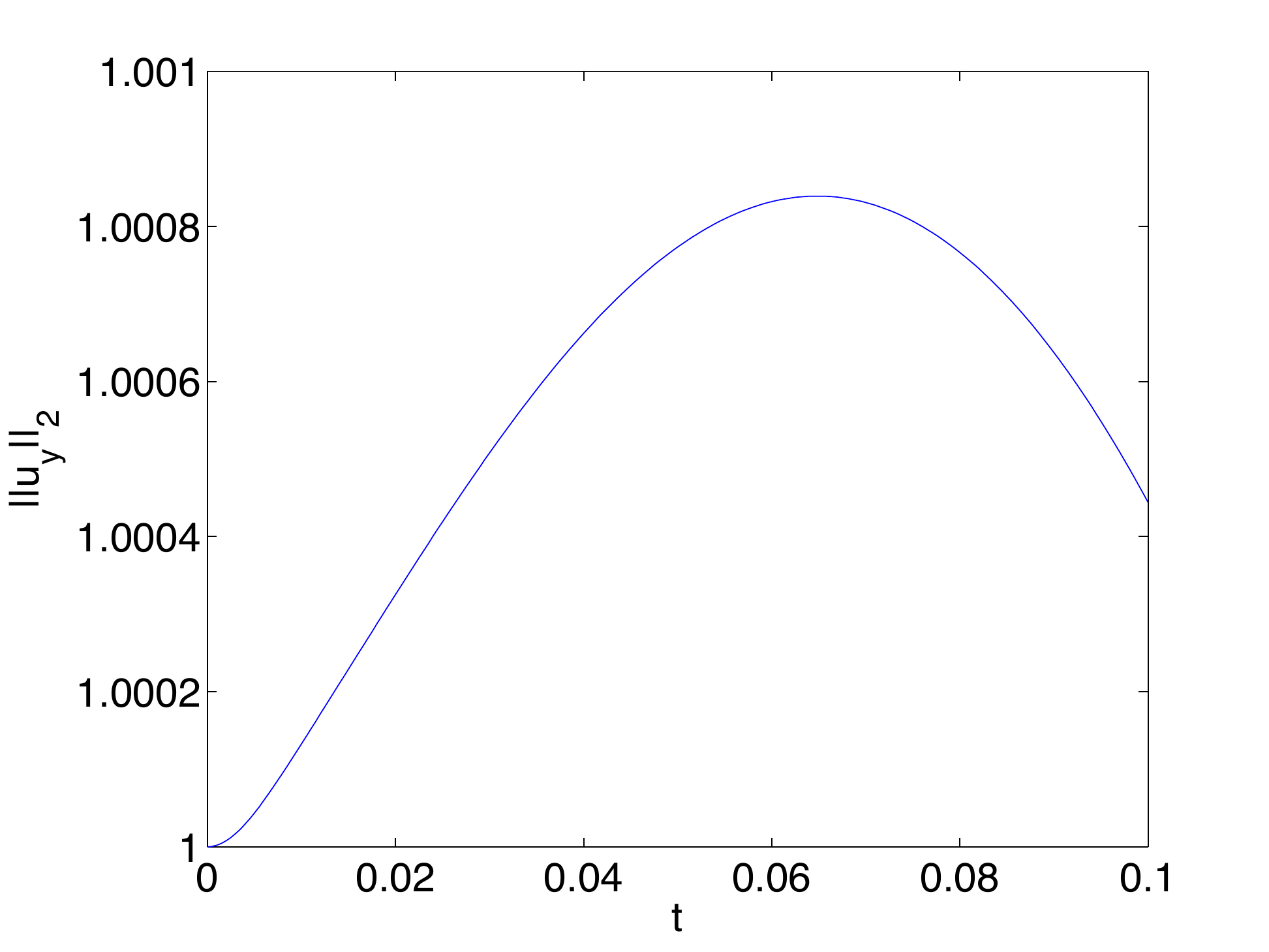}
   \includegraphics[width=0.49\textwidth]{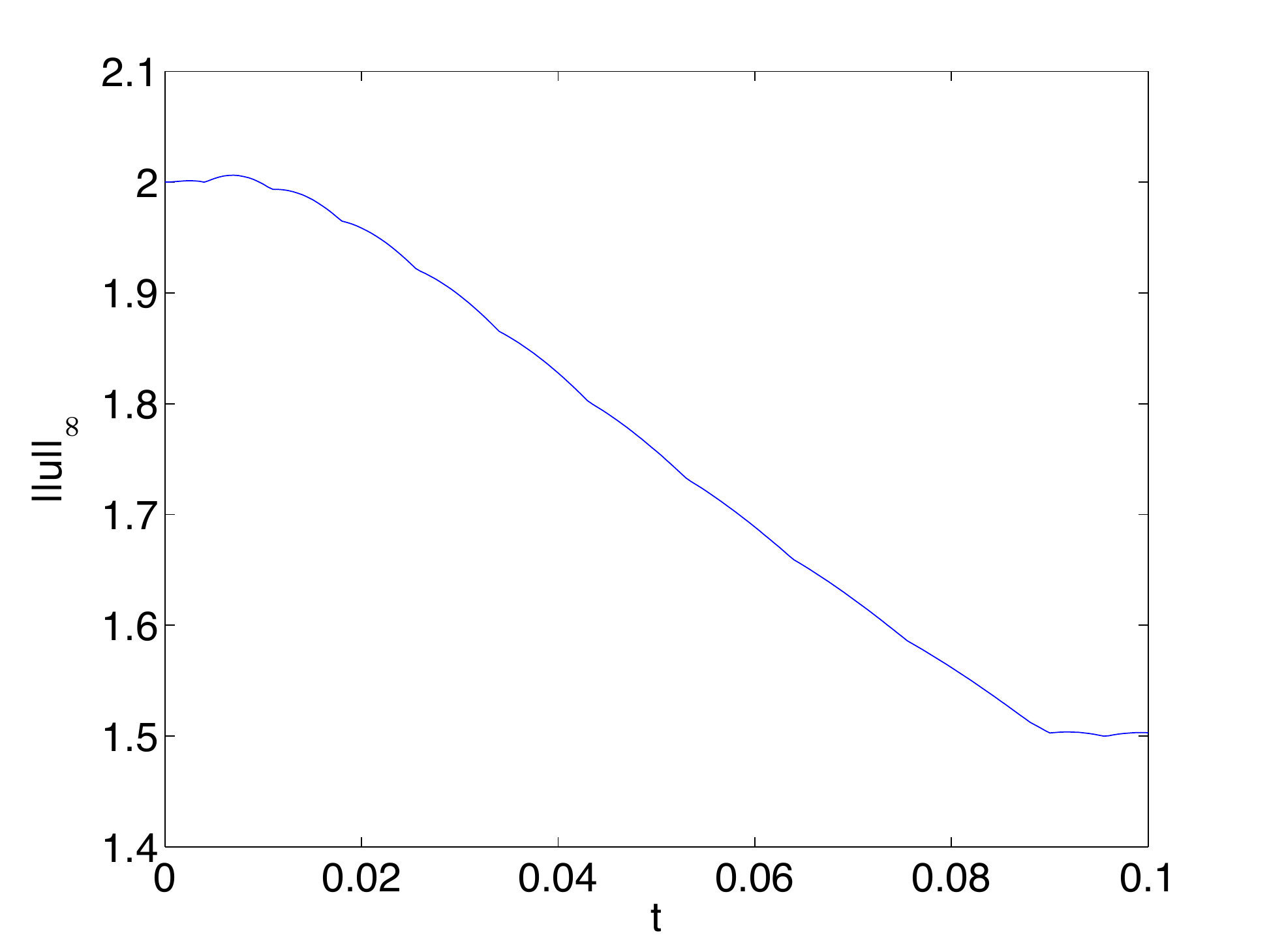}
   \caption{Norms of the solution to the gKP I equation (\ref{gKP}) with 
   $n=2$ for the 
   initial data $u_{0}=\partial_{xx}\exp(-(x^{2}+y^{2}))$; the 
   $L_{2}$ norm of $u_{y}$ on the left, the
   $L_{\infty}$ norm of $u$ on right.}
   \label{gKPn2gaussnorm}
\end{figure}

This is confirmed by the norms of the solution shown in 
Fig.~\ref{gKPn2gaussnorm}. Both the $L_{2}$ norm of $u_{y}$ and the 
$L_{\infty}$ norm of $u$ appear to decrease monotonically for large 
$t$. There is no indication of blow-up in this case. 

\begin{figure}[ht]
   \centering
   \includegraphics[width=\textwidth]{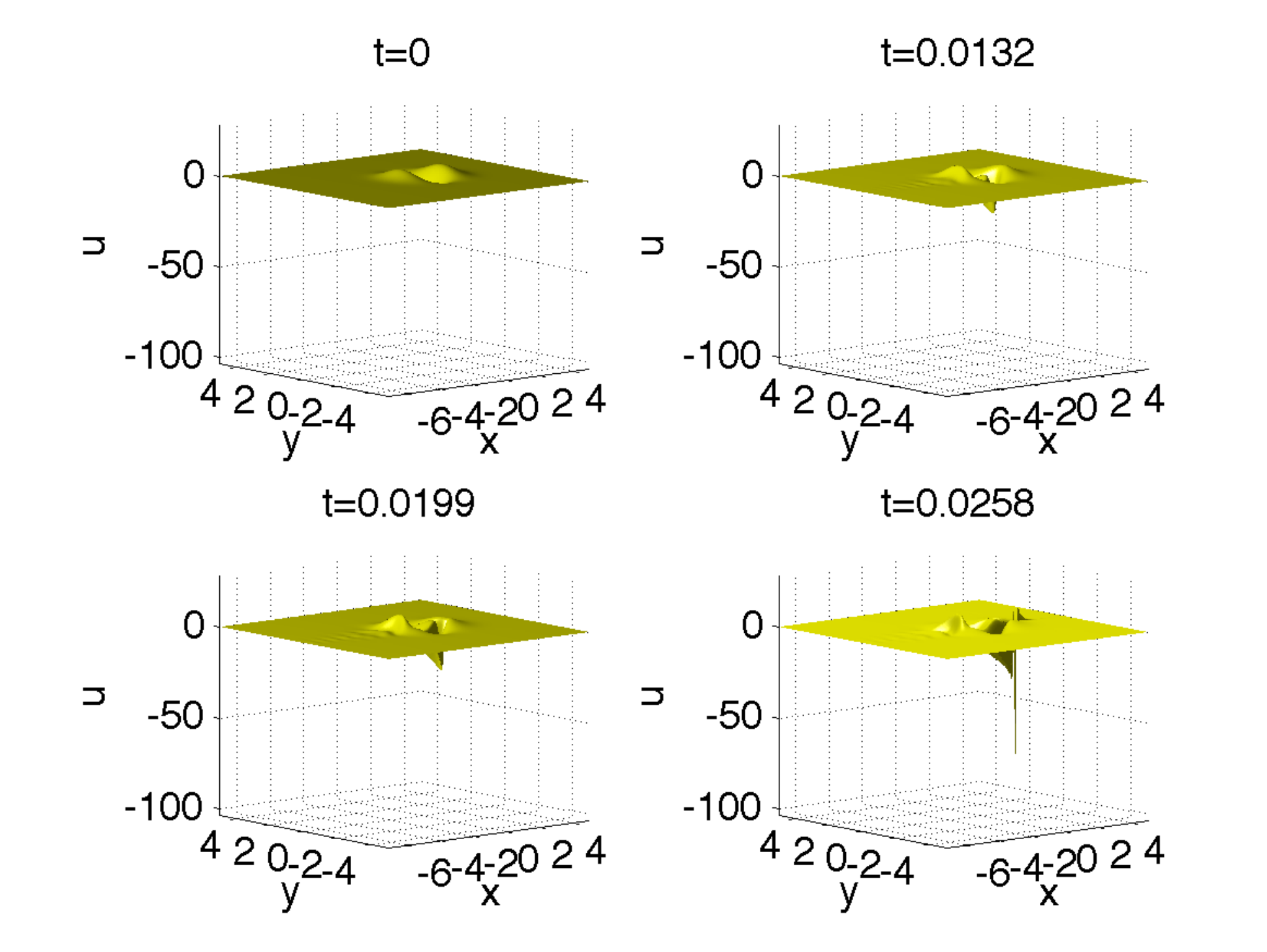}
   \caption{Solution to the gKP I equation (\ref{gKP}) with 
   $n=2$ for the 
   initial data $u_{0}=6\,\partial_{xx}\exp(-(x^{2}+y^{2}))$ for 
   several times. }\label{fig:n2_u}
\end{figure}

The situation is completely different for initial data  $u_0(x,y) = 6\,\del_{xx} \exp\bigl(-(x^2 + y^2)\bigr)
$ for which the energy is  negative. The computation is carried out 
with  $L_x = L_y = 5$, $N_x = 2^{11}$,
 $N_y = 2^{13}$ and           $N_t = 50000$ time steps for 
 $t\in[0,0.0265]$. As can be seen in Fig.~\ref{fig:n2_u}, the 
 solution develops the same dispersive oscillations propagating to 
 the left as in Fig.~\ref{gKPn2gausst01} (only part of the 
 computational domain is shown), but this cannot stop the 
 negative peak from becoming more and more compressed. This peak finally 
 appears to blow up in a point. The code is stopped at the time $t = 
 0.0258375$ when $\Delta>10^{-3}$. Note that the final stage before 
 blow-up is happening on much shorter time scales than in the 
 previous section for $n=4/3$. 

\begin{figure}[ht]
   \centering
   \includegraphics[width=0.49\textwidth]{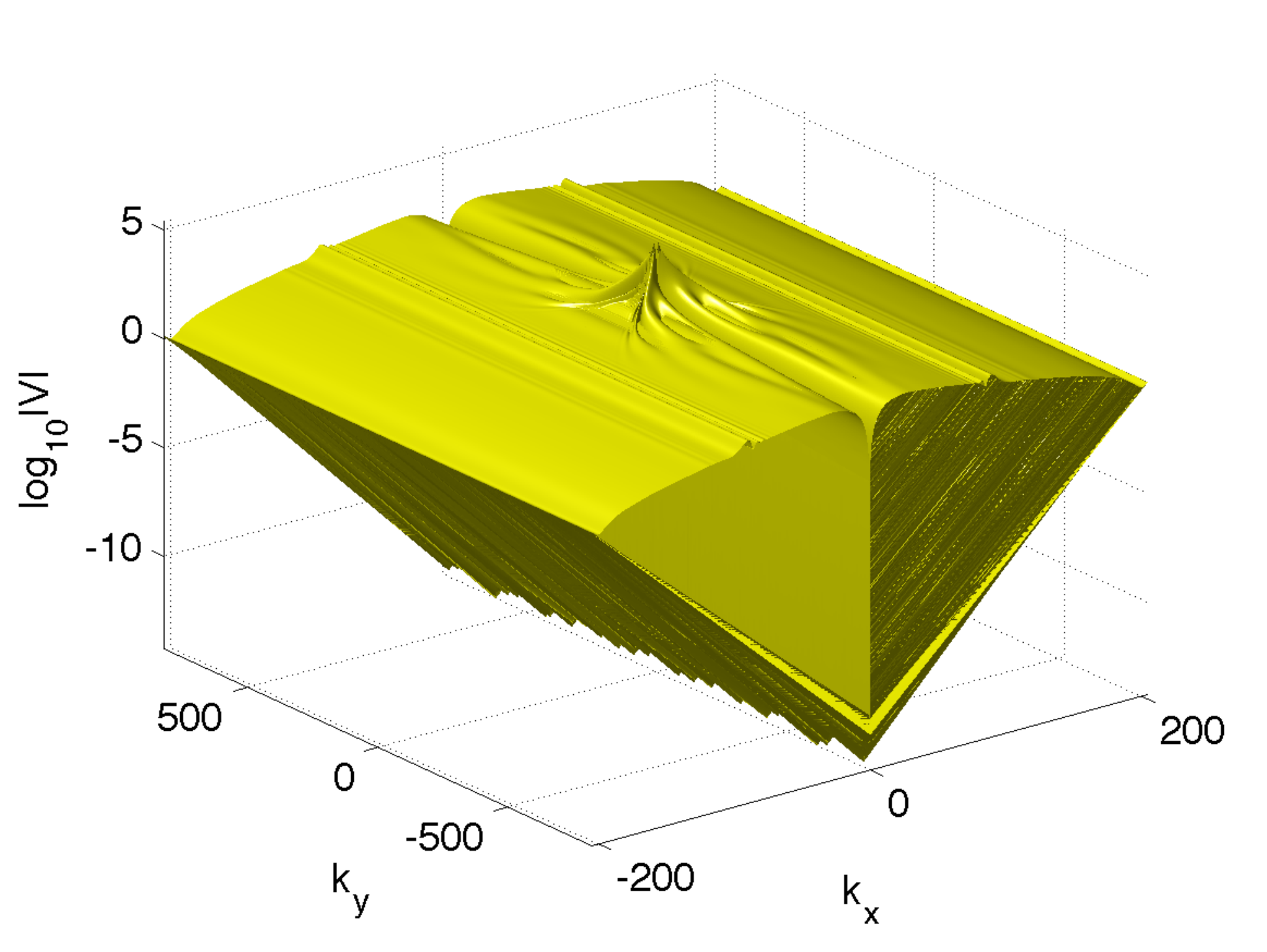}
   \caption{Modulus of the Fourier coefficients of the solution 
   $u$ in Figure \ref{fig:n2_u} at $t = 0.0258$.}\label{fig:n2_u1}
\end{figure}

It can be seen in Fig.~\ref{fig:n2_u1} that we run out of resolution 
in Fourier space. As expected from the rescaling (\ref{gKP4}), there 
are strong gradients especially in $y$ which cannot be addressed with even 
higher resolution. 

\begin{figure}[ht]
   \centering
   \includegraphics[width=0.49\textwidth]{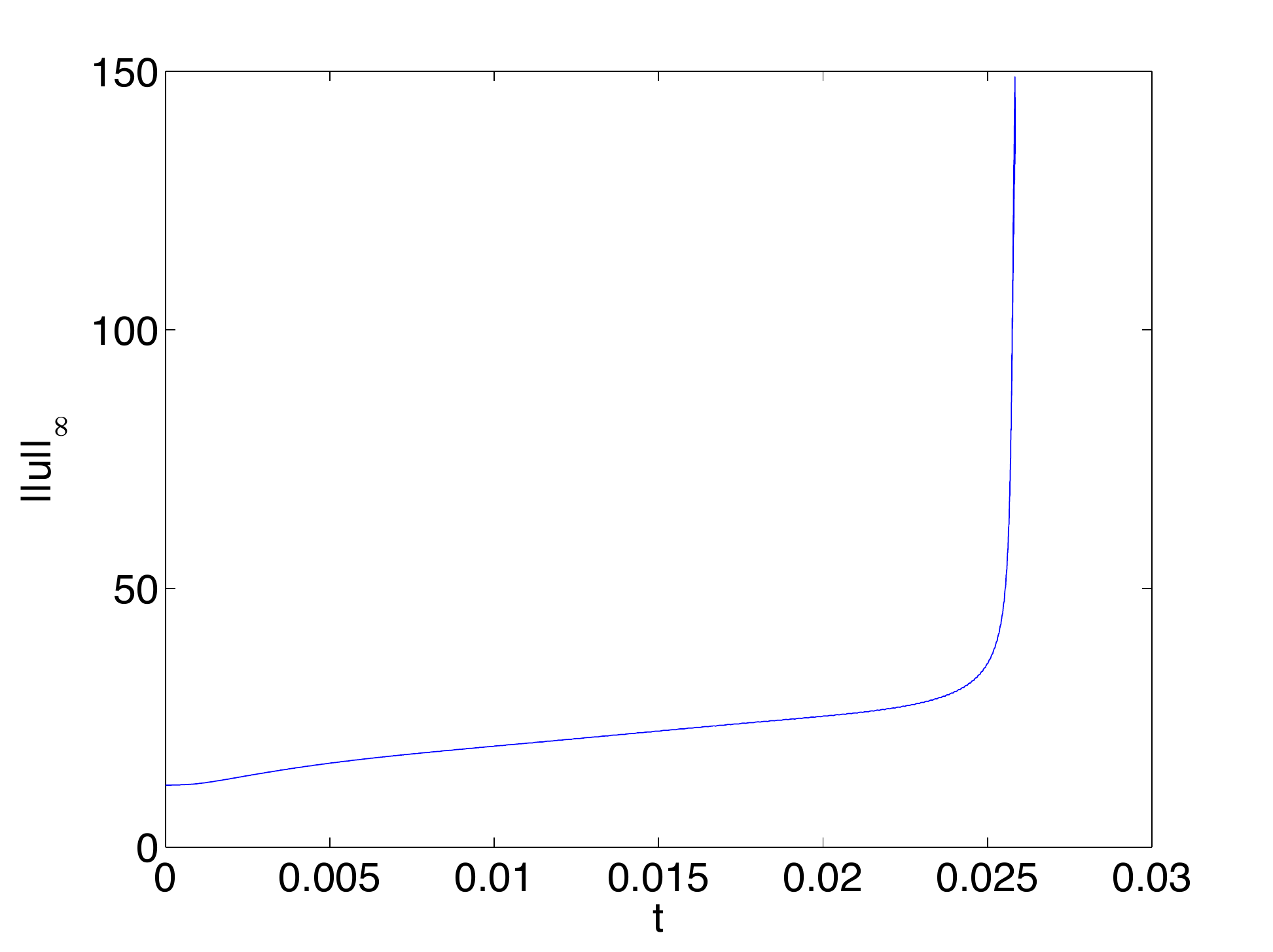}
   \includegraphics[width=0.49\textwidth]{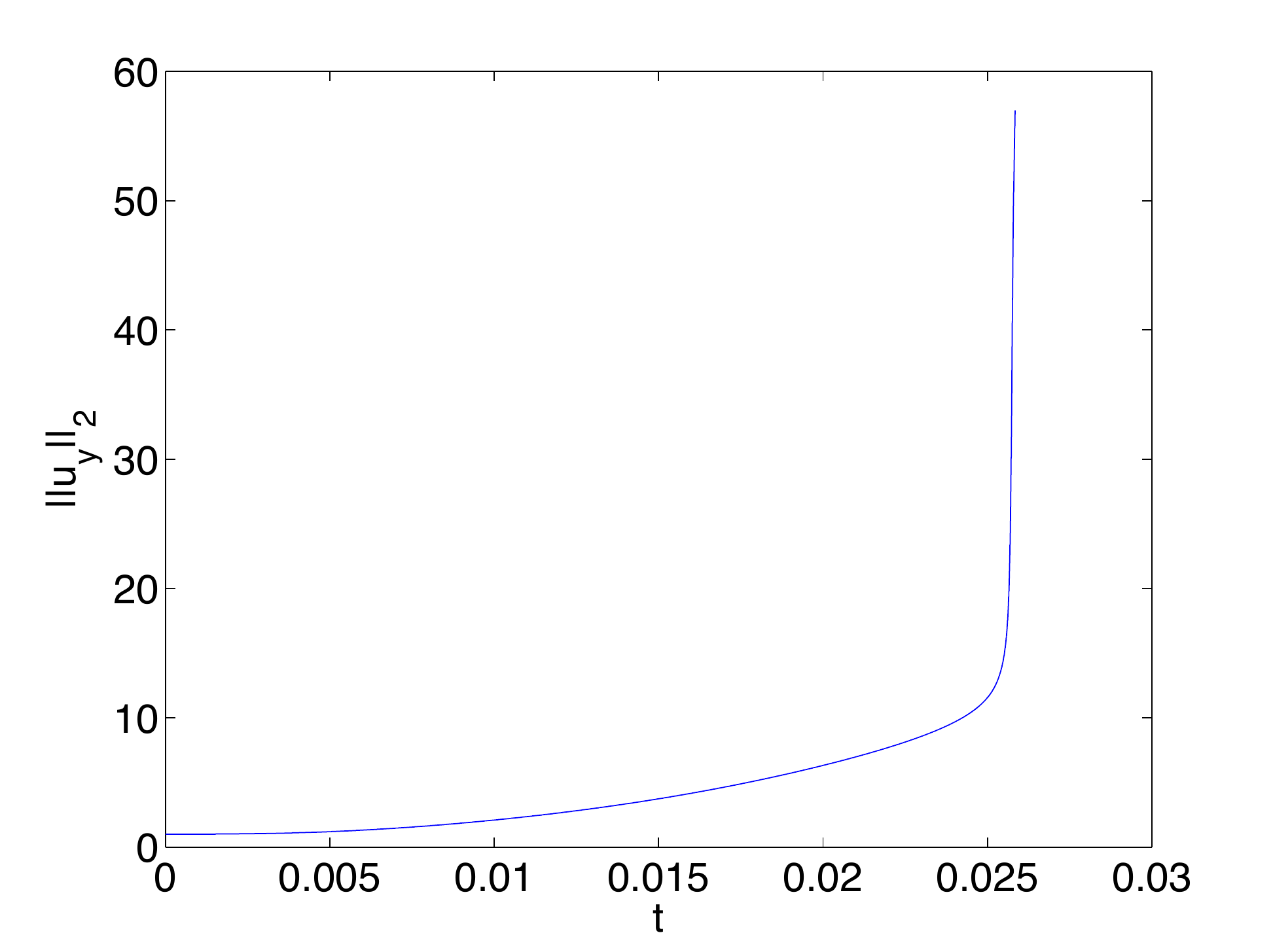}
   \caption{Norms of the solution to the gKP I equation (\ref{gKP}) with 
   $n=2$ for the 
   initial data $u_{0}=6\,\partial_{xx}\exp(-(x^{2}+y^{2}))$) in 
   dependence of time; on the left the $L_{\infty}$ norm  of the 
   solution $u$, on the right the $L_{2}$ norm of $u_{y}$ normalized 
   to 1 at $t=0$.}\label{fig:n2_LinfNorm_u}
\end{figure}

Both the $L_{2}$ norm of $u_{y}$ and the $L_{\infty}$ norm of $u$ 
indicate a blow-up as is clear from Fig.~\ref{fig:n2_LinfNorm_u}. 

\begin{figure}[ht]
   \centering
   \includegraphics[width=0.49\textwidth]{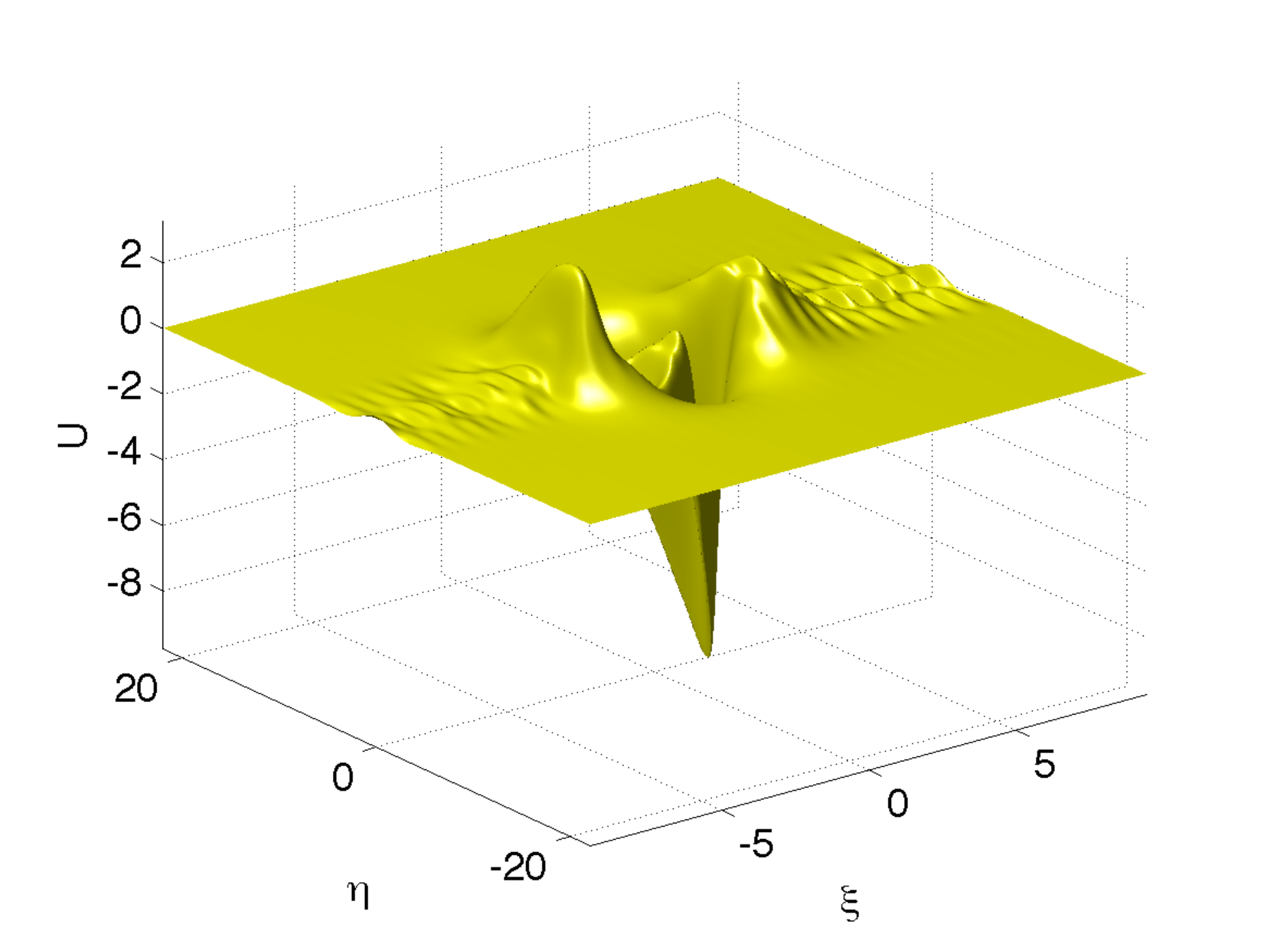}
   \includegraphics[width=0.49\textwidth]{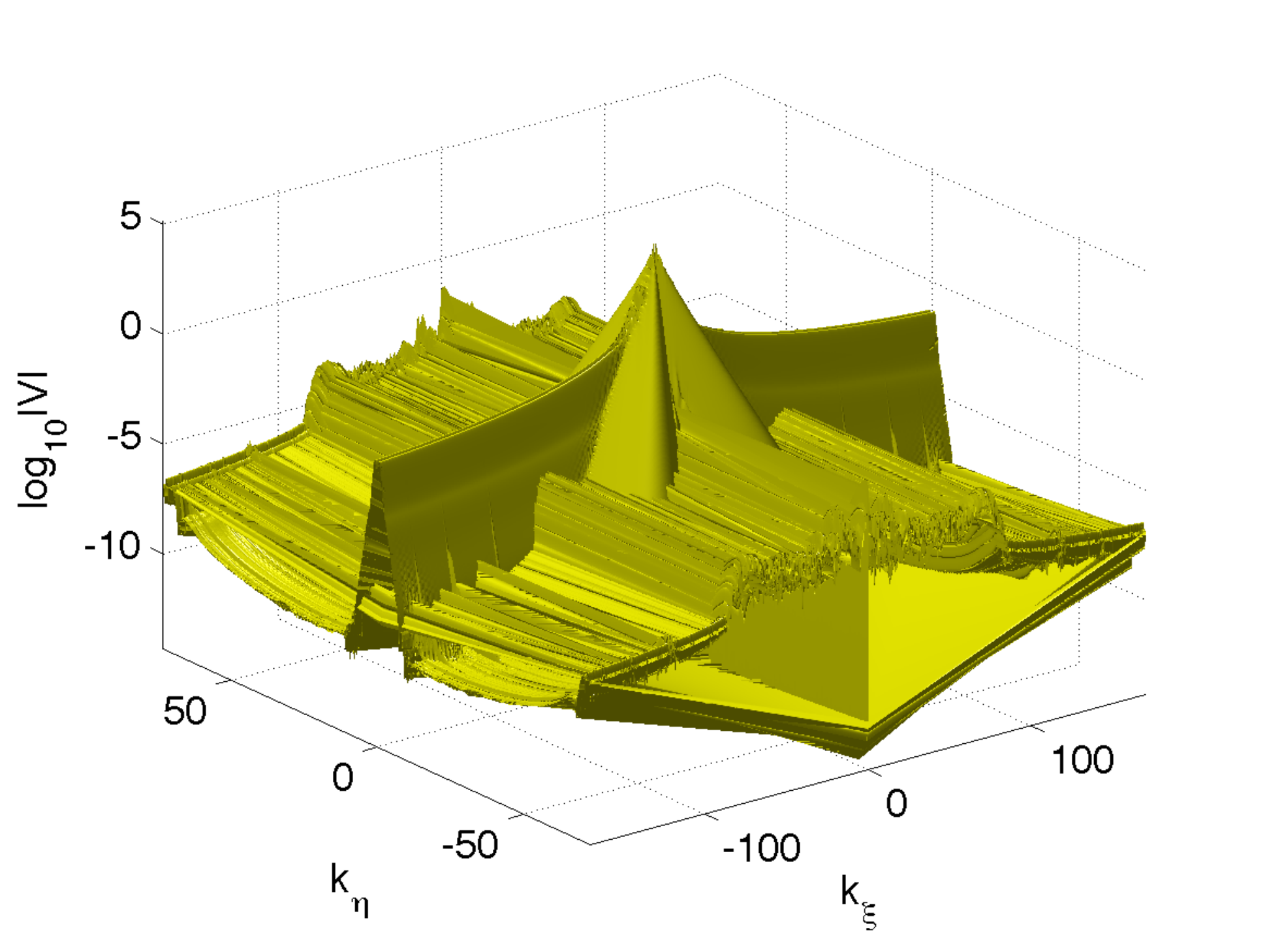}
   \caption{Solution to the rescaled gKP I equation (\ref{gKP5}) with $n=2$ for the 
   initial data $U_{0}=6\,\partial_{\xi\xi}\exp(-(\xi^{2}+\eta^{2}))$ at 
   $\tau=0.5$ on the left, and the modulus of the corresponding 
   Fourier coefficients on the right.}
   \label{gKPn26gausstau05}
\end{figure}

In an attempt to identify the type of blow-up, we  solve the 
rescaled gKP equation (\ref{gKP5}) for the initial data  
$U_{0}=6\,\partial_{\xi\xi}\exp(-(\xi^{2}+\eta^{2}))$. The computation is carried out 
with $L_{\xi}=3$, $L_{\eta}=7$,
$N_{\xi}=N_{\eta}=2^{10}$ and $N_{\tau}=10^{4}$ for $\tau\leq 0.5$ resulting 
in a relative conservation of the computed mass of the order of 
$10^{-3}$ (again energy conservation is much worse due to Gibbs 
phenomena at the boundary of the computational domain). The solution at the final time is shown in 
Fig.~\ref{gKPn26gausstau05}. It can be seen that the `zooming in' 
effect of the rescaling close to blow-up is at work. But due to the 
dispersive oscillations with a slowly decaying amplitude and due to 
the algebraic fall off of the solution, this eventually leads to a 
Gibbs phenomenon.
This is clearly reflected in the Fourier coefficients of the solution 
in the same figure. The instabilities due to the 
increase in the high wave numbers will eventually break the code. 
This problem does not disappear on a larger computational domain, the 
code tends to become unreliable even earlier due the above described 
problems. 

\begin{figure}[ht]
   \centering
   \includegraphics[width=0.49\textwidth]{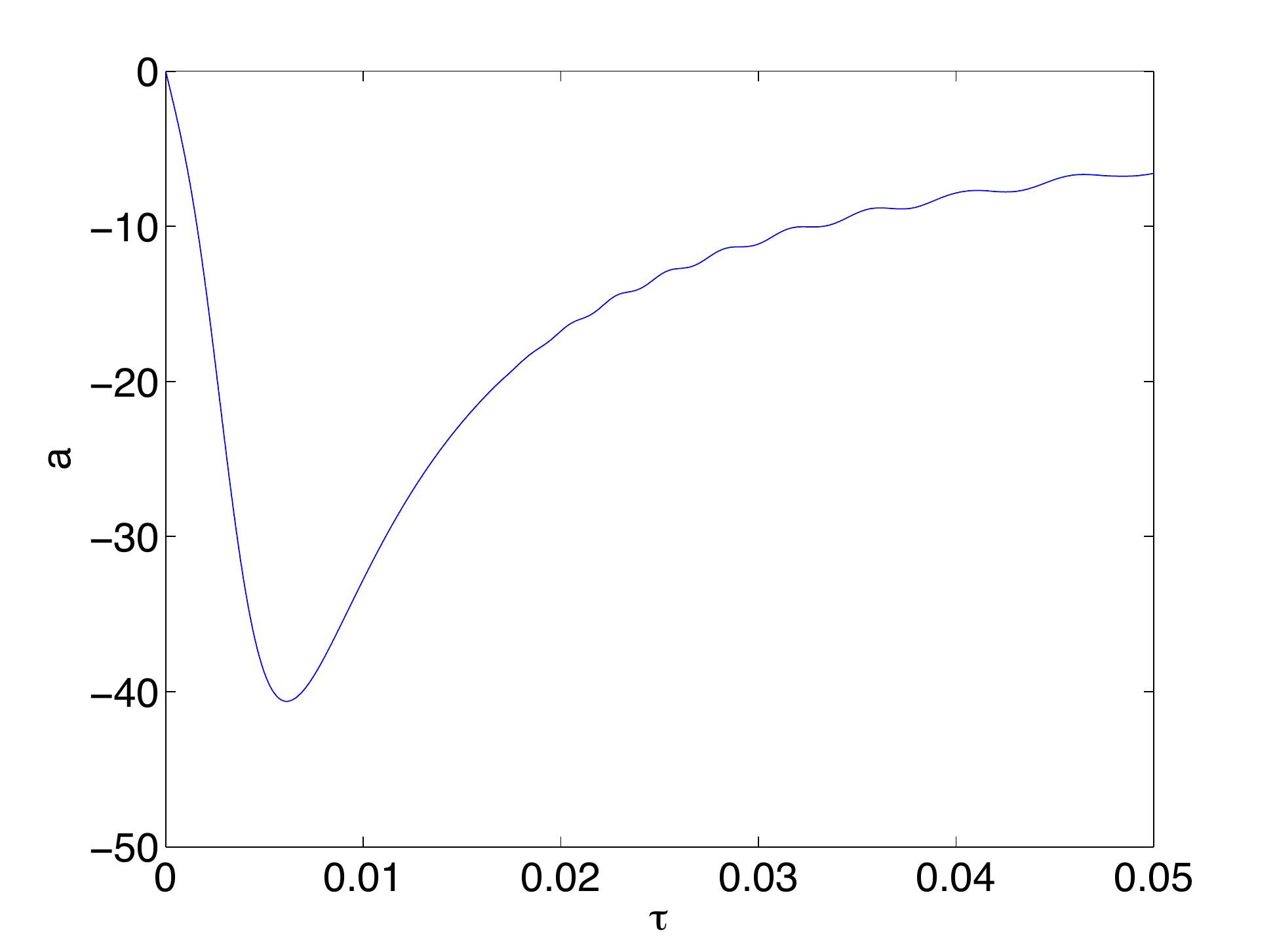}
   \includegraphics[width=0.49\textwidth]{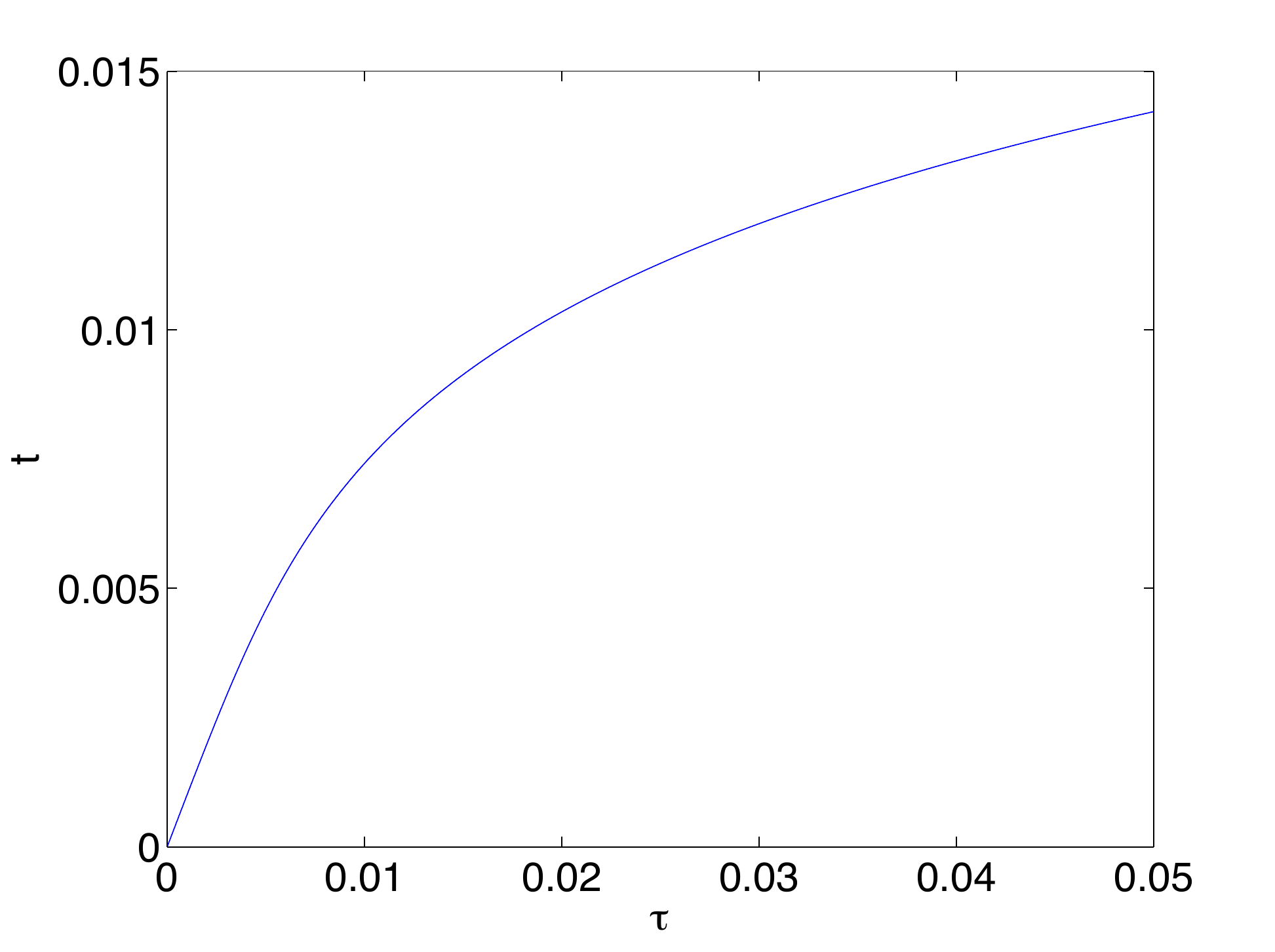}
   \caption{Solution to the equation (\ref{gKP5}) with $n=2$ for the 
   initial data $U_{0}=6\,\partial_{\xi\xi}\exp(-(\xi^{2}+\eta^{2}))$; the quantity 
   $a$ (\ref{a}) in dependence of $\tau$ on the left, the physical time on the 
   right.}
   \label{gKPn26gaussta}
\end{figure}

In Fig.~\ref{gKPn26gaussta} we show the logarithmic derivative $a$ in 
(\ref{a}) of the scaling factor $L$. It decreases rapidly to strongly 
negative values which implies that most of the dynamical evolution 
happens for comparatively small $\tau$. The quantity $a$ appears to 
approach a negative constant which would indicate an exponential  
$\tau$-dependence of $L$. The oscillations in $a$ are due to the 
imposed periodicity forcing the dispersive radiation to reenter the 
computational domain on the right side. Still we do not get close 
enough to the blow-up to be able to decide which type of blow-up is 
realized here. This is also clear from the physical time we show in 
the same figure which is not close enough to the time in 
Fig.~\ref{fig:n2_u} where the direct integration breaks. 

\begin{figure}[ht]
   \centering
   \includegraphics[width=0.49\textwidth]{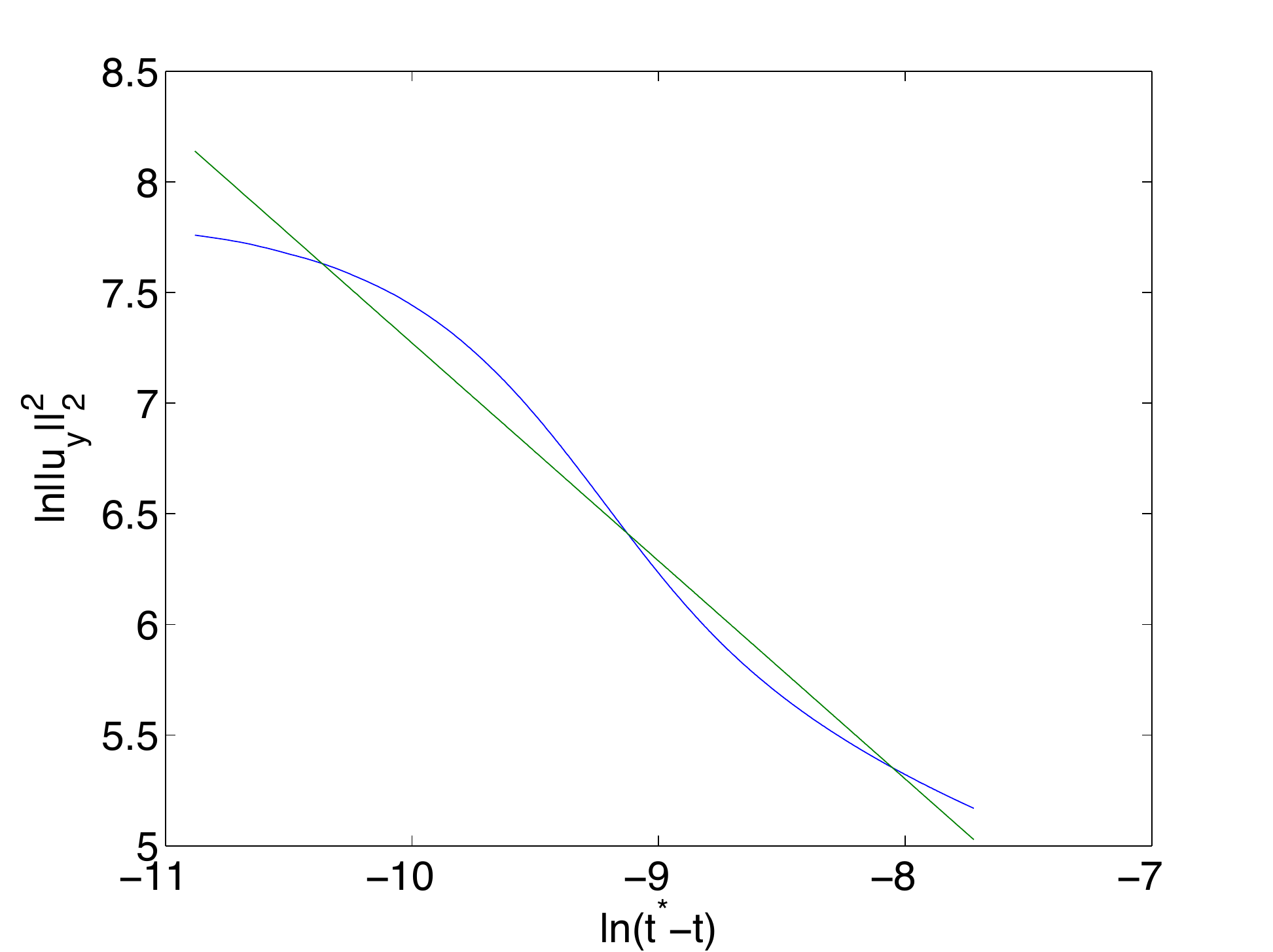}
   \includegraphics[width=0.49\textwidth]{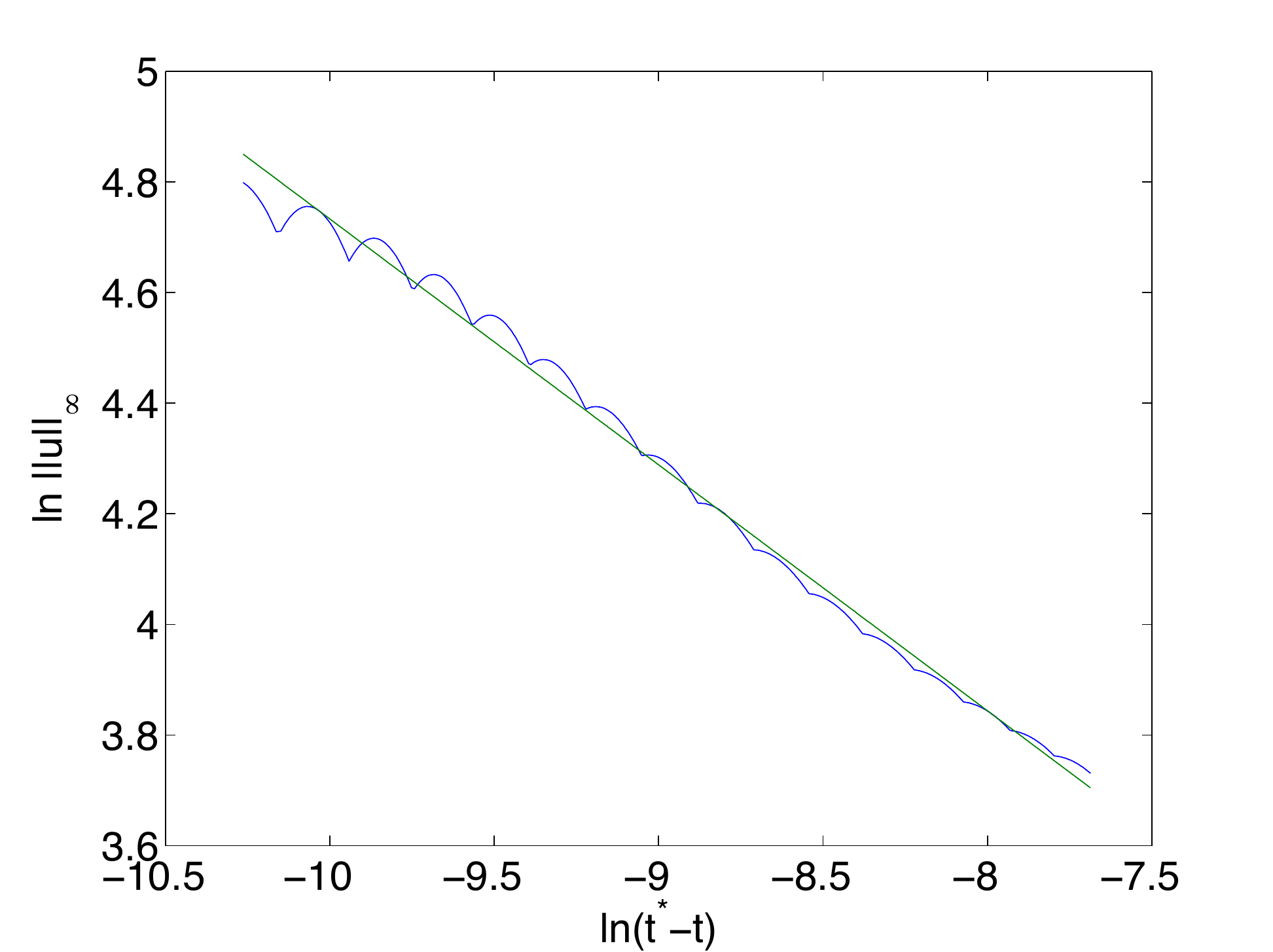}
   \caption{Fitting of $\ln||u_{y}||_{2}^{2}$ (left) and 
 $\ln||u||_{\infty}$ (right)  in blue to $c \ln(t^{*}-t)+C$ in green 
 for the situation shown in Fig.~\ref{fig:n2_u}.}\label{fig:Linf_u_fit}
\end{figure}

Due to the described instabilities of the code we did not succeed to 
run it with higher resolutions on larger domains even with a very 
high time resolution. Thus it appears again more promising to 
trace the norms of the solution obtained via a direct integration of 
the gKP I equation as in Fig.~\ref{fig:n2_u}. The norms indicate as Fig.~\ref{gKPn26gaussta} an 
exponential $\tau$ dependence of $L$ which implies with (\ref{gKP4}) 
and (\ref{L2yKP})
\begin{equation}
    ||u_{y}||_{2}\propto (t^{*}-t)^{-(1+4/n)/6},\quad ||u||_{\infty}\propto 
    (t^{*}-t)^{-2/(3n)}
    \label{n2scal}.
\end{equation}
As before we fit $\ln ||u||_{\infty}$ for the solution shown in 
Fig.~\ref{fig:n2_u} to $\ln C_1+c_{1}\ln(t^* - t)$ and 
$\ln||u_{y}||_{2}^{2}$ to $ \ln C_2+c_{2}\ln(t^* - t)$. 
Since strong gradients appear 
close to blow-up, the $L_{\infty}$ norm of the solution 
 appears here to be 
asymptotically a more reliable indicator. The fitting is performed for 
times close to the last recorded time since the asymptotic behavior 
can be only expected for $t\sim t^{*}$. To have enough points for a 
reliable fitting, the time interval is chosen to include the 800 last 
computed points up to $t=0.0257845$.  We get for              
$||u_{y}||_{2}^{2}$ the fitting parameters $ C_2 =   -2.5795$, 
$c_2 =  -0.9851$ and  $t^* =  0.0258$. Thus 
$c_{2}$ is compatible with the theoretically expected $-1$. 
The quality of the fitting can be seen in Fig.~\ref{fig:Linf_u_fit}.              
For $||u||_{\infty}$ we get $ C_1 =   0.2878$, $c_1 = -0.4445$ and  $t^* =   0.0258$. Despite the oscillations in 
$||u||_{\infty}$ in Fig.~\ref{fig:Linf_u_fit}, the value of 
$c_{1}$ is very close to the expected $-1/3$, and the value for 
$t^{*}$ coincides within numerical precision with the value found 
from $||u_{y}||_{2}$ which shows the consistency of the approach.



\begin{figure}[ht]
   \centering
   \includegraphics[width=0.49\textwidth]{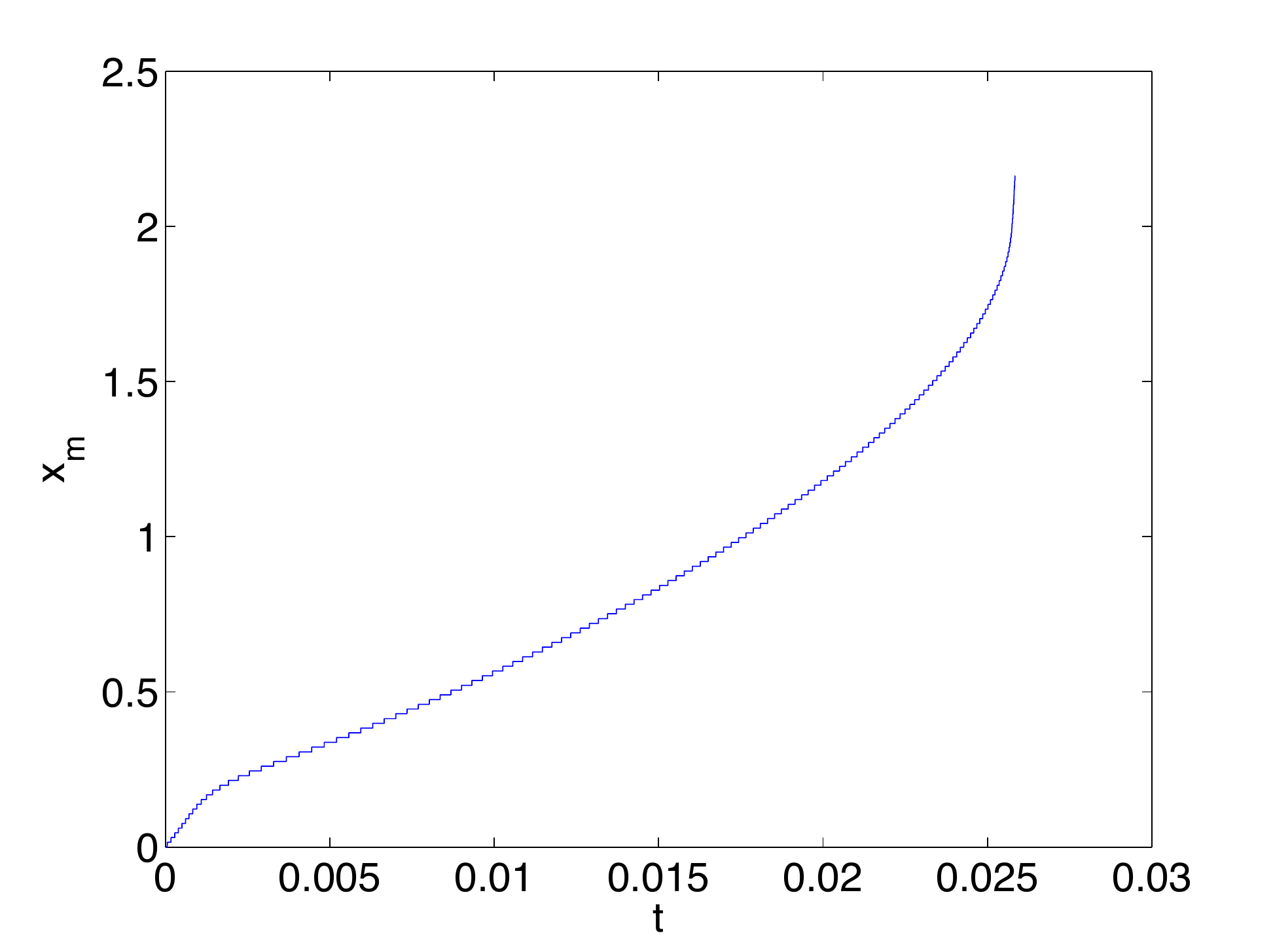}
   \includegraphics[width=0.49\textwidth]{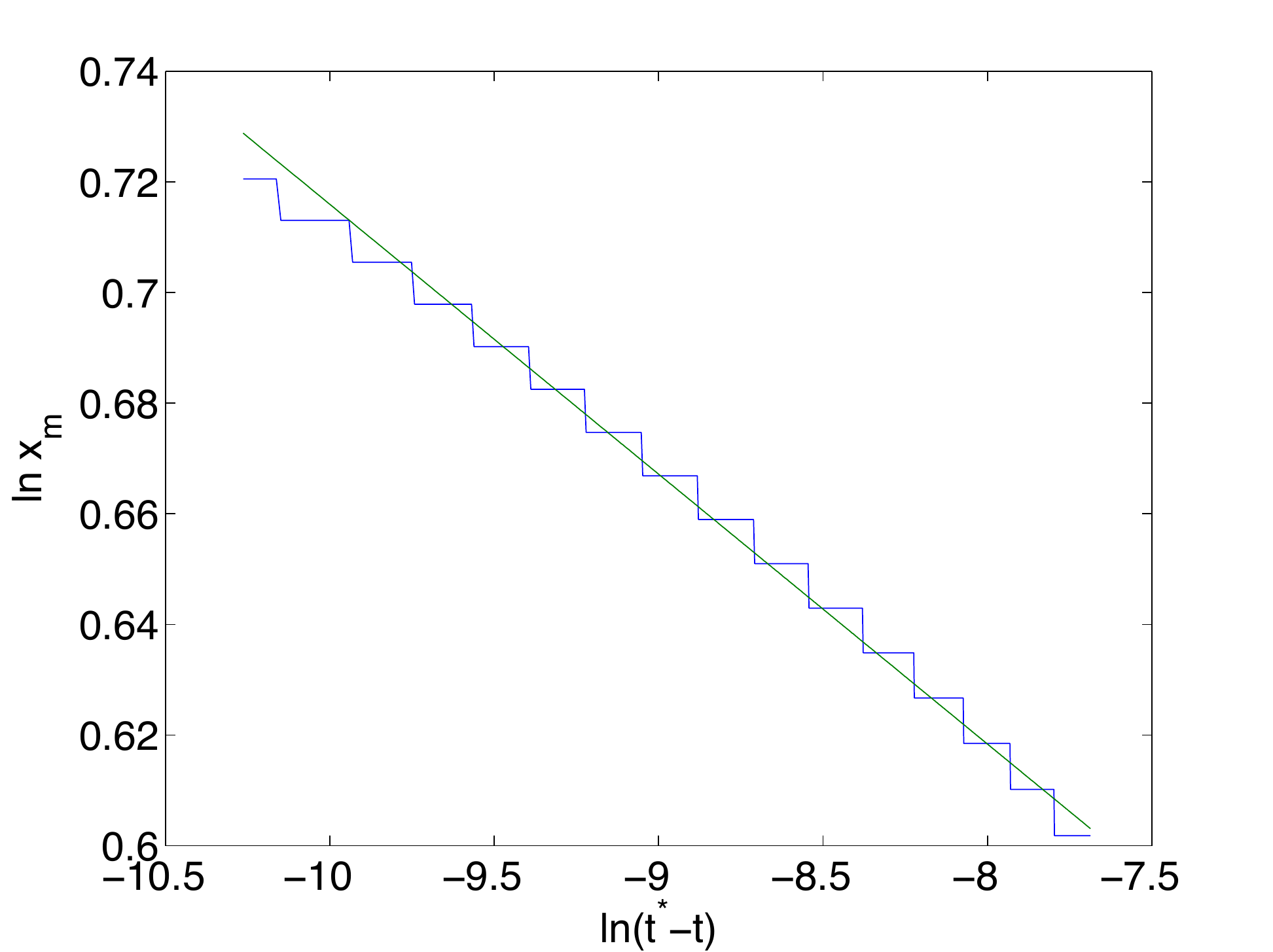}
   \caption{The location 
   $x_{m}$ of the global minimum of the 
   solution   
   $u$ in Figure \ref{fig:n2_u} in dependence of 
   $t$ on the left, and a fit of $\ln x_{m}$ in blue for $t\sim 
   t^{*}$ to 
   $\alpha_{1}\ln(t^{*}-t)+\alpha_{2}$ in green on the right. }
   \label{gKPn26gaussxm}
\end{figure}

The location of the global minimum can be seen as a function of time in 
Fig.~\ref{gKPn26gaussxm}. The behavior is very similar to the case 
$n=4/3$ in Fig.~\ref{gKPn4312gaussxm}. Again it is difficult to decide whether there will be a blow-up   
at finite $x^{*}$, whereas again $y_{m}=0$ for symmetry 
reasons. If we fit $\ln 
x_{m}$ as the norms of $u$ above for the last 800 time steps to 
$\alpha_{1}\ln (t^{*}-t)+\alpha_{2}$, where $t^{*}$ is one of the 
values determined above by fitting the norms of $u$, we get  the 
figure shown in Fig.~\ref{gKPn26gaussxm} and the values  $\alpha_{1}=-0.0488$ 
and $\alpha_{2}=
   0.2281$. Thus the negative value for $\alpha_{1}$ could again 
   imply a blow-up, but it is probable that we simply did not get 
   close enough to the blow-up to be able to decide. The small value 
   of $|\alpha_{1}|$ is also consistent with a finite value of 
   $x^{*}$. 

The above example clearly indicates a blow-up for solutions to the 
gKP I equation for $n=2$ for initial data with sufficiently small 
energy. An interesting question is again to identify the condition on the
initial data for blow up. For initial data of the form (\ref{initial}) we 
find that the solution stays regular for $\beta\leq 2$. But there appears to 
be blow-up for initial data with $\beta\geq3$ for which the energy is 
positive. Again it is numerically difficult to 
get closer to the actual threshold. 

\begin{figure}[ht]
   \centering
    \includegraphics[width=0.7\textwidth]{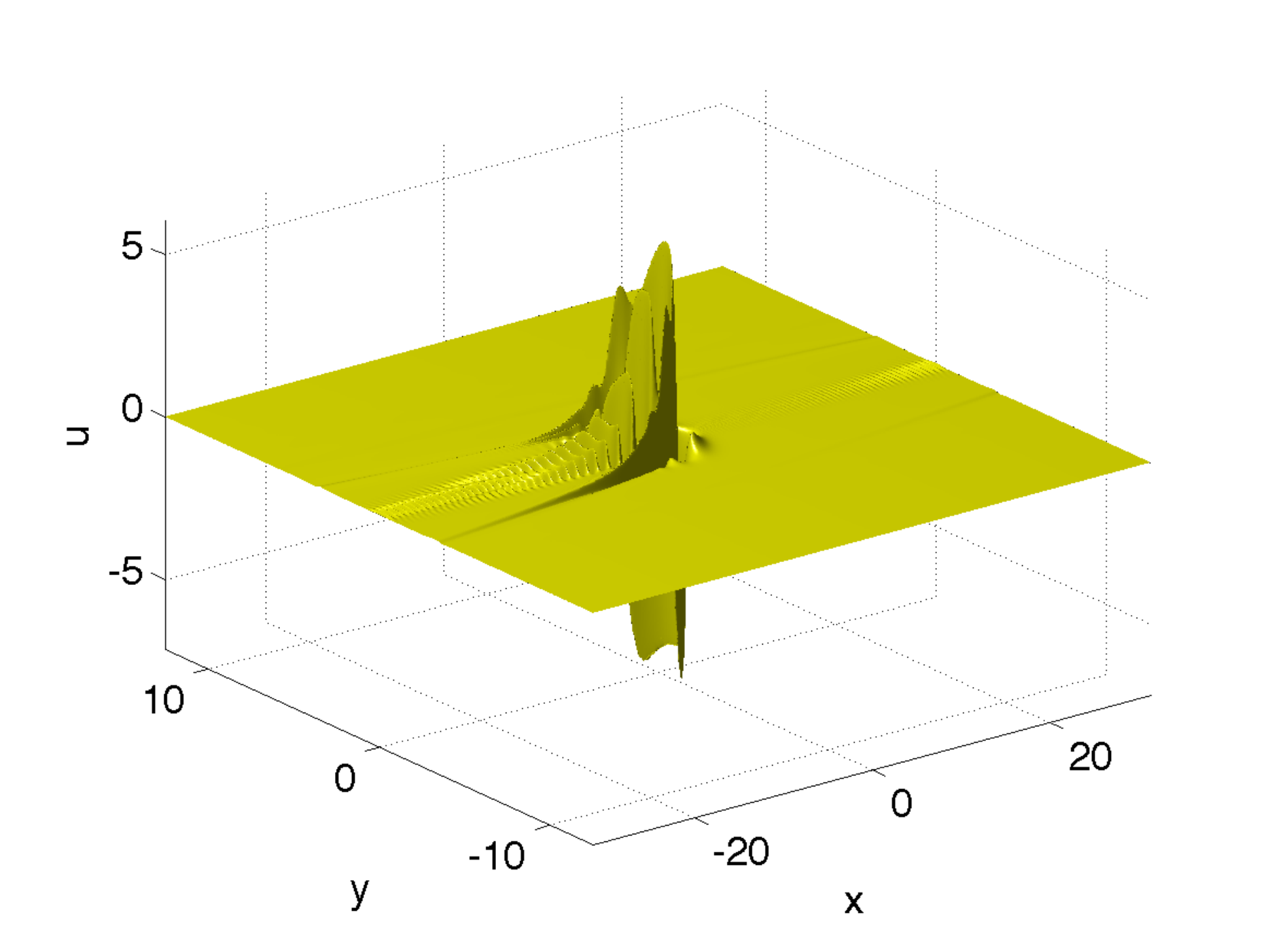}
   \caption{Solution to the gKP II equation (\ref{gKP}) with $\lambda=1$, 
   $n=4/3$ for the 
   initial data $u_{0}=6\,\partial_{xx}\exp(-(x^{2}+y^{2}))$ at 
   $t=2$.}
   \label{gKPIIn26gauss}
\end{figure}

Blow-up does not seem to appear though for the corresponding gKP II 
equation. If we consider the same initial data that led to blow-up 
for gKP I, the energy is again negative, but there is no indication 
of blow-up in this case as can be seen in Fig.~\ref{gKPIIn26gauss}. 
The computation is carried out with $L_{x}=10$, $L_{y}=4$, 
$N_{x}=N_{y}=2^{10}$ and $N_{t}=10^{3}$ time steps for $t<0.1$. The 
typical tails, for gKP II to the left, can be clearly seen as well as 
the dispersive oscillations propagating in the same direction. 
Because of the periodicity, they reenter on the right. But it appears 
that generic localized initial data will be just radiated away to 
infinity for gKP II for all $n$.

\begin{figure}[ht]
   \centering
   \includegraphics[width=0.49\textwidth]{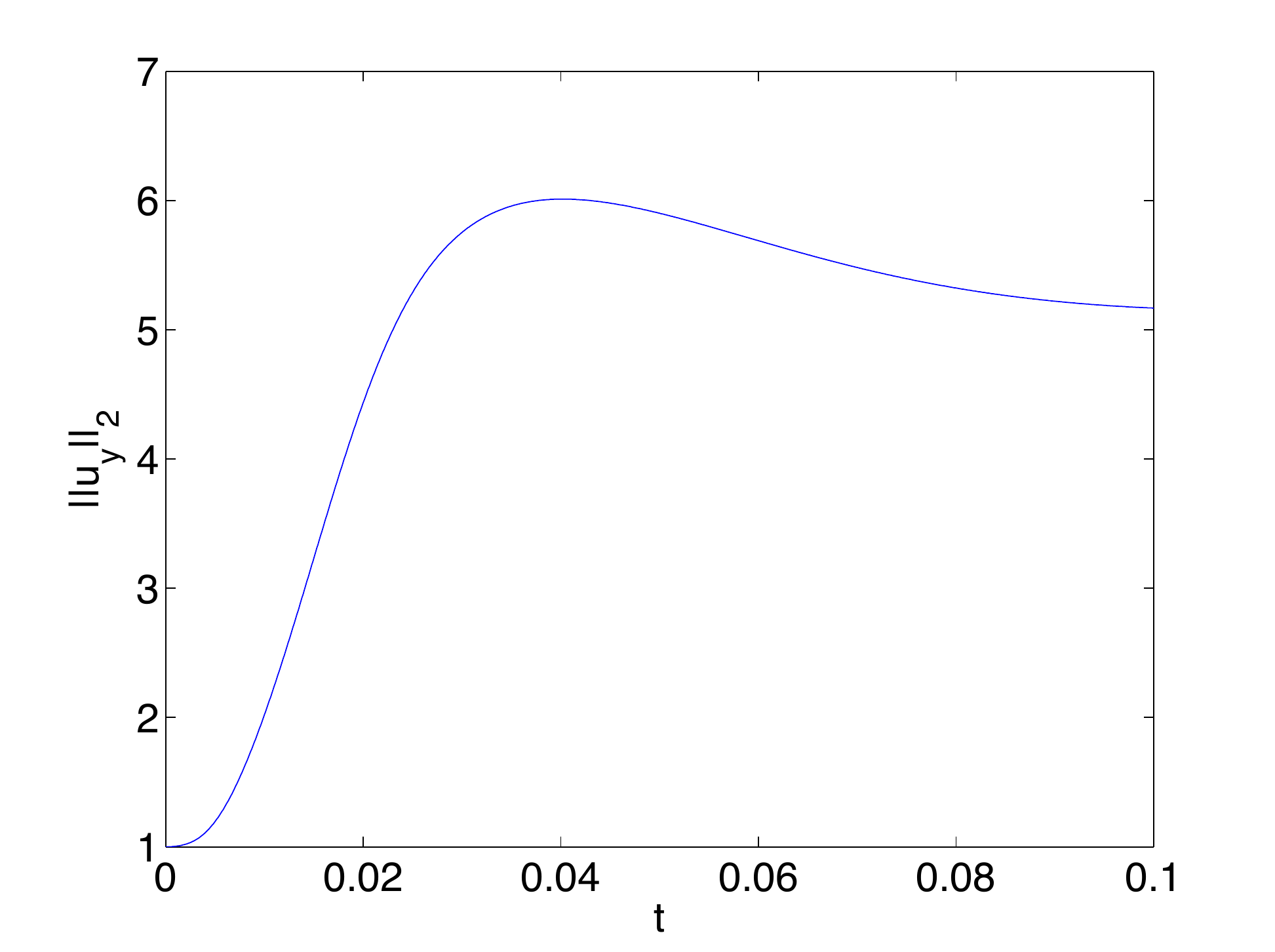}
   \includegraphics[width=0.49\textwidth]{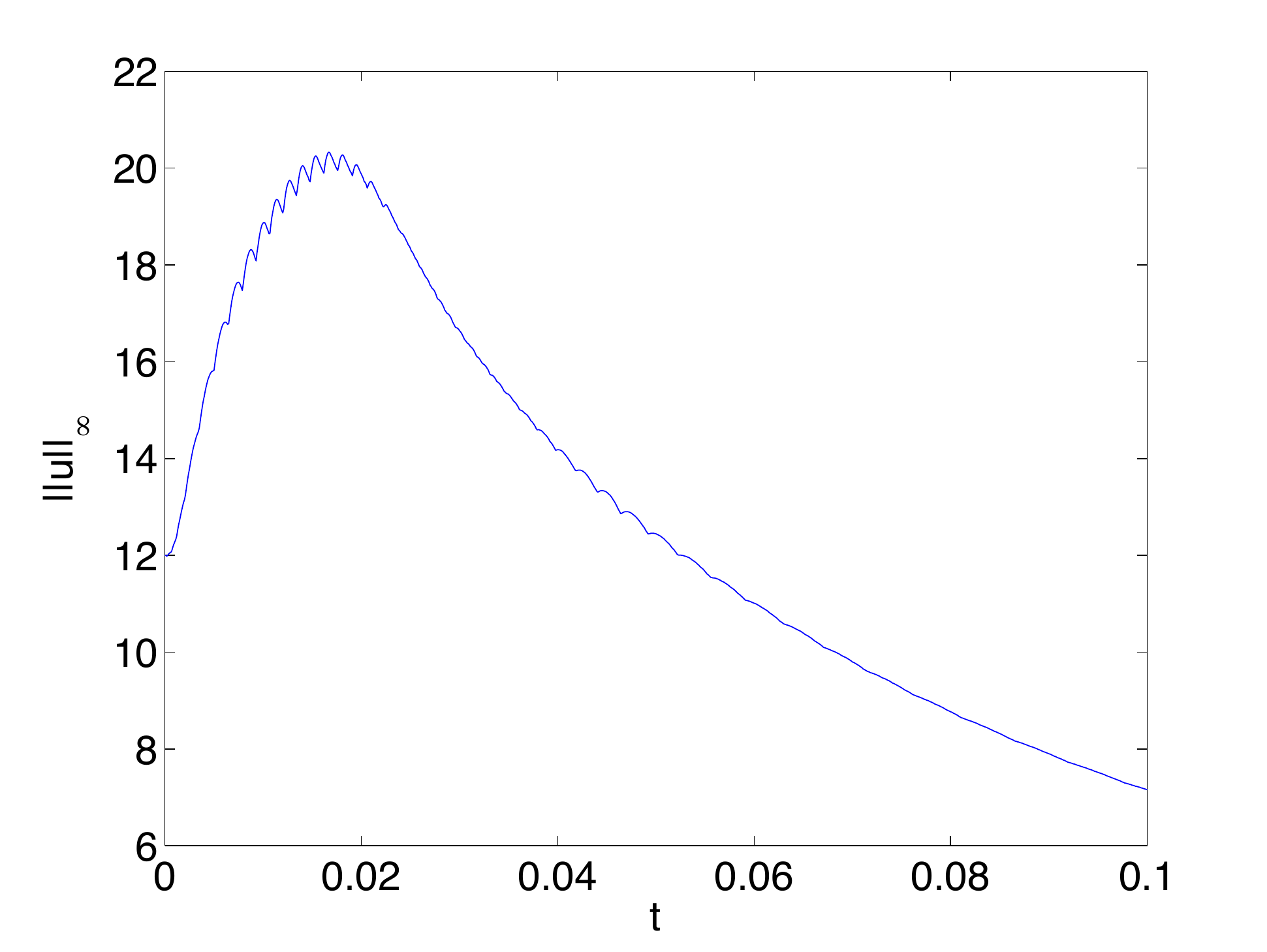}
   \caption{Norms of the solution to the equation (\ref{gKP}) with $\lambda=1$, 
   $n=2$ for the 
   initial data $u_{0}=6\,\partial_{xx}\exp(-(x^{2}+y^{2}))$; the 
   $L_{2}$ norm of $u_{y}$ on the left, 
   $L_{\infty}$ norm of $u$ on right.}
   \label{gKPIIn26gaussnorm}
\end{figure}

This is even more obvious from certain norms of the solution as shown 
in  Fig.~\ref{gKPIIn26gaussnorm}. Both the $L_{2}$ norm of $u_{y}$ 
and the $L_{\infty}$ norm of $u$ appear to decrease monotonically. 
Note that the same qualitative behavior is obtained for twice the 
initial data considered in Fig.~\ref{gKPIIn26gauss}. Thus we did not 
find an indication for blow-up in gKP II solutions for $n=2$. 
\section{The supercritical cases $n=3$ and $n=4$}
The numerical experiments of the previous sections indicate blow-up 
for gKP I solutions, but not for gKP II solutions. Whereas this does 
not exclude blow-up for the latter, it definitely does not appear at 
comparable energies or masses as for gKP I, presumably due to the 
defocusing effect of gKP II. In this section we will study whether 
there is blow-up in gKP II solutions for $n=3$ and $n=4$. Since the 
behavior of gKP I solutions appears to be similar to the case $n=2$, 
we concentrate here on gKP II solutions. 

\begin{figure}[ht]
   \centering
   \includegraphics[width=\textwidth]{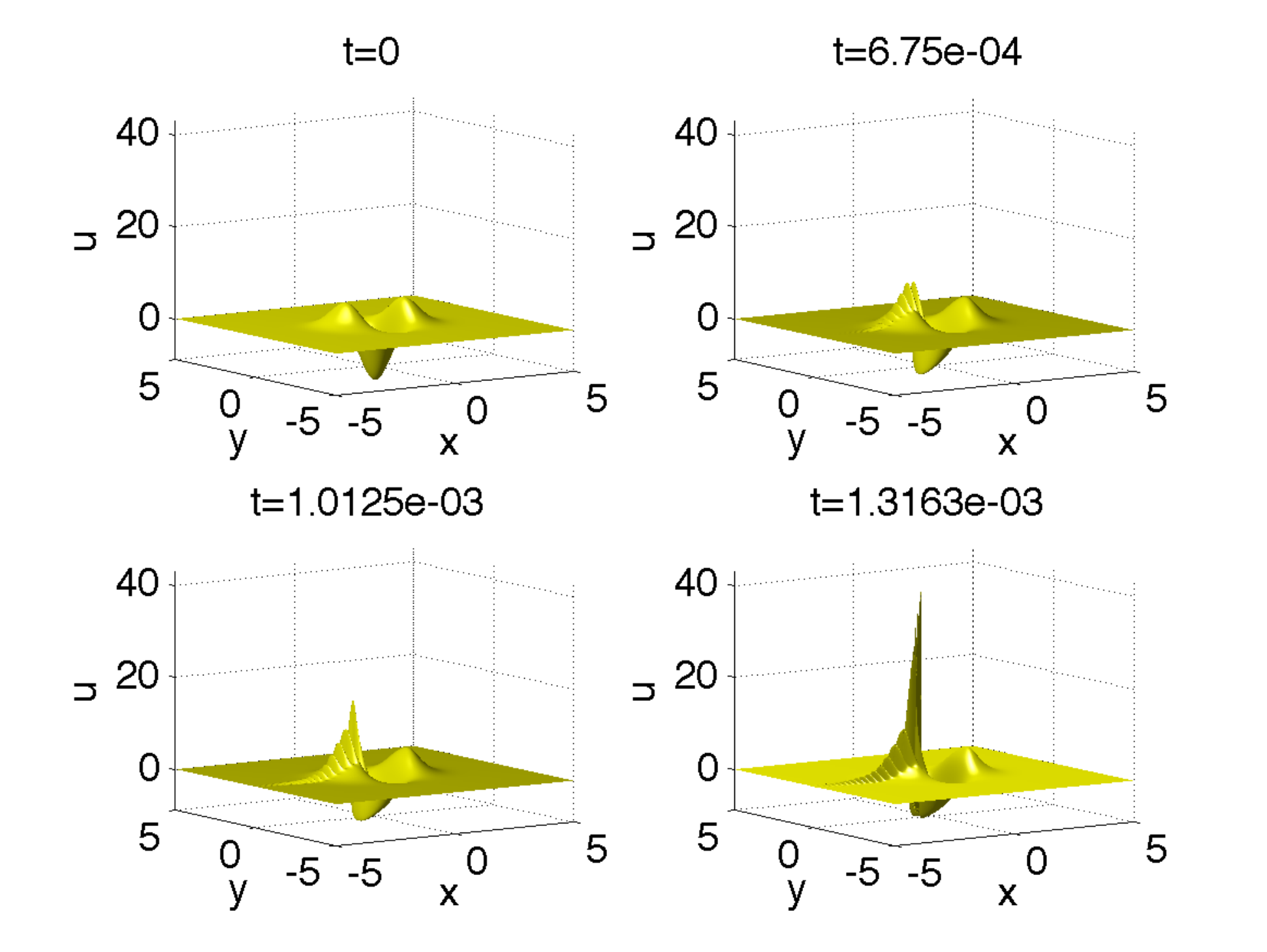}
   \caption{Solution to the gKP II equation (\ref{gKP}) with 
   $n=3$ for the 
   initial data $u_{0}=3\,\partial_{xx}\exp(-(x^{2}+y^{2}))$ for 
   several times. }\label{gKPIIn36gaussu}
\end{figure}

For initial data of small enough mass, the gKP II solutions appear 
again to be simply radiated away. But for initial data of the form 
(\ref{initial}), the situation changes 
for large enough $\beta$. In Fig.~\ref{gKPIIn36gaussu} we show the 
gKP II solution for $\beta=6$ for several values of $t$. The 
computation is carried out with $L_{x}=5$, $L_{y}=4$ and 
$N_{x}=N_{y}=2^{12}$ Fourier modes and $N_{t}=20000$ time steps for $t<0.0014$. The 
code is stopped once the numerically computed relative energy is conserved 
to less than $10^{-3}$. It can be seen that this time the maximum of 
the solution appears to blow up. 

\begin{figure}[ht]
   \centering
   \includegraphics[width=0.49\textwidth]{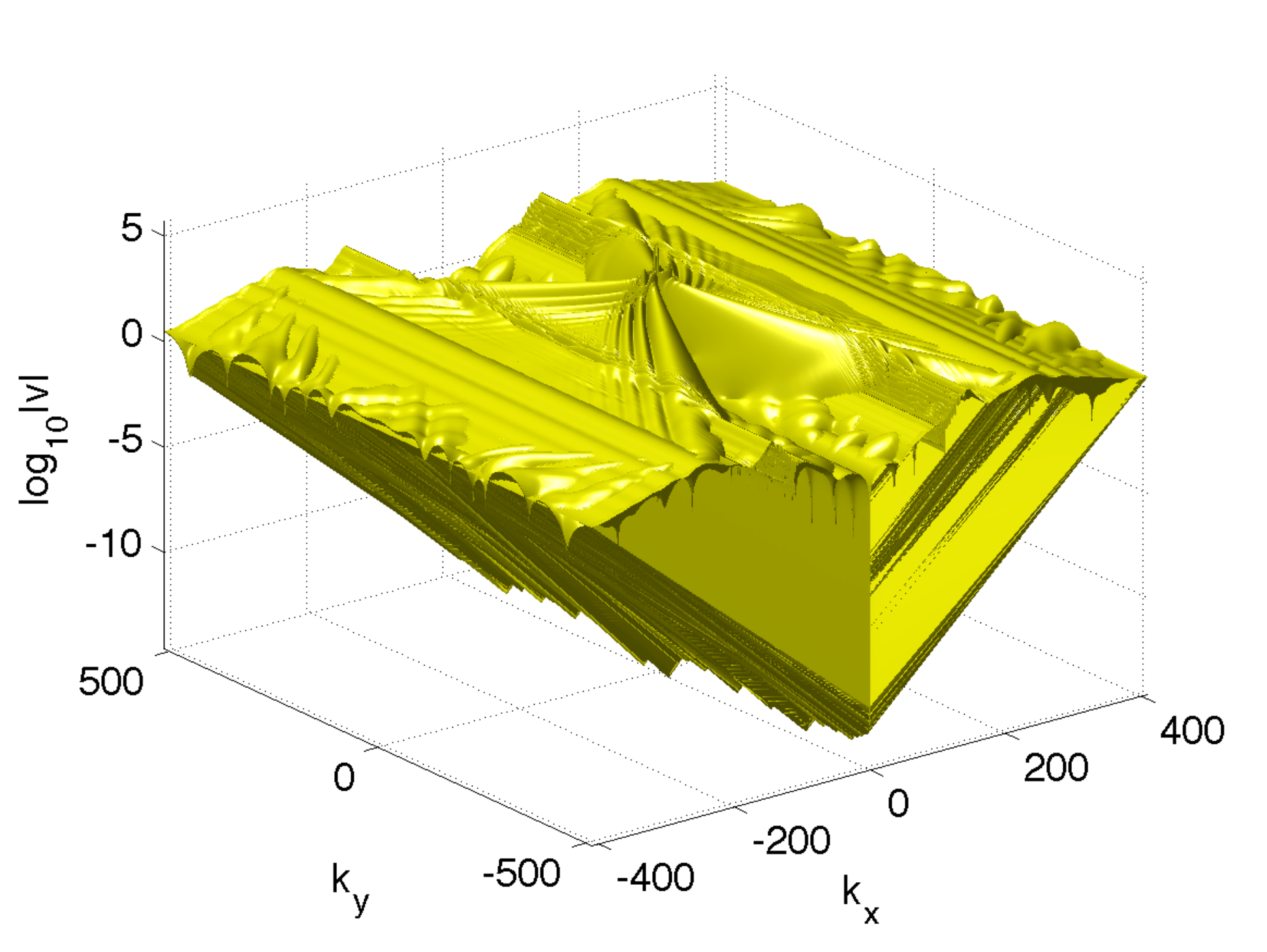}
   \caption{Modulus of the Fourier coefficients of the solution 
   $u$ in Figure \ref{gKPIIn43gaussu} at $t = 
   1.3163*10^{-3}$.}\label{gKPIIn36gaussfourier}
\end{figure}

The resolution in Fourier space at the last recorded time can be seen 
in Fig.~\ref{gKPIIn36gaussfourier}. Again due to the behavior 
indicated by the rescaling (\ref{gKP4}), there 
are strong gradients especially in $y$ which imply that the loss in 
resolution in Fourier space is predominantly in $y$-direction. 

\begin{figure}[ht]
   \centering
   \includegraphics[width=0.49\textwidth]{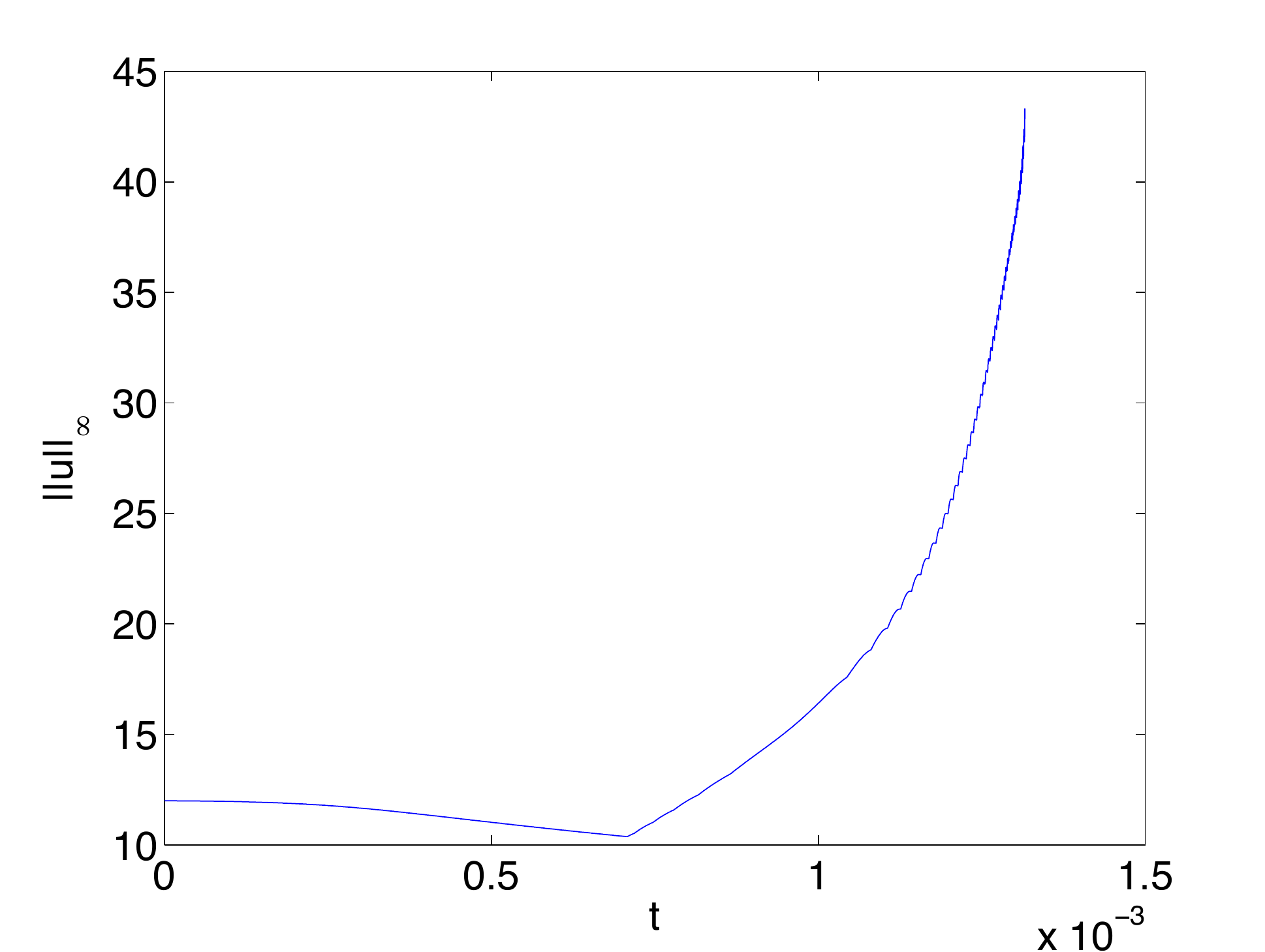}
   \includegraphics[width=0.49\textwidth]{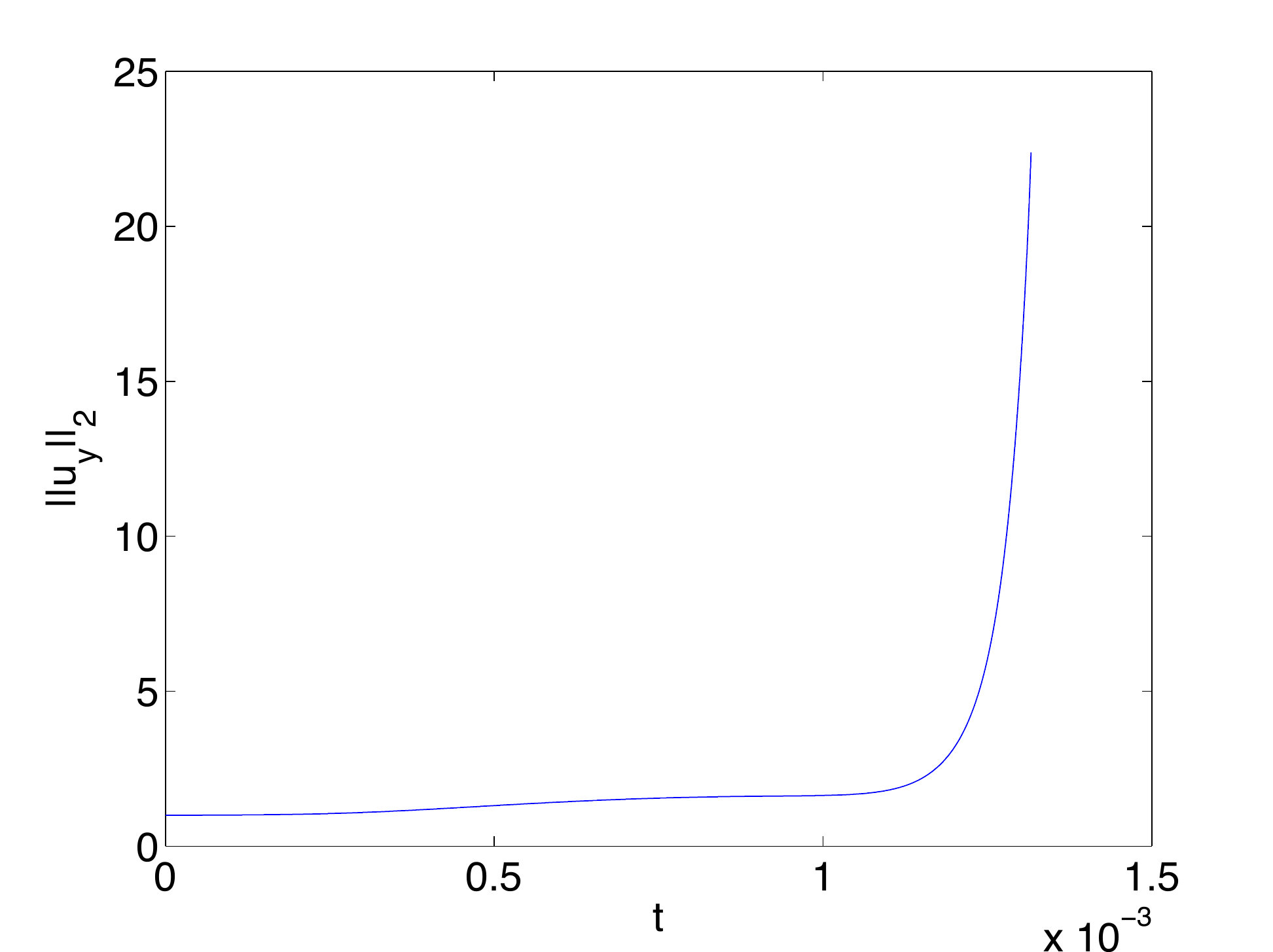}
   \caption{Norms of the solution to the gKP I equation (\ref{gKP}) with 
   $n=3$ for the 
   initial data $u_{0}=6\,\partial_{xx}\exp(-(x^{2}+y^{2}))$) in 
   dependence of time; on the left the $L_{\infty}$ norm  of the 
   solution $u$, on the right the $L_{2}$ norm of 
   $u_{y}$.}\label{gKPIIn36gaussuinf}
\end{figure}

As for the gKP I blow-ups in the previous sections, 
both the $L_{2}$ norm of $u_{y}$ and the $L_{\infty}$ norm of $u$ 
indicate a blow-up as is clear from Fig.~\ref{gKPIIn36gaussuinf}. 

\begin{figure}[ht]
   \centering
   \includegraphics[width=0.49\textwidth]{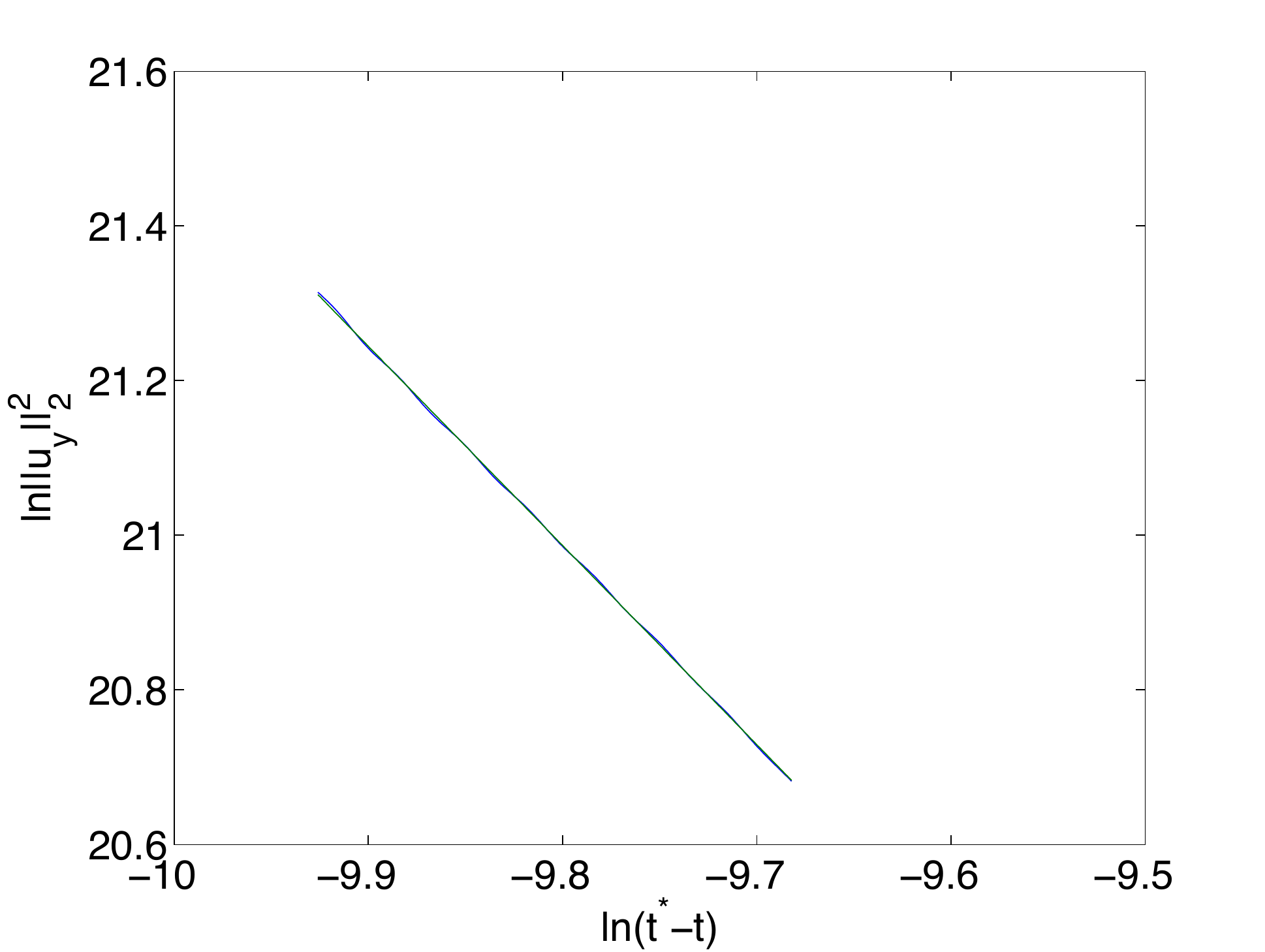}
   \includegraphics[width=0.49\textwidth]{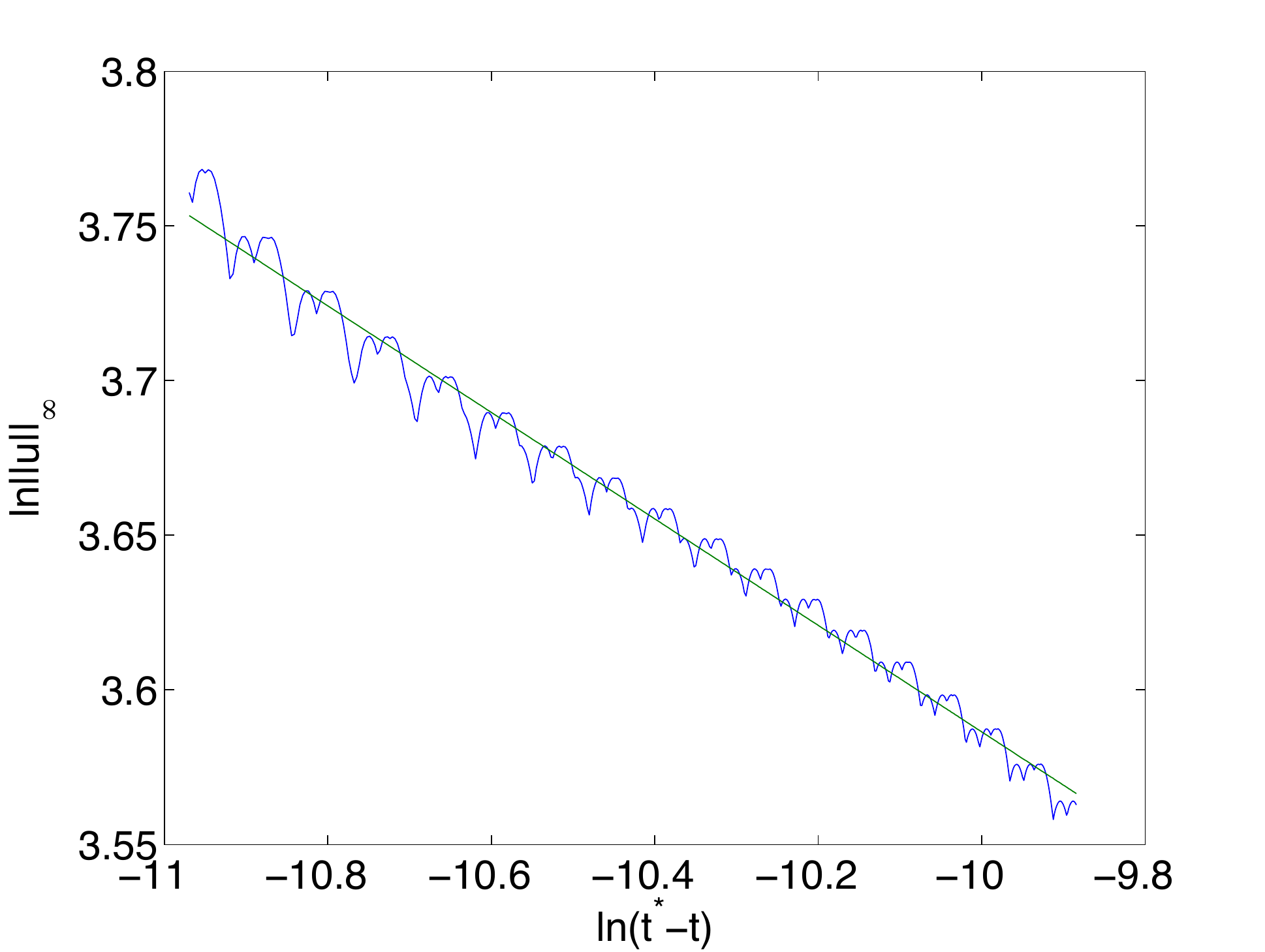}
   \caption{Fitting of $\ln||u_{y}||_{2}^{2}$ (left) and 
 $\ln||u||_{\infty}$ (right)  in blue to $c \ln(t^{*}-t)+C$ in green 
 for the situation shown in 
 Fig.~\ref{gKPIIn36gaussu}.}\label{gKPIIn36gaussfit}
\end{figure}

This already indicates a similar type of blow-up as for gKP I. This 
is confirmed by an asymptotic analysis of the norms in 
Fig.~\ref{gKPIIn36gaussuinf} close to the blow-up. 
As before we fit $\ln ||u||_{\infty}$ for the solution shown in 
Fig.~\ref{gKPIIn36gaussu} to $\ln C_1+c_{1}\ln(t^* - t)$ and 
$\ln||u_{y}||_{2}$ to $ \ln C_2+c_{2}\ln(t^* - t)$. 
For the fitting of $||u||_{\infty}$ we use the 500 last        
computed points and get $ C_1 =   1.8650$, $c_1 = -0.1721$ and  $t^* 
=   1.3335*10^{-3}$. 
Despite the oscillations in 
$||u||_{\infty}$ in Fig.~\ref{gKPIIn36gaussfit}, the value of 
$c_{1}$ is close to the expected $-2/9$. The fitting is less 
convincing for the $L_{2}$ norm of $u_{y}$. For the last 200 time 
steps we get   $ C_2 =  -4.256 $, $c_2 = -2.576$ and  $t^* 
=1.365*10^{-3}$. The value of $t^{*}$ is clearly too large for what is to 
be expected by the breaking of the code. This indicates that the 
$L_{2}$ norm of $u_{y}$ is again affected by errors in the gradient 
close to the boundaries of the computational domain where we do not 
have enough resolution. But the results 
are compatible with what was obtained in the previous sections for 
gKP I, and the fitting for the $L_{\infty}$ of $u$ appears to be 
conclusive. 

\begin{figure}[ht]
   \centering
   \includegraphics[width=\textwidth]{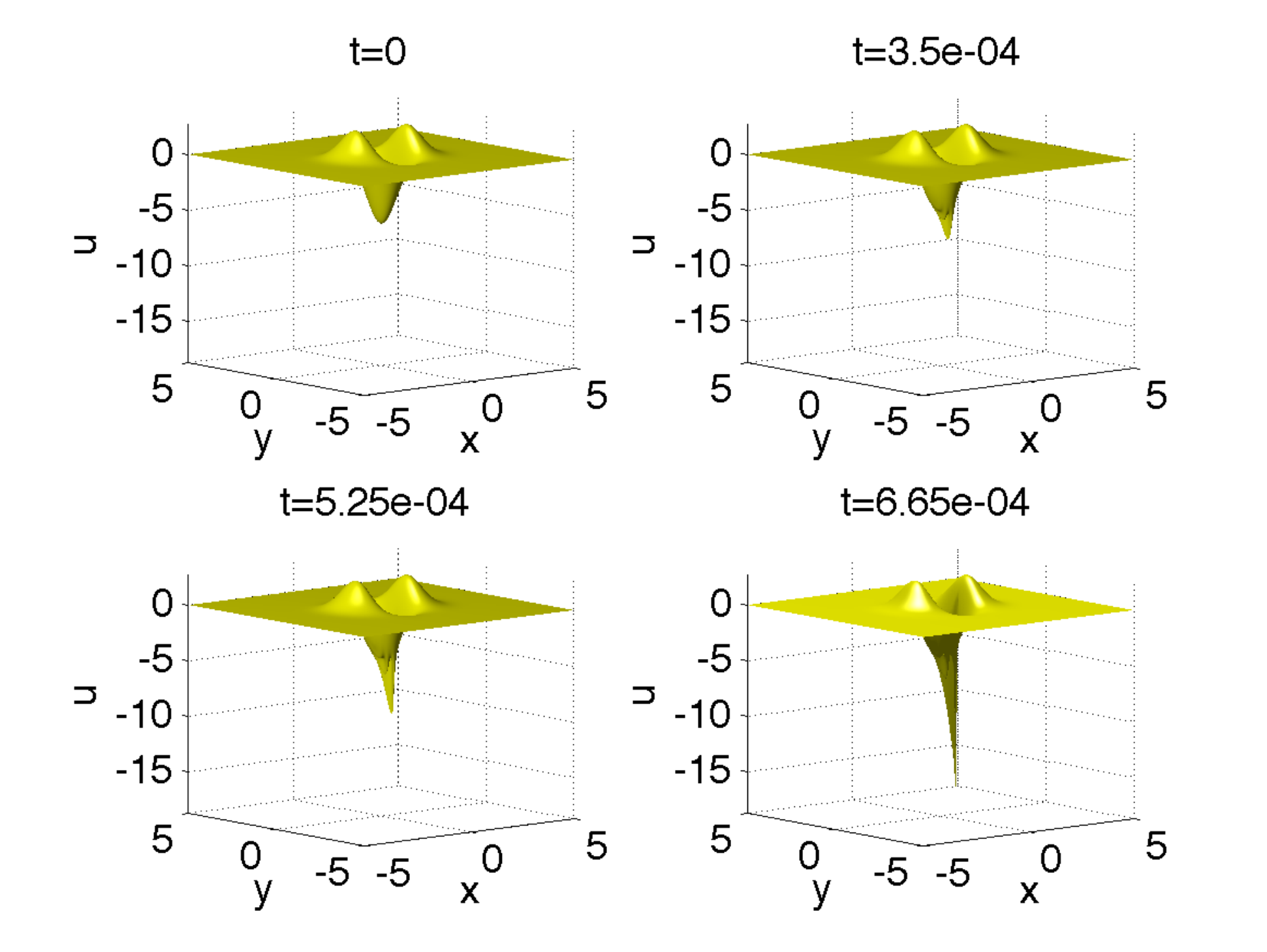}
   \caption{Solution to the gKP II equation (\ref{gKP}) with 
   $n=4$ for the 
   initial data $u_{0}=6\,\partial_{xx}\exp(-(x^{2}+y^{2}))$ for 
   several times. }\label{gKPIIn43gaussu}
\end{figure}

The same experiment will be carried for $n=4$. For small enough 
$\beta$ in (\ref{initial}), the initial data will again just be radiated away. This 
changes for $\beta=3$ as can be seen in  Fig.~\ref{gKPIIn43gaussu}. 
The 
computation is carried out with $L_{x}=L_{y}=5$ and 
$N_{x}=N_{y}=2^{12}$ Fourier modes and $N_{t}=20000$ time steps for 
$t<0.0007$. The 
code stops shortly after the numerically computed relative energy is conserved 
to less than $10^{-3}$, and we record only the data for which 
$\Delta<10^{-3}$. Since $n$ is again even, the minimum blows up as 
before in these cases. 

\begin{figure}[ht]
   \centering
   \includegraphics[width=0.49\textwidth]{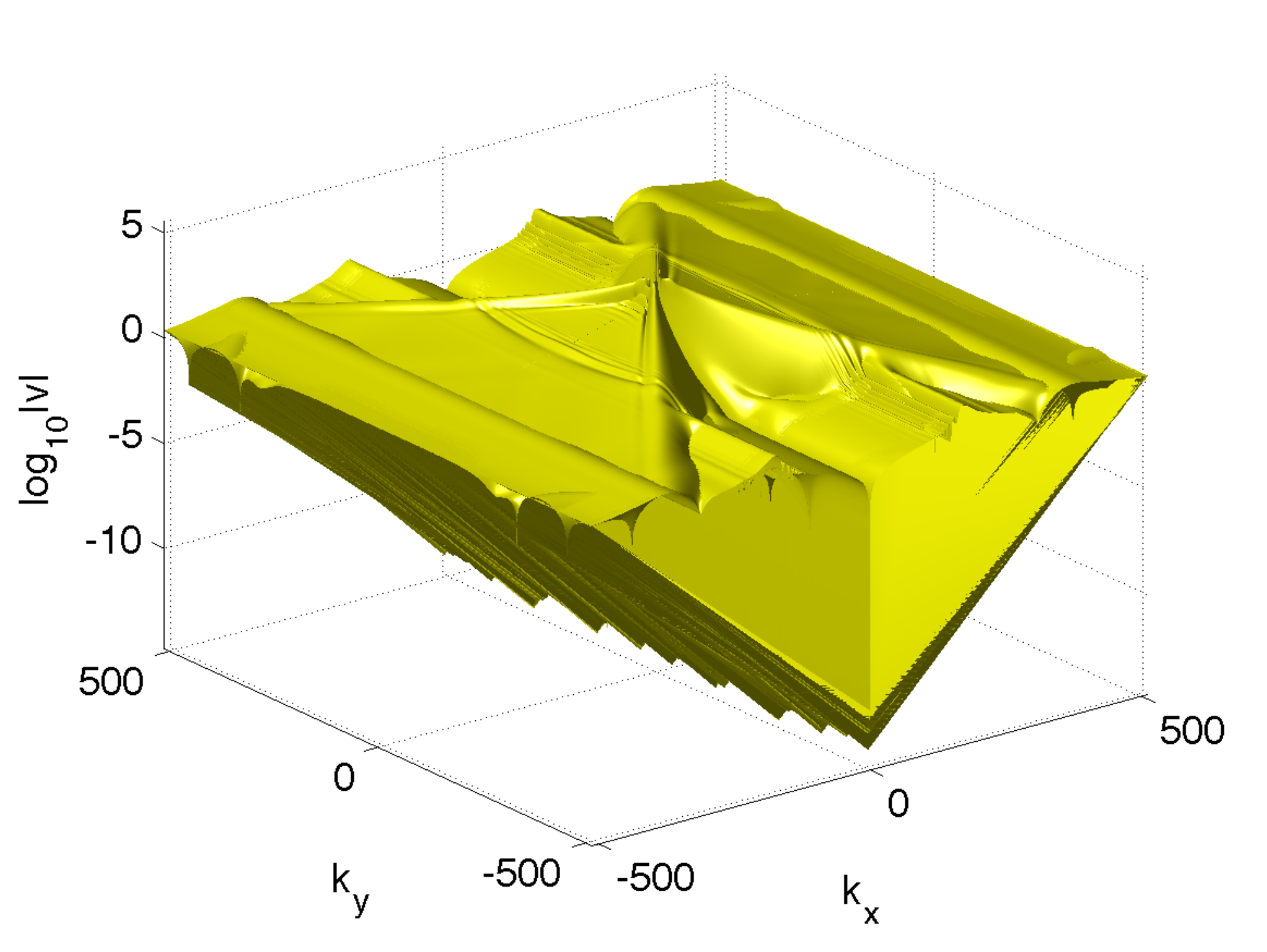}
   \caption{Modulus of the Fourier coefficients of the solution 
   $u$ in Figure \ref{gKPIIn43gaussu} at $t = 
   6.65*10^{-4}$.}\label{gKPIIn43gaussfourier}
\end{figure}

The resolution in Fourier space at the last recorded time can be seen 
in Fig.~\ref{gKPIIn43gaussfourier}. 

\begin{figure}[ht]
   \centering
   \includegraphics[width=0.49\textwidth]{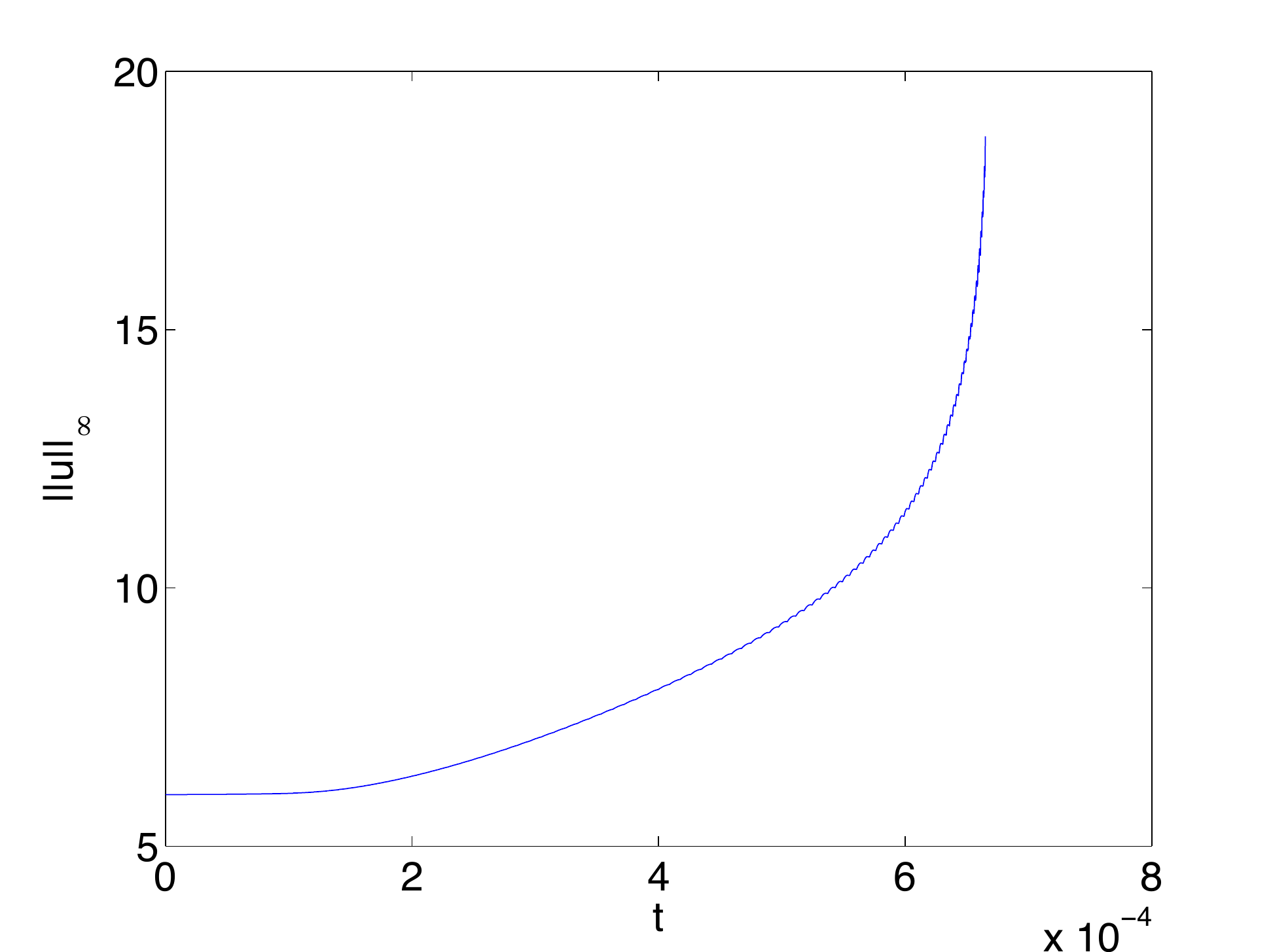}
   \includegraphics[width=0.49\textwidth]{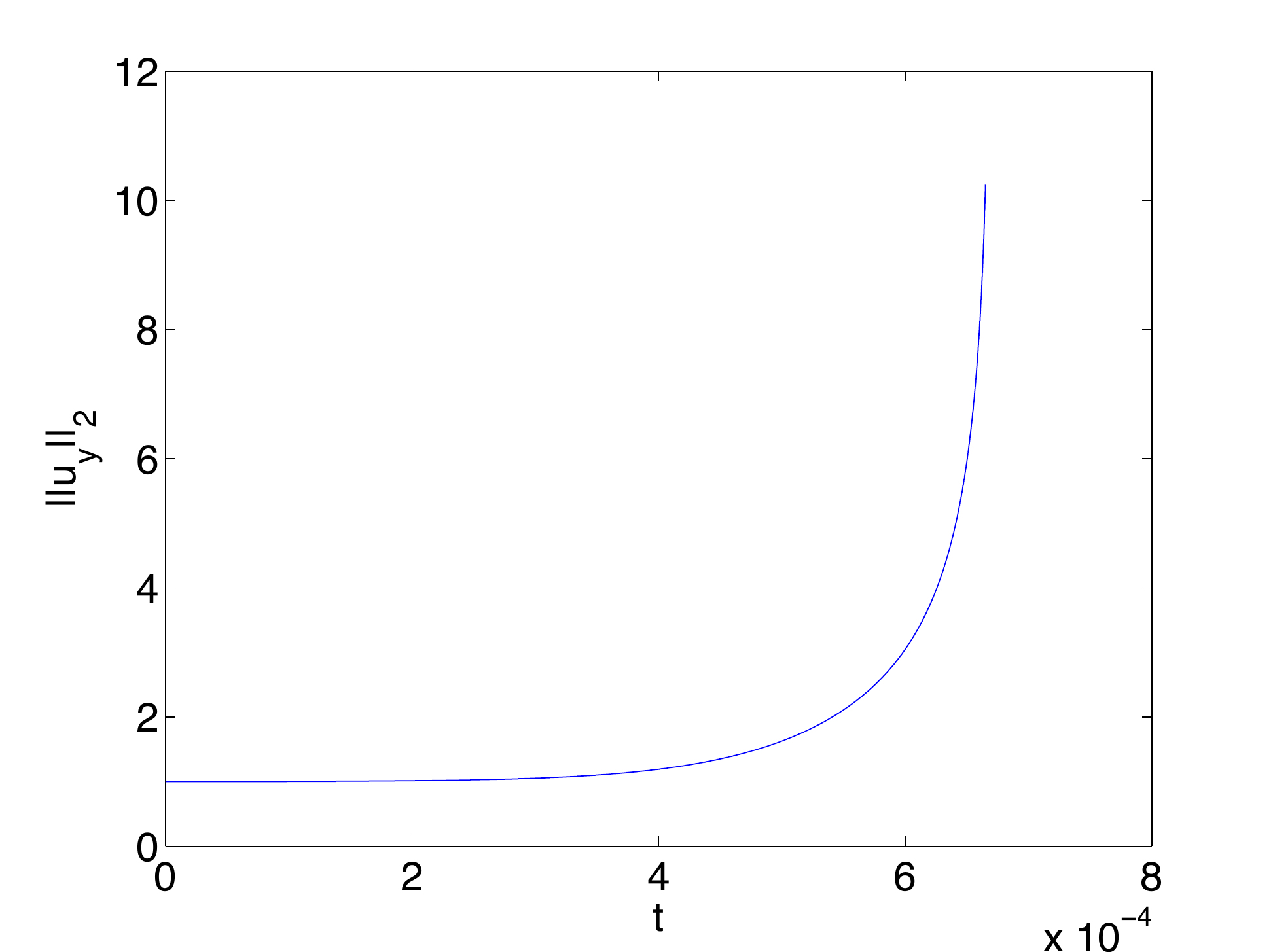}
   \caption{Norms of the solution to the gKP I equation (\ref{gKP}) with 
   $n=4$ for the 
   initial data $u_{0}=3\,\partial_{xx}\exp(-(x^{2}+y^{2}))$) in 
   dependence of time; on the left the $L_{\infty}$ norm  of the 
   solution $u$, on the right the $L_{2}$ norm of 
   $u_{y}$.}\label{gKPIIn43gaussuinf}
\end{figure}

Both the $L_{2}$ norm of $u_{y}$ and the $L_{\infty}$ norm of $u$ 
clearly show a blow-up as can be seen in Fig.~\ref{gKPIIn43gaussuinf}. 

\begin{figure}[ht]
   \centering
   \includegraphics[width=0.49\textwidth]{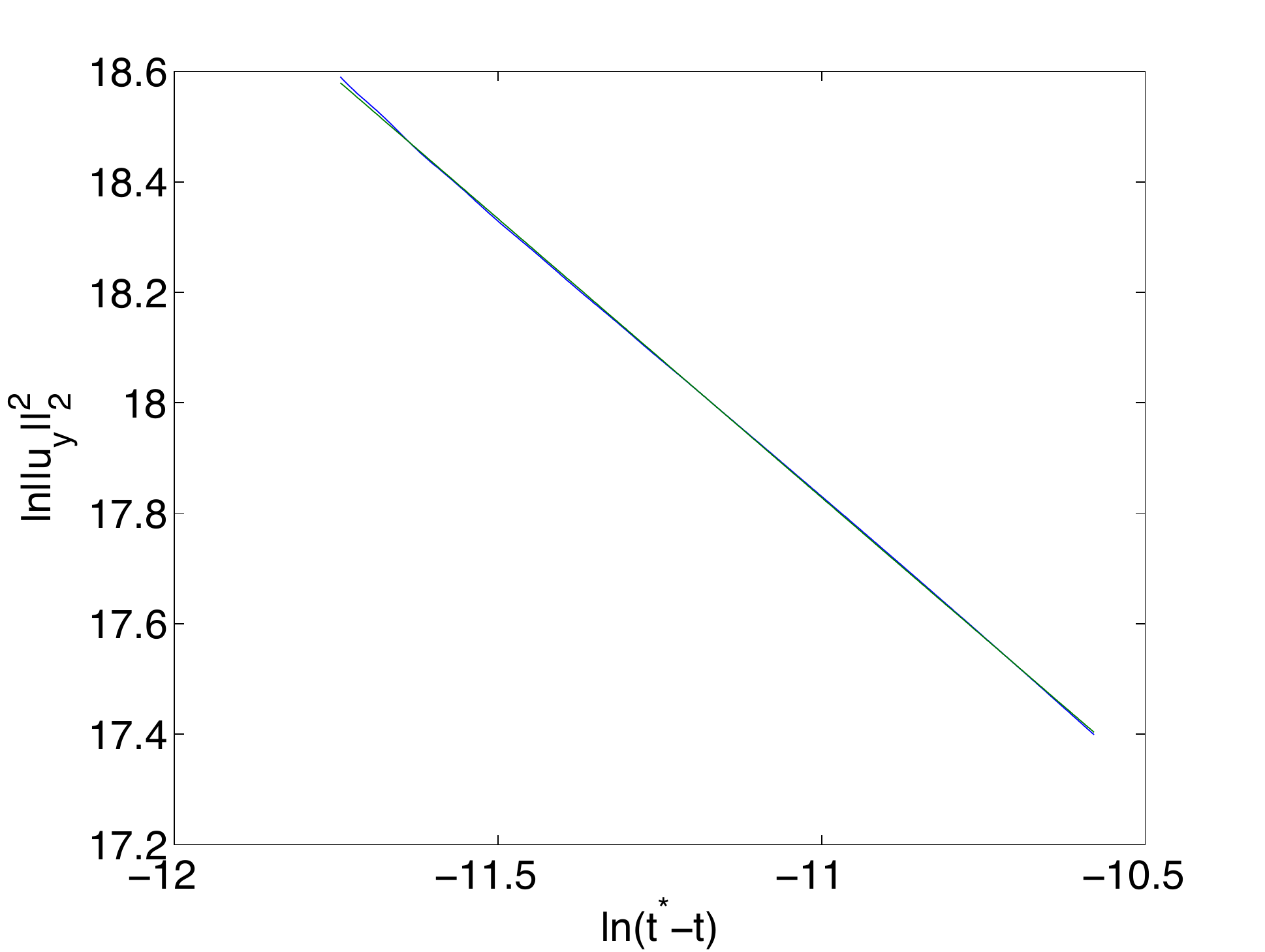}
   \includegraphics[width=0.49\textwidth]{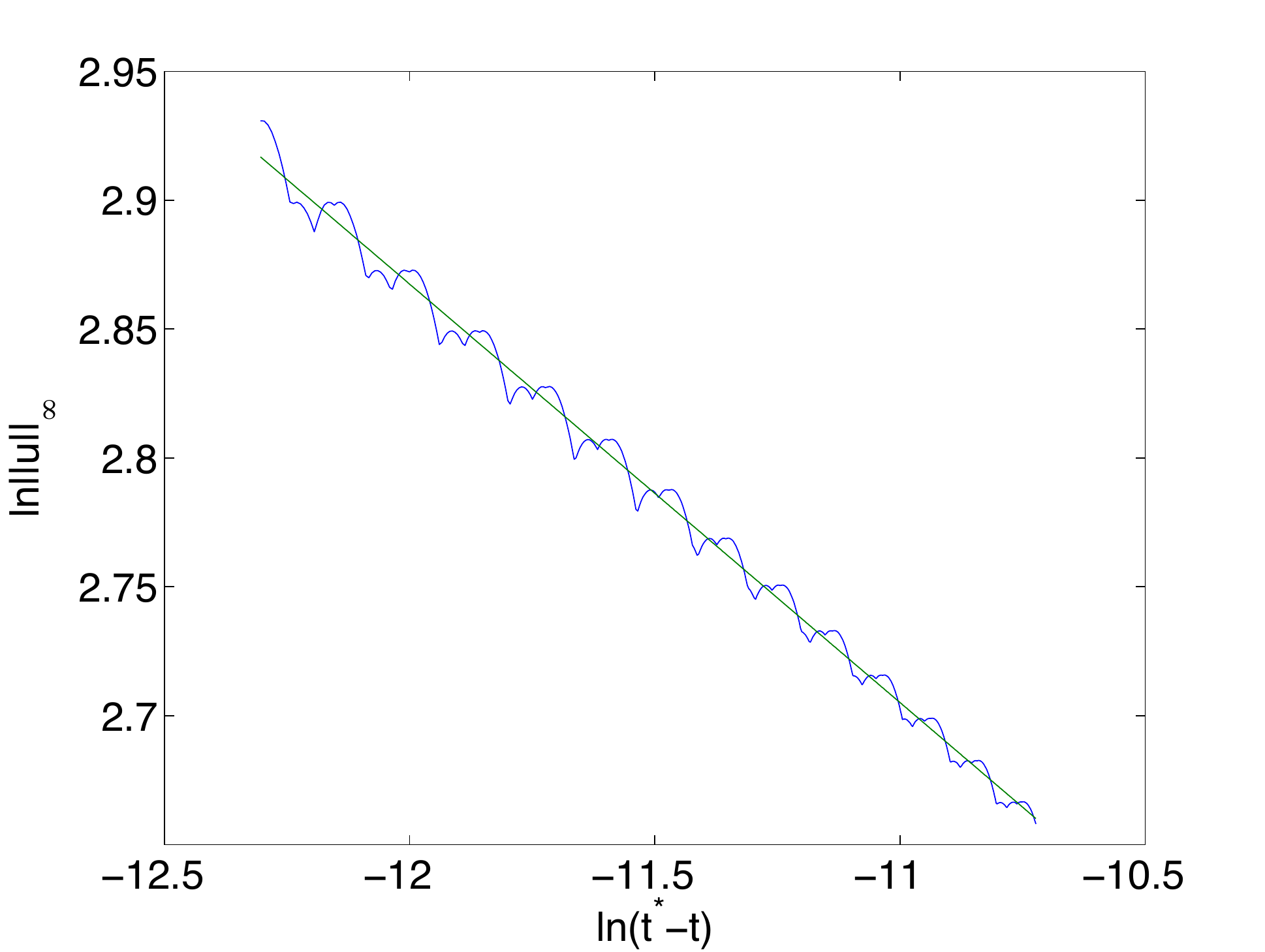}
   \caption{Fitting of $\ln||u_{y}||_{2}^{2}$ (left) and 
 $\ln||u||_{\infty}$ (right)  in blue to $c \ln(t^{*}-t)+C$ in green 
 for the situation shown in 
 Fig.~\ref{gKPIIn43gaussu}.}\label{gKPIIn43gaussfit}
\end{figure}

To identify the type of blow-up we again fit $\ln ||u||_{\infty}$ for the solution shown in 
Fig.~\ref{gKPIIn43gaussu} to $\ln C_1+c_{1}\ln(t^* - t)$ and 
$\ln||u_{y}||_{2}^{2}$ to $ \ln C_2+c_{2}\ln(t^* - t)$. 
For the fitting of $||u||_{\infty}$ we use the 500 last           
computed points and get $ C_1 =   0.91937$, $c_1 = -0.1623$ and  $t^* 
=   6.6953*10^{-4}$. 
The quality of the fitting can be seen in  
Fig.~\ref{gKPIIn43gaussfit}. The value of 
$c_{1}$ is very close to the expected $-1/6$. The fitting is again less 
convincing for the $L_{2}$ norm of $u_{y}$ for the last 200 time steps. We       
get   $ C_2 =  8.3069 $, $c_2 = -0.858$ and  $t^*      
=6.7126*10^{-4}$, which is, however, in accordance with expectations 
($c_{1}=-2/3$) and 
what was found  for the $L_{\infty}$ norm of $u$.

\section{Conclusion}
In this paper we have numerically solved the Cauchy problem for the 
gKP equations for smooth localized initial data with a single minimum 
and single maximum. 
The results  can be summarized in the following \\
\textbf{Conjecture:}
Consider the Cauchy problem for the gKP equations (\ref{gKP}) with $n\in 
\mathbb{Q}$, i.e., $n=p/q$ with $p,q\in \mathbb{N}$ and $p,q$ 
coprime,  with initial data $u_{0}(x,y)$ with a single global 
minimum (for $p$ even) or a single global maximum (for $p$ odd) such that $\partial_{x}^{-1}u\in 
\mathcal{S}(\mathbb{R}^{2})$. Then 
\begin{itemize}
    \item  for $n<4/3$, the solution is smooth for all $t$.

    \item  for gKP II, the solution is smooth for all $t$ for  
    $n\leq 2$.

    \item  for gKP I with $n=4/3$, initial data with sufficiently 
    small energy and sufficiently large mass
    lead to blow-up at $t^{*}<\infty$; asymptotically for $t\sim 
    t^{*}$, the solution is given by a rescaled (via (\ref{gKP4})) soliton 
    (\ref{soliton}) where the scaling factor $L\propto1/\tau$ for 
    $\tau\to\infty$. This implies the blow-up is characterized by 
    \begin{equation}
        ||u||_{\infty}\propto \frac{1}{(t^{*}-t)^{3/4}},\quad 
	||u_{y}||_{2}\propto \frac{1}{t^{*}-t}.
        \label{scal43}
    \end{equation}

    \item  for gKP I with $n>4/3$ and gKP II with $n>2$,  initial data with sufficiently 
    small energy and sufficiently large mass
    lead to blow-up at $t^{*}<\infty$; asymptotically for $t\sim 
    t^{*}$, the solution is given by a localized solution to 
    (\ref{ODE}), which is conjectured to exist and to be unique, 
    rescaled (via (\ref{gKP4})) where the scaling factor $L\propto\exp(\kappa 
    \tau)$ for 
    $\tau\to\infty$ with $\kappa$ a negative constant.  This implies the blow-up is characterized by 
    \begin{equation}
        ||u||_{\infty}\propto \frac{1}{(t^{*}-t)^{2/(3n)}},\quad 
	||u_{y}||_{2}\propto \frac{1}{(t^{*}-t)^{(1+4/n)/6}}.
        \label{scal2}
    \end{equation}
\end{itemize}
Obviously numerical experiments can just give indications on possible 
theorems, but the results presented in this paper are conclusive 
within the limitations discussed in the text. 

An important open question is the condition on the initial data to 
lead to blow-up. It is unclear whether the relevant quantity  is simply the energy, or 
whether this is related to the mass and energy of solitons in the 
cases where these exist (recall that there are no lumps for gKP I 
for $n\geq 4$ and for gKP II in general). For 
gKdV it was proven in \cite{MMR2012_I} in the $L_{2}$ critical case 
$n=4$ that for initial data in the vicinity of the soliton, the criterion 
for blow-up is that the energy of the initial data is smaller than 
the soliton initial data, whereas the mass has to be greater than the 
soliton mass. Numerical experiments indicate, see for instance 
\cite{KP13}, that this is also the case for more general localized initial 
data. Thus it would be important to study solutions of 
(\ref{soliton}), which are known to be unstable \cite{BS97a,BS97b}, 
and to see what type of perturbation if any leads to blow-up.  To this 
end, these solitons, which are not known analytically, have to 
be constructed first. Note, 
however, that there is a whole family of soliton solutions 
parametrized by the speed $c$. For gKdV with $n=4$ the mass of the 
solitons is 
independent of $c$, and the same holds for the energy which vanishes 
for all $c$. 

Of similar interest as the solitons would be to solve the PDE 
(\ref{ODE}) numerically which appears in the asymptotic description 
of blow-up and to check whether this gives the expected asymptotic 
description of the blow-up. A further interesting question is related to the 
asymptotic behavior of the location of the  blow-up 
which presumably stays finite. But it is unclear how such a finite 
value is approached
These questions  will be the subject of further work.

\section*{Acknowledgments}
We thank J.-C.~Saut for helpful remarks and hints. This work has been supported by the project FroM-PDE funded by the European
Research Council through the Advanced Investigator Grant Scheme 
and the ANR via the program ANR-09-BLAN-0117-01.

\bibliographystyle{siam}
\bibliography{Num_gKP_bib}

\begin{thebibliography}{10}

\bibitem{AF}
{\sc M.~Ablowitz and A.~Fokas}, {\em On the inverse scattering and direct
  linearizing transforms for the {K}adomtsev-{P}etviashvili equation}, Phys.
  Lett. A, 94(2) (1983), pp.~67--70.

\bibitem{BPP}
{\sc M.~Boiti, F.~Pempinelli, and A.~K. Pogrebkov}, {\em Solutions of the {KPI}
  equation with smooth initial data}, Inverse Problems, 10 (1994),
  pp.~504--519.

\bibitem{CoxMatthews2002}
{\sc S.~M. Cox and P.~C. Matthews}, {\em Exponential time differencing for
  stiff systems}, J. Comp. Phys., 176 (2002), pp.~430--455.

\bibitem{BS97a}
{\sc A.~de~Bouard and J.-C. Saut}, {\em Solitary waves of the generalized {KP}
  equations}, Ann. Inst. Henri Poincar{\'e}, Anal. Non LinŽaire, 14 (1997),
  pp.~211--236.

\bibitem{BS97b}
\leavevmode\vrule height 2pt depth -1.6pt width 23pt, {\em Symmetry and decay
  of the generalized {K}adomtsev-{P}etviashvili solitary waves}, SIAM J. Math.
  Anal., 28 (1997), pp.~1064--1085.

\bibitem{FT85}
{\sc G.~Falkovitch and S.~Turitsyn}, {\em Stability of magnetoelastic solitons
  and self-focusing of sound in antiferromagnet}, Sov. Phys. JETP, 62 (1985),
  pp.~146--152.

\bibitem{FS}
{\sc A.~Fokas and L.~Sung}, {\em The inverse spectral method for the {KP I}
  equation without the zero mass constraint}, Math. Proc. Camb. Phil. Soc., 125
  (1999), pp.~113--138.

\bibitem{Gra08}
{\sc P.~Gravejat}, {\em Asymptotics of the solitary waves for the generalised
  {K}adomtsev-{P}etviashvili equations}, Discrete Contin. Dyn. Syst., 21
  (2008), pp.~835--882.

\bibitem{HO}
{\sc M.~Hochbruck and A.~Ostermann}, {\em Exponential {R}unge-{K}utta methods
  for semilinear parabolic problems}, SIAM J. Numer. Anal., 43 (2005),
  pp.~1069--1090.

\bibitem{KP}
{\sc B.~B. Kadomtsev and V.~I. Petviashvili}, {\em On the stability of solitary
  waves in weakly dispersing media}, Sov. Phys. Dokl., 15 (1970), pp.~539--541.

\bibitem{KassamTrefethen2005}
{\sc A.-K. Kassam and L.~N. Trefethen}, {\em Fourth order time-stepping for
  stiff pdes}, SIAM J. Sci. Comput., 26 (2005), pp.~1214--1233.

\bibitem{Klein2008}
{\sc C.~Klein}, {\em Fourth order time-stepping for low dispersion
  {K}orteweg-de {V}ries and nonlinear {S}chr{\"o}dinger equations}, ETNA, 29
  (2008), pp.~116--135.

\bibitem{KP13}
{\sc C.~Klein and R.~Peter}, {\em Numerical study of blow-up in solutions to
  generalized {K}orteweg-de {V}ries equations},  (2013).
\newblock Preprint available at:
  \href{http://arxiv.org/abs/1307.0603}{\texttt{arXiv:1307.0603}}.

\bibitem{KleinRoidot2011}
{\sc C.~Klein and K.~Roidot}, {\em Fourth order time-stepping for
  {K}adomtsev-{P}etviashvili and {D}avey-{S}tewardson equations}, SIAM J. Sci.
  Comput., 33 (2011), pp.~3333--3356.

\bibitem{KR11}
{\sc C.~Klein and K.~Roidot}, {\em Fourth order time-stepping for
  {K}adomtsev-{P}etviashvili and {D}avey-{S}tewartson equations}, SIAM Journal
  on Scientific Computing, 33 (2011), p.~DOI: 10.1137/100816663.

\bibitem{KS12}
{\sc C.~Klein and J.~Saut}, {\em Numerical study of blow up and stability of
  solutions of generalized {K}adomtsev-{P}etviashvili equations}, J. Nonl.
  Sci., 22 (2012), pp.~763--811.

\bibitem{KSM}
{\sc C.~Klein, C.~Sparber, and P.~Markowich}, {\em Numerical study of
  oscillatory regimes in the {K}adomtsev-{P}etviashvili equation}, J. Nonl.
  Sci., 17 (2007), pp.~429--470.

\bibitem{fminsearch}
{\sc J.~C. Lagarias, J.~A. Reeds, M.~H. Wright, and P.~E. Wright}, {\em
  Convergence properties of the {N}elder-{M}ead simplex method in low
  dimensions}, SIAM Journal of Optimization, 9 (1998), pp.~112--147.

\bibitem{Liu2001}
{\sc Y.~Liu}, {\em Blow-up and instability of solitary-wave solutions to a
  generalized {K}adomtsev-{P}etviashvili equation}, Tamsui Oxf. J. Manag. Sci.,
  353 (2001), pp.~191--208.

\bibitem{MMR2012_I}
{\sc Y.~Martel, F.~Merle, and P.~Rapha{\"e}l}, {\em Blow up for the critical
  {gKdV} equation {I}: {D}ynamics near the solition},  (2012).
\newblock Preprint available at:
  \href{http://arxiv.org/abs/1204.4625/}{\texttt{arXiv:1204.4625}}.

\bibitem{MST}
{\sc L.~Molinet, J.~C. Saut, and N.~Tzvetkov}, {\em Remarks on the mass
  constraint for {KP} type equations}, SIAM J. Math. Anal., 39 (2007),
  pp.~627--641.

\bibitem{Schme}
{\sc T.~Schmelzer}, {\em The fast evaluation of matrix functions for
  exponential integrators}, PhD thesis, Oxford University, 2007.

\bibitem{MZBM}
{\sc L.~B. S.V.~Manakov, V.E.~Zakharov and V.~Matveev}, {\em Two-dimensional
  solitons of the {K}adomtsev-{P}etviashvili equation and their interaction},
  Phys. Lett. A, 63 (1977), pp.~205--206.

\end{thebibliography}
\end{document}